\documentclass[a4paper]{amsart}
\usepackage{amsthm,amsmath,amssymb}
\usepackage{bm}
\usepackage[square,numbers,sectionbib]{natbib}
\usepackage{color}
\usepackage[dvipsnames]{xcolor} 

\usepackage{times}
\usepackage[shortlabels]{enumitem}
\usepackage{enumitem}
\newcommand{\subscript}[2]{$#1 _ #2$}

\newtheorem{teor}{Theorem}[section]

\newtheorem{lemm}[teor]{Lemma}
\newtheorem{osse}[teor]{Remark}
\newtheorem{prop}[teor]{Proposition}

\newtheorem{coro}[teor]{Corollary}

\def\H{\boldsymbol{H}}
\def\W{\boldsymbol{W}}
\def\A{\boldsymbol{A}}
\def\L{\boldsymbol{L}}
\def\ff{\textbf{\textit{f}}}
\def\uu{\textbf{\textit{u}}}

\def\ddt{\frac{\d}{\d t}}

\def\n{\boldsymbol{n}}
\def\E{\boldsymbol \varepsilon}

\newcommand{\numberset}{\mathbb}
\newcommand{\N}{\numberset{N}}
\newcommand{\R}{\numberset{R}}

\def\divg{\mathrm{div}_{\Gamma}}

\def\P{\boldsymbol{P}} 
\def\v{\boldsymbol{v}}

\def\u{\boldsymbol{u}}
\def\vphi{\varphi}

\def\A{\boldsymbol{A}}

\def\phit_{\phi_t}

\def\tT{{T}}
\def\H{\boldsymbol{H}}

\def\XT{X_{{T}}}
\def\YT{Y_{{T}}}
\def\norm#1{\left\Vert#1\right\Vert }
\def\norma#1{\left\vert#1\right\vert }

\def\f{\boldsymbol f}

\def\ddt{\frac{d}{dt}}

\def\LLQ(#1,#2){L^{#1}(0,\widetilde{T};\L^{#2}(\Gamma_0))}

\def\gam{\Gamma}
\def\gt{(\gam)}

\def\z{\boldsymbol z}

\def\A{\boldsymbol{A}}

\def\Pk{\P_{\mathcal{K}}}
\def\Lts{{\L^2_\sigma(\Gamma)}}
\def\I{\boldsymbol I}
\def\KK{\mathcal K}
\def\BB{\mathcal B}
\def\AA{\mathcal A}
\def\BBB{\mathbb B}
\def\tA{\widetilde{\AA}}
\def\JJ{\mathcal J}

\newcommand{\T}{{\rm T}}
\newcommand{\Rz}{{\mathbb R}}
\newcommand{\zero}{\boldsymbol 0}

\numberwithin{equation}{section}
\begin{document}

	\title[Long-Time behavior of the tangential surface
        Navier--Stokes equation]{Long-Time behavior of the \\ tangential surface Navier--Stokes equation}

\author[A.~Poiatti]{Andrea Poiatti}
\address[A.~Poiatti]{Faculty of Mathematics, University of Vienna, Oskar-Morgenstern-Platz
	1, A-1090 Vienna, Austria}
\email[A.~Poiatti]{andrea.poiatti@univie.ac.at}

\author[U.~Stefanelli]{Ulisse Stefanelli}
\address[U.~Stefanelli]{ Faculty of Mathematics, University of
  Vienna, Oskar-Morgenstern-Platz 1, A-1090 Vienna, Austria, Vienna
  Research Platform on Accelerating Photoreaction Discovery,
  University of Vienna, W\"ahringerstrasse 17, A-1090 Wien, Austria,
  \& Istituto di Matematica Applicata e Tecnologie Informatiche {\it E. Magenes}, via Ferrata 1, I-27100 Pavia, Italy.}
\email[U.~Stefanelli]{ulisse.stefanelli@univie.ac.at}

\begin{abstract}
   We investigate the initial-value problem for the incompressible
  tangential Navier--Stokes equation {with variable viscosity} on a given
  {two-dimensional} surface without boundary. Existence of global weak and strong solutions under
  inhomogeneous forcing is proved by a fixed-point and continuation
  argument. Continuous dependence on data, backward uniqueness, and
  instantaneous regularization are also
  discussed. Depending on the effect of the inhomogeneous forcing on
  the dissipative and the nondissipative components of the system, we
  investigate 
  the long-time behavior of solutions. We prove the existence and properties of the $\sigma$-{global}
  attractor, in the case of bounded trajectories, and of the so-called
  unbounded attractor, for unbounded trajectories.  
\end{abstract}


\maketitle

%
%

\section{Introduction}
\label{sec:intro} 
This paper  is concerned with the initial-value problem for the  {\it incompressible tangential surface Navier--Stokes equation} \cite{V1} {with variable viscosity}, namely,  
\begin{align}
	&\label{t1}\partial_t\u+(\u\cdot
   \nabla_\Gamma)\u-2\P_\Gamma\divg(\nu(\cdot)\E_\Gamma(\u))+\nabla_\Gamma p=\f(\cdot,\u),\\&
	\divg\u=0,\label{t12}\\&
	\u(0)=\u_0.
	\label{tangential}
\end{align}
This equation describes the motion of a viscous incompressible fluid
on a given, sufficiently smooth,  and connected 
two-dimensional surface $\Gamma\subset \Rz^3$ without boundary. The state of the fluid
is determined by its {\it velocity} $\u : \Gamma
\times \Rz_+ \to \T \Gamma $ and its {\it pressure}
$p:\Gamma\times \Rz_+ \to \Rz$.  Here, $\T \Gamma$ denotes the tangent bundle
of $\Gamma$. 

Fluid flows on closed surfaces  arise  in
connection with  models at different scales, from 
biological membranes, to the dynamics and coating of microfluidic droplets, to soap bubbles.
Atmospheric flows around planets, as well as oceanic flows on
completely fluid-covered bodies, are also examples of surface flows. At an even larger
scale, magnetohydrodynamics is the basis of different astrophysical
models for surface flows of charged
fluids, like plasma, to model star dynamics. Surface flows have in fact
already attracted mathematical attention and the reader is referred to
\cite{Arnaudon,Chan,V1} and \cite{AGP2024,AGP2024b,Koba,Koba2} for some discussion of the case of 
stationary and evolving surfaces, respectively. See also
\cite{Brandner} for an overview and a comparison on the different derivations of the
surface Navier--Stokes equation.

In
equation \eqref{t1}, the symbol $\partial_t$ denotes the partial
derivative with respect to time and $\P_\Gamma:=\I - \n \otimes \n$ is
the projection on $\T\Gamma$, where $\I$ is the identity and $\n$
is the outward pointing normal to $\Gamma$. For given
differentiable fields $\u$ and $p$
defined on a neighborhood of $\Gamma$, we indicate by $\nabla_\Gamma
\u = {\P_\Gamma}\nabla \u \P_\Gamma= \nabla \u- (\n\otimes \n)\nabla \u - \nabla \u (\n\otimes \n) + \n\otimes \n$   and
$\nabla_\Gamma p  = \nabla p \P_\Gamma = \nabla p - (\nabla p \cdot \n)\n$ the {covariant derivative  and the scalar surface gradient, respectively}, by ${\rm div}_\Gamma
\u = \P_\Gamma {:} \nabla \u - {\rm div} \u - \n \cdot (\nabla \u) \n$
the surface divergence, and by $\E_\Gamma(\u):=
(\nabla_\Gamma\u+\nabla_\Gamma^T\u)/2 $ the surface rate-of-strain
tensor (the symbols $\cdot$, $:$, and $\otimes$ denote the scalar,
the contraction, and the tensor product, respectively). We refer to
Chapter 2 in \cite{V1}, \cite{greenbook}, or the
Appendix of \cite{Simonett} for  additional information and
details  on these
objects. The   function
$\nu$ in \eqref{t1} is the strictly positive {\it interface shear viscosity} which we assume to be space
dependent to allow for the modeling of  binary or multi-phase fluids (see, e.g., \cite{AGG,AGP}), {or even of nonisothermal fluids (e.g., \cite{Bernardi,GPPV})}.
Eventually, $\f:\Gamma\times \T\Gamma\to \T\Gamma$ is a
suitable divergence-free tangential forcing term, which  we also
allow to  
depend on $\u$. 

The aim of this paper is to investigate the well-posedness and the
long-time behavior of solutions to the initial-value problem
\eqref{t1}--\eqref{tangential}, especially in the nonhomogeneous case $\f
\not =\zero$. Note that in most settings,  
nontrivial forcings $\f$ naturally occur. Referring to the above list
of examples, we remark that microdroplets are
subject to friction and surface tension, biological membranes are driven by cell
electrochemistry, and gravity influences soap bubbles. At another
scale, atmospheric and oceanica flows are driven by thermal,
chemical, geophysical, and tidal processes, while astrophysical flows are subject to
electromagnetic, thermal, and gravitational effects. 

Existence results for the Navier--Stokes equation in three dimensions
are a mainstay of nonlinear PDE theory. The reader is referred to the
monographs \cite{Constantin88,Galdi94,Girault86,Temam84} for a collection of material. Existence
for the surface
Navier--Stokes case has also been considered.
The surface {stationary} Stokes equation has been considered in
\cite{V1} while the surface {evolutionary} Navier--Stokes equation, {with $\ff=\mathbf 0$ and constant viscosity $\nu$}, has been
studied in \cite{Simonett} and checked to be locally well-posed in
time. {Additionally, solutions to the same equation on two dimensional surfaces, departing from a sufficiently regular divergence-free initial datum, has been shown to exist globally in \cite{Simonett2}.}

Our first result concerns the global well-posedness of
problem \eqref{t1}--\eqref{tangential} in the general inhomogeneous
case $\f \not = \zero$ {and variable viscosity}. We prove that
weak solutions to problem \eqref{t1}--\eqref{tangential} globally
exist, regardless of the size of the initial datum. In case the
initial datum is in $\H^1$ weak solutions are actually strong.
Moreover, we show continuous
dependence on data, entailing uniqueness, as well as instantaneous
regularization, see Theorems \ref{thm1} and \ref{insta}. In the setting of the
tangential equation, this
extends the results in \cite{Simonett,Simonett2} to the general, global-in-time, inhomogeneous case.

In addition, we show the so-called {\it backward uniqueness} property,
namely, 
that two
trajectories coinciding at some point in time necessarily coincide for
all times, see Theorem \ref{backuni}. This property, which is well known in the case of the
Navier--Stokes problem in a three-dimensional domain
\cite{Escauriaza,Iskauriaza,Seregin}, was apparently not yet
investigated for surface flows.

The main focus of this paper is however on the long-time behavior of
solutions. The asymptotic behavior of a solution to the Navier--Stokes
equation for large times in three-dimensional domains is classical,
see the monographs \cite{WeakGlobalNSAttractor,Foias2,Temam} and the
references therein. The long-time behavior for the surface Navier--Stokes equation was initiated in \cite{Simonett,Simonett2}, {concerning the equilibria of a single trajectory}. In this setting a key role is played by the so-called {\it Killing fields} on $\Gamma$,
namely, fields $\u_K \in C^\infty(\Gamma;\T\Gamma)$ with the property
that $\E_\Gamma(\u_K)=\zero$. Note that the occurrence of 
Killing fields is completely determined by the 
geometry of the given $\Gamma$. In particular, Killing fields
correspond to rotational symmetries of $\Gamma$. Hence, the space of Killing field is
$n$-dimensional, with $n=0,1,2$, or $3$.

In case $\f=\zero$, Killing fields are
stationary solutions to \eqref{t1}--\eqref{t12}. In particular, the
tangential surface Navier--Stokes equation is nondissipative on
Killing fields, which in fact correspond to
equilibria. 
In \cite{Simonett}, killing fields are identified as equilibria and are
proved to be stable. Eventually, trajectories starting close enough to
equilibria are shown to exist globally and to converge to equilibria
exponentially fast in time. {Also, in \cite{Simonett2} it is shown that any weak solution, departing from a sufficiently smooth divergence-free initial datum, converges exponentially fast to an equilibrium, that is, to a Killing field.}

{Inspired by these results,}  we turn our attention on sets of trajectories instead {of single ones}.  Moving  from our global well-posedness results, we  discuss  the long-time behavior in terms of the dynamical system
generated by equation \eqref{t1}--\eqref{t12}. One decomposes
$\u = \u_K + \u_{NK}$ in its Killing and non-Killing components by
projecting on the finite-dimensional space of Killing fields. As
{ only }the component $\u_{NK}$ dissipates {thanks to the tensor $\boldsymbol\varepsilon_\Gamma$}, one is asked to investigate the
properties of $\f$ with respect to the dynamics of the Killing
component $\u_{K}$.

In Theorems \ref{thm:323} and \ref{thm:325}, we consider the case of
bounded trajectories. Different assumptions entailing such boundedness are
presented, including, for instance,
$$\int_\Gamma\f_{NK}(\cdot,\u)\cdot \u \leq \delta \|\u_{NK}\|^2_{\L^2} + C \|\u_K\|^2_{\L^2},
\quad \int_\Gamma\f_{K}(\cdot,\u)\cdot \u \leq 0 \quad \forall \u \in \L^2_\sigma(\Gamma)$$
for some $C>0$ and $\delta >0 $ small, where $\f = \f_K + \f_{NK}$ is
the decomposition of $\f$ in its Killing and non-Killing components
and we have used the space $\L^2_\sigma(\Gamma) :=\{\v
\in \L^2(\Gamma)\::\: {\rm div}\, \v=0 \ \text{a.e. in} \ \Gamma\}$.

In the case of bounded trajectories, we prove that the generated
semigroup restricted to $\{\u\in \L^2_\sigma(\Gamma)\::\:
\|\u_{K}\|_{\L^2}\leq r\}$ for $r>0$ admits
a global attractor $\mathcal A_r$ (nonempty, compact, connected,
invariant, and attracting) of finite fractal dimension and that
$\cup_{r>0}\mathcal A_r{=\cup_{m\in \mathbb N}\mathcal A_m}$ is minimal among closed sets attracting all
bounded sets {of the phase space $\L^2_\sigma(\Gamma)$}. {Since  this is a countable union, such} {set is called {\it $\sigma$-global attractor.}} Moreover, we prove that for all $r>0$ there exists an
exponential attractor $\mathcal M_r$ (nonempty, compact, positively
invariant, and exponentially attracting) of finite fractal dimension,
so that {$\cup_{m\in \mathbb N}\mathcal M _m$, which we}  term  {\it $\sigma$-exponential attractor}, exponentially attracts all bounded
sets {of $\L^2_\sigma(\Gamma)$}. 

The case of possibly unbounded trajectories is discussed in Theorems
\ref{thm:332} and \ref{t1b}, instead. If $$\exists R>0: \quad \int_\Gamma\f_K(\cdot,\u)\cdot \u \geq 0 \quad
\forall \u \in \L^2_\sigma(\Gamma) \  \ \text{such that} \ \  \| \u_K\|_{\L^2} \geq
R, $$
one can consider the notion of 
{\it unbounded attractor}, \cite{Chepyzov}, which consists of the values of all
bounded-in-the-past complete trajectories of the dynamical system. This {set} has many of the
properties of a global attractor. In particular, it is nonempty,
closed, and invariant. It is moreover bounded compact: its intersection
with any closed bounded set is compact. On the other hand, its attractivity can proved just for bounded
trajectories departing from any given bounded set.

The paper is organized as follows. In Section \ref{sec:strong}, we
introduce notation and recall some basic facts on Killing fields. All
statements are collected in Section \ref{main}. In particular,
well-posedness (Theorem \ref{thm1}) and backward uniqueness (Theorem
\ref{backuni}) are in Section \ref{sec:well}. Some preliminaries on
long-time behavior, including the decomposition $\u = \u_K + \u
_{NK}$, the bounds on these components, and instantaneous regularization
(Theorem \ref{insta}), are collected in Section \ref{longtt}. The
notion of attractor for the underlying dynamical system is discussed
in Section \ref{notion}. Eventually, the statements on the existence of attractors for
the case of bounded (Theorems \ref{thm:323} and \ref{thm:325}) and unbounded (Theorems \ref{thm:332} and \ref{t1b}) trajectories are given in Sections
\ref{sec:attr} and \ref{unbd}, respectively. The proof of the
statements from Sections
\ref{sec:well}, \ref{longtt}, \ref{sec:attr}, and \ref{unbd} is then
given in Sections \ref{sec:pr1},
\ref{sec:prel}, \ref{bdda}, and \ref{unbdd}, respectively. Eventually,
Appendix \ref{sec:appendix} contains a lemma on the existence of the unbounded attractor.

 

 \section{Notation and functional setting} \label{sec:strong}
 Let us consider a smooth, compact, connected, embedded hypersurface $\Gamma\subset \R^3$ without boundary.
 In the following, the classical Sobolev spaces are denoted as usual by $W^{k,p}(\Gamma)$%
 , where $k\in \mathbb{N}$ and $1\leq p\leq \infty $, with norm $\Vert \cdot
 \Vert _{W^{k,p}(\Gamma)}$. The Hilbert space $W^{k,2}\gt$ is denoted
 by $H^{k}\gt$ with norm $\Vert \cdot \Vert _{H^{k}\gt}$. We then set
 $H^{-1}\gt=(H^1\gt)'$. Given a vector space $X$ of functions defined
 on $\Gamma$, we denote by $\boldsymbol{X}$ the generic space of
 tangential vectors or matrices, with each component in $X$. We recall
 that a tangential vector $\v$ is such that $\v\cdot\n=0$ on $\gam$,
 i.e., $\P_\Gamma\v=\v$. In a similar way, a tangential matrix
 $\boldsymbol M$ is a matrix such that $\P_\gam\boldsymbol M\P_\Gamma=\boldsymbol
 M$. An example of a tangential matrix is the covariant derivative of
 a vector field $\v$, denoted by $\nabla_\Gamma\v$.  The symbol
  $\vert \boldsymbol{v} \vert$  stands for  the Euclidean norm
 of $\boldsymbol{v}\in  \boldsymbol{X}$, i.e., $\vert\boldsymbol{v}\vert^2=\sum_{j}\|v_j\|_X^2$.
 We then denote by $(\cdot,\cdot )$ the inner product in $\boldsymbol{L}^{2}\gt$ and by $\Vert \cdot \Vert $
 the induced norm. We also denote by $(\cdot ,\cdot )_{H}$ and $\Vert \cdot
 \Vert _{H}$ the natural inner product and its corresponding induced norm in the Hilbert
 space $H$. Moreover, given a Banach space $X$, $T>0$, and $q\in[2,\infty]$, we denote by $L^q (0,T; X )$ the Bochner space of $X$-valued $q$-integrable (or essentially bounded) functions. We then denote by $BC([0,T];X)$ the Banach space of bounded continuous functions on $[0,T]$, equipped with the supremum norm. The space $BUC([0,T];X)$ is then its subspace of bounded and uniformly continuous functions. Finally, $W^{1,p} (0, T ; X)$, $1 \leq p < \infty$, is the space of functions $f$ such that $\partial_t f\in L^p(0,T;X)$ and $f\in L^p(0,T;X)$, where $\partial_t$ is the vector-valued distributional derivative of $f$. We set $H^1 (0, T ; X) =W^{1,2}(0,T;X)$.
 
   In order to handle  divergence-free vector fields, we introduce the spaces 
 \begin{align*}
 	&\L^2_\sigma\gt:=\{\v\in \L^2\gt: \divg \v=0\text{ a.e. on } \Gamma \},
 	\\
 	&\H^1_\sigma\gt:=\{\v\in \H^1\gt: \divg \v=0\text{ a.e. on } \Gamma \},
 \end{align*}
 and observe that $\H^1_\sigma\gt,\L^2_\sigma\gt ,\H^1_\sigma\gt')$
  forms a Hilbert triplet.  We also introduce the surface 
 Helmholtz  orthogonal projector on $\Gamma$, denoted by
  \begin{align}
  	\P_0:\L^2(\Gamma)\to \L^2_\sigma(\gam),
 	\label{Leray}
 \end{align}
 see, e.g., \cite{Simonett2} for  details.  
 
 We now introduce the notion of Killing fields as 
 \begin{align*}
 	\mathcal K:=\{\u\in \L^2_\sigma(\Gamma)\cap \boldsymbol
   C^\infty( \Gamma;\T\Gamma ):\ \E_\Gamma(\u)\equiv \boldsymbol 0\},
 \end{align*}
 where we recall that the divergence free condition is already
 entailed by the condition $\E_\Gamma(\u)\equiv \boldsymbol 0$, by $\divg
 \u=\text{tr}(\E_\Gamma(\u))=\boldsymbol 0$. Moreover, as noticed in
 \cite{Simonett,Simonett2}, one can show that any vector field $\u\in
 \W^{1,q}(\Gamma)$ such that $\E_\Gamma(\u)\equiv \boldsymbol 0$ is smooth, see
 for instance \cite[Lemma 3]{Priebe}. In particular, the Killing
 fields on a Riemannian manifold form a  Lie sub-algebra  of
 the Lie algebra of all tangential fields. We also recall that the Killing fields of a Riemannian manifold
 $(M, g)$ are the infinitesimal generators of the isometries $I(M, g)$
 on $(M, g)$, i.e., the generators of flows that are isometries on
 $(M, g)$. Furthermore,  if  $(M, g)$ is
 complete,  which is  our case, the Lie algebra of Killings
 fields is isometric to the Lie algebra of $I(M, g)$, (see for example
 \cite[Corollary III.6.3]{Sakai}). From this observation we can deduce
 that $\KK$ is a finite dimensional vector space, which is thus closed
 in $\L^2_\sigma(\Gamma)$. Indeed, in the case of a surface embedded
 in $\R^{d+1}$, from \cite[Proposition III.6.5]{Sakai} one deduces
 that $\text{dim}\ \KK \leq d(d+1)/2$, with  the  equal sign characterizing
  the case where $\Gamma$ is isometric to  an  Euclidean sphere. In our case $d=2$, and thus $\text{dim}\ \KK \leq 3$. Additionally, if $(M, g)$ is compact and the Ricci tensor is negative-definite everywhere, then any Killing field on $M$ is equal to zero and $I (M, g)$ is a finite group,
 see, for example \cite[Proposition III.6.6]{Sakai}. Namely, if $(M, g)$ is a two-dimensional Riemannian
 manifold with negative Gaussian curvature then any Killing field is $\boldsymbol 0$. Note that, among two-dimensional compact closed surfaces only those of genus 0 and 1 may have nonzero Killing
 fields  \cite[Thm.~6]{Myers}.  Some possible examples for which the Killing fields are nontrivial are then the following (see, e.g., \cite{Simonett, Olshsup}):
 \begin{itemize}
 	\item $ \Gamma= \mathbb S^2$. Then, $\text{dim } \KK = 3$ and any $\u\in \KK$ is a
 	rotation about an axis spanned by a vector $\omega = (\omega_1,\omega_2,\omega_3) 	\in \R^3$, so that
 	$\u \in \KK$ is given by
 	$$\u(x) = \omega \times x,\quad  x \in \mathbb S^2,$$
 	for some $\omega \in\R^3$.
 	\item $\Gamma=\mathbb T^2$ with  parametrization 
    \begin{align*}
 &	x_1 = (R + r \cos \phi) \cos \theta,\\&
 	x_2 = (R + r \cos \phi) \sin \theta,\\&
 	x_3 = r \sin \phi,
 	\end{align*}
 	where $\phi, \theta \in [0, 2\pi)$ and $ 0 < r < R$. Then, the velocity
 	field $\u = \omega\boldsymbol e_3 \times x$, with $	\omega \in \R$, is an element of $\KK$, i.e., the fluid on the torus
 	rotates about the $x_3$-axis with angular velocity $\omega\boldsymbol e_3$, and these rotations can be proven to be the only isometries on $\mathbb T^2$. This means that $\text{dim}\ \KK=1$. 
 	\item As observed in \cite{Olshsup}, any surface $\Gamma$ of revolution supports a nonzero Killing
 	field. Moreover, it looks plausible (although apparently 
        still not discussed in the  literature) that among closed
 	compact smooth surfaces isometrically embedded in $\R^3$ only
        surfaces of revolution support  nontrivial  Killing fields. In the case of a surface of revolution,  one has KIlling fields  $\u\in \KK$ is of the form 
 	$$
 	\u=\omega\boldsymbol e_z \times x,
 	$$
 	where $\omega\in \R$ and $\boldsymbol e_z$ is the unit vector in
        the direction of the axis of rotation.  This  corresponds to a rotation about the $z$-axis with angular velocity $\omega\boldsymbol e_z$. 
 \end{itemize}
 
 The importance of the Killing fields $\KK$ lies in the lack of
 dissipativity properties of the tangential Navier--Stokes system on
 this vector space. In order  to take  this issue  into
 account,  solutions need to be decomposed in a Killing component
 in $\KK$ and a non-Killing one. More precisely, we  
 define the orthogonal $\L^2$ projector  on the (closed) Killing vector space $\mathcal K$:  
 $$
 \Pk: \Lts\to \mathcal{K}.
 $$
 Since the space $\mathcal K$ is finite dimensional and at most of
 dimension $n=3$, we can always find an  orthonormal  $\L^2$-basis
 $\{\v_n\}$, where $n$ can be $0$,   $1$, $2$, or $3$.  Therefore, any element $\v\in \L^2_\sigma(\Gamma)$ can be uniquely decomposed into 
 $$
 \v=\Pk \v+(\boldsymbol I-\Pk)\v:=\v_{K}+\v_{NK},
 $$
 where $\boldsymbol I$ is the identity operator on $\L^2_\sigma(\Gamma)$.
 
 Note that, as shown in \cite[Remark 4.10]{Simonett2}, it also holds
 $$
 \H^1_\sigma(\Gamma)=\KK\oplus {(}(\I-\Pk)\L^2_\sigma(\Gamma)\cap \H^1(\Gamma){)}.
 $$
In conclusion, we recall the following Korn's inequality, which is fundamental to deal with the lack of dissipativity properties of the system (see, e.g. \cite[(4.7)]{V1}): 
\begin{align}
	\norm{\v}_{\H^1(\Gamma)}\leq C(\norm{\v}+\norm{\E_\Gamma(\v)}),\quad \forall \v\in \H^1(\Gamma),
	\label{Korn1}
\end{align}
which entails that 
\begin{align}
	\norm{\v}_{\H^1(\Gamma)}\leq C\norm{\v},\quad \forall \v\in \KK.
	\label{Korn2}
	\end{align}
	Furthermore, we also have  (see, for instance, \cite[Theorem A.3]{Simonett2})
	\begin{align}
		\norm{\v}_{\H^1(\Gamma)}\leq C_P\norm{\E_\Gamma(\v)},\quad \forall \v\in (\I-\Pk)\L^2_\sigma(\Gamma)\cap \H^1(\Gamma),
		\label{Korn3}
	\end{align}
	for some $C_P>0$.
        \section{Main results}\label{main}
        
\subsection{Well-posedness and backward uniqueness}
\label{sec:well}
We state here a well-posedness result concerning the tangential
surface Navier--Stokes equation.  This is an extension of the
former \cite{Simonett2}, where    $\f=\boldsymbol
0$ and  a  constant viscosity  $\nu$ are assumed.  
The proof of this result is postponed to Section \ref{sec:pr1}. We have:
\begin{teor}[Well-posedness]
	\label{thm1}
Let $\nu\in W^{1,\infty}(\Gamma)$ be such that $\nu\geq \nu_*>0$, and
$\f:\Gamma\times \T\Gamma\to \T\Gamma$ be a  Carath\'eodory  function such that there exist $C_1>0$ so that
\begin{align}
\norm{\f(x,\boldsymbol 0)}\leq C_1, \quad  \text{for a.e.} \ x \in
  \Gamma, 
\label{C1}
\end{align}
and there exists $C_2>0$ so that
  \begin{align}
	\norm{ \f(x, \v_1(x))-\f(x,\v_2(x))}\leq C_2\norm{
    \v_1-\v_2},\quad \forall \v_1,\v_2\in \L^2(\Gamma), \quad  \text{for a.e.} \ x \in
  \Gamma .
	\label{C2}
\end{align}
Let $\u_0\in \L^2_\sigma(\Gamma)$. Then there exists a unique global weak solution $\u:\Gamma\times[0,\infty)\to \T\Gamma$ to 
\eqref{t1}--\eqref{tangential} such that, for any $T>0$, 
\begin{align}
	&\u\in C([0,T]; \L^2_\sigma(\Gamma))\cap L^2(0,T;\H^1_\sigma(\Gamma))\cap H^1(0,T;\H_\sigma^1(\Gamma)'),\label{regga}
\end{align}
and, for almost  every  $t>0$,
\begin{align*}
	\langle \partial_t\u, \boldsymbol\eta\rangle_{\H^1_\sigma(\gam)',\H^1_\sigma(\Gamma)}+\int_\Gamma (\u\cdot\nabla_\Gamma)\u\cdot \boldsymbol \eta+\int_\Gamma 2\nu(x)\E_\Gamma(\u):\E_\Gamma(\boldsymbol{\eta})=\int_\Gamma  \f(x,\u)\cdot \boldsymbol{\eta},\quad \forall \boldsymbol{\eta}\in \H^1_\sigma(\gam).
\end{align*}

  Moreover,  given two weak solutions $\u_1,\u_2$ as above,
 departing from initial data $\u_{0,1},\u_{0,2}\in
 \L^2_\sigma(\Gamma)$,  for all $T>0 $ the following 
 continuous-dependence estimate  holds  
  \begin{align}
 	\sup_{t\in[0,T]}\norm{\u_1(t)-\u_2(t)}^2+2\nu_*\int_0^T\norm{\mathcal
    E_S(\u_1(t)-\u_2(t))}^2d t\leq  C  \norm{\u_{0,1}-\u_{0,2}}^2,
 	\label{contdep2}
 \end{align}
where the constant $C>0$ depends on $T$ and $\norm{\u_{0,i}}$, $i=1,2$.

 Additionally,   if $\u_0\in \H^1_\sigma(\Gamma)$  the unique
weak solution is also a global strong solution.  Namely, for all
$T>0$ we have  
\begin{align}
	&\u\in C([0,T]; \H^1_\sigma(\Gamma))\cap L^2(0,T;\H^2(\Gamma))\cap H^1(0,T;\L^2_\sigma(\Gamma))\label{strong}
\end{align}
 and 
\eqref{t1}--\eqref{tangential} are satisfied  almost everywhere on $\Gamma\times(0,\infty)$.

 Eventually,   the unique
global weak solution departing from  a given $\u_0\in
\L^2_\sigma(\Gamma)$  instantaneously regularizes.  In
particular, it is  
such that, for any $T>0$ and any {$\tau\in(0,T)$}, there exists $ C>0$
 depending on $T$, $\norm{\u_0}$, $\f$, $\Gamma$, and  on 
the parameters of the problem, but not on $\tau$, so that
\begin{align}
	\norm{\u}_{C([\tau,T];\H^1_\sigma
  (\Gamma))}+\norm{\u}_{L^2(\tau,T;\H^2(\Gamma))}+\norm{\partial_t\u}_{L^2(\tau,T;\L^2_\sigma(\Gamma))}\leq
  \frac{ C }{\sqrt \tau}.
	\label{regulari}
\end{align} 
\end{teor}
\begin{osse}\rm (Divergence-free condition)
	Without loss of generality  one  can directly assume that $\f(\cdot,\v)$ is divergence free for any suitably regular divergence free vector field $\v$. Indeed, only the divergence-free part of the function $\f$ has an active role in the system.
\end{osse}
\begin{osse}\rm (Pressure)
	By  the classical   De Rahm's theorem, in the case of strong solutions we can also deduce the existence of a unique pressure $p$ such that $\int_\Gamma p\equiv 0$	and $p\in L^2(0,T;H^1(\Gamma))$ for any $T>0$.
\end{osse}
 We additionally point out that the  system
\eqref{t1}--\eqref{tangential} has the property of backward uniqueness
for strong solutions, i.e., if two strong solutions coincide at some
instant of time, then they coincide for all previous times. This
property has been shown for the 2D Navier--Stokes equations on bounded
domains (see, e.g., \cite{Tartar}), but, as far as we know, it has not
been shown in the case of closed surfaces. Indeed, in the latter case
one cannot exploit the full dissipativity nature of the equations as
in the bounded domains case. This is to some extent similar to what
has been adressed in \cite{Foias} in the case of semi-dissipative
Boussinesq equations, where a part of the system is not
dissipative. In our problem, the main ingredients will be a careful
use of Korn's inequality to deal with the lack of dissipativity,
together with the decomposition of the solution in its Killing and
non-Killing components.  An additional  difficulty comes from
the presence of a variable viscosity $\nu$. 
We thus have the following,  which is also proved in Section
\ref{sec:pr1}. 
\begin{teor}[Backward uniqueness]\label{backuni}
	Let assumptions \eqref{C1}--\eqref{C2} hold, together with $\nu\in W^{1,\infty}(\Gamma)$ and $\nu\geq \nu_*>0$. For $i=1, 2$, let $\u_i$ be two {global} strong solutions to \eqref{t1}--\eqref{tangential} according to Theorem \emph{\ref{thm1}}. If there exists a time $T^*>0$ such that $\u_1(T^*) =\u_2(T^*)$ almost  everywhere  on $\Gamma$, then for any $t\in[0,T^*]$, we have $\u_1(t) =\u_2(t)$ almost everywhere on $\Gamma$.
\end{teor}
\begin{osse}\rm (Solutions do not cross)\label{coincidenza}
	Notice that, as a consequence of  Theorem  \emph{\ref{backuni}}
         and the uniqueness result of Theorem \emph{\ref{thm1}}, if two (global) strong solutions coincide at some instant of time $T^*\geq0$, then they must coincide for any $t\geq0$.  
\end{osse}

\subsection{Long-Time behavior: Preliminaries}
\label{longtt}

Observe that the main peculiarity of the surface tangential Navier--Stokes equation is related to the fact that the system is not dissipative in general. Indeed, in case the Killing field space $\mathcal{K}\not=\{\boldsymbol 0\}$, the dissipation due to the stress tensor $\E_\Gamma$ does not act on the component of the velocity in the direction of a Killing field (since $\E_\Gamma(\v)=\boldsymbol 0$ for any $\v\in \mathcal K$). 

In order to study the long-time behavior of the solutions to
\eqref{t1}--\eqref{tangential},  we consider the dynamical system
in the space $\Lts$.  For the sake of readability,  all proofs
of the statements  of this section are postponed to Section \ref{sec:prel}.

Under the same assumptions of Theorem \ref{thm1}, we can now define a dynamical system $(\Lts,S(t))$ where
\begin{equation*}
	S(t):\Lts\rightarrow \Lts, \quad
	S(t)\boldsymbol{u}_{0}=\boldsymbol{u}(t),\quad \forall \, t\geq 0.
\end{equation*}%
Observe that $S(t)$ is a continuous semigroup, since it satisfies the following properties:

\begin{itemize}
	\item $S(0)=\boldsymbol I$;
	
	\item $S(t+\tau)=S(t)S(\tau)$, for every  $t,\, \tau \geq
          0$;   
	
	\item $t\mapsto S(t)\boldsymbol{u}_0\in C([0,\infty);\Lts)$,
	for every $\boldsymbol{u}_0\in\Lts$;
	
	\item $\boldsymbol{u}_0\mapsto S(t)\boldsymbol{u}_0\in C(\Lts;%
	\Lts)$, for any $t\in[0,\infty)$.
\end{itemize}
In particular, the last property is a direct consequence of the continuous dependence estimate \eqref{contdep2}.

Observe also that the backward uniqueness property, together with the
well-posedness result (see Theorems \ref{thm1}  and 
\ref{backuni})  implies  that if any two (sufficiently smooth) trajectories intersect, then they must be identical.

As already anticipated, due to the lack of dissipation of the entire
velocity, it is useful to decompose the semigroup into the sum of two
other operators. First, we introduce $\f_K:\Gamma\times
\L^2_\sigma(\Gamma)\to \L^2_\sigma(\Gamma)$  via 
\begin{align*}
	\f_K(\cdot,\u):=\P_{\mathcal{K}}\f(\cdot,\u)=\sum_{j=1}^n\left(\int_\Gamma \f(x,\u(x))\cdot \v_j(x)\right)\v_j,
\end{align*}
where $\{\v_j\}$ is  a given  orthonormal basis of $\mathcal K$.
We also define $\f_{NK}:=\f-\f_K$.

\noindent Now, given $\u_0\in \Lts$, setting $\overline{\u}:=S(t)\u_0$, we define $\u_{NK}(t)$ as the solution to 
	\begin{align}
		&\label{t11}\partial_t\u_{NK}+\overline{\u}\cdot \nabla_\Gamma\overline{\u}-2\P\divg(\nu\boldsymbol\varepsilon_\Gamma(\u_{NK}))+\nabla_\Gamma p=\f_{NK}(\cdot,\overline{\u}),\\&
		\divg\u_{NK}=0,\\&
		\u_{NK}(0)=(\boldsymbol I-\P_{\mathcal{K}})\u_0,
		\label{tangential1}
\end{align} 
 and $\u_K(t)$ to be  the solution to 
\begin{align}
&\label{K1}\partial_t\u_K=\f_K(\cdot,\overline{\u}),\\&
	\label{K2}\u_K(0)=\P_{\mathcal{K}}\u_0.
\end{align}
We then have the following.

\begin{prop}[Decomposition of $\u$]
	Let $\u_0\in \Lts$ and set $\overline{\u}=S(t)\u_0$.  Under
         assumptions \eqref{C1}--\eqref{C2} it holds  that 
	\begin{align*}
		\u_{NK}=(\boldsymbol I-\P_{\mathcal{K}})\overline{\u},\quad \u_{K}=\P_{\mathcal{K}}\overline{\u},
	\end{align*}
	where $\u_{NK}$ and $\u_{K}$ are defined in \eqref{t11}--\eqref{tangential1} and \eqref{K1}--\eqref{K2} respectively.
	\label{dec}
\end{prop}

	 In order to analyze  the asymptotic behavior of the
        Killing projection $\u_K$, we need to further distinguish 
        some cases,  depending on the dissipative nature of the
        forcing term $\f_K$.  We collect some observations in the
        following proposition, which is proved in Section
        \ref{sec:prel}.  
	
	\begin{prop}[Bounds on the Killing component $\u_K$]
          \label{lemmadec}
		Under assumptions \eqref{C1}--\eqref{C2}, let us
                assume, additionally, that there exists $ C_3 >0$ such that 
		\begin{align}
			\int_\Gamma \f_K(x,\u)\cdot \u\leq  C_3 \norm{\P_{\mathcal{K}}\u}^2+ C_3 ,\quad \forall \u\in \Lts. 
			\label{uk1}
		\end{align}
		Then there exists $ C_4 >0$ such that 
		\begin{align}
			\norm{\u_K(t)}^2\leq e^{2  C_4  t}(\norm{\u_K(0)}^2+ C_4 ),\quad \forall t\geq 0,
			\label{exp}
		\end{align}
		for any $\u_0\in \Lts$.
		Moreover, we distinguish three more specific cases:
		\begin{itemize}
			\item[(i)] If $\f_K$ is independent of $\u$, i.e.,$\f_K:\Gamma\to \T\Gamma$, then
			\begin{align}
				\int_{\Gamma}\f_K(x)\cdot \u_{K}(x,t)=	\int_{\Gamma}\f_K(x)\cdot \u_{K}(x,0)+t\norm{\f_K}^2,\quad \forall t\geq 0,
				\label{fd}
			\end{align}
			whereas 
			\begin{align}
				\norm{\u_K(t)}^2-	\left(\int_{\Gamma}\f_K(x)\cdot \u_{K}(x,t)\right)^2=	\norm{\u_K(0)}^2-	\left(\int_{\Gamma}\f_K(x)\cdot \u_{K}(x,0)\right)^2,\quad \forall t\geq 0.
				\label{uk}
			\end{align}
			As a consequence, 
			\begin{align}
				\norm{\u_K(t)}^2=\norm{\u_K(0)}^2+t^2\norm{\f_K}^4+2t\norm{\f_K}^2\int_{\Gamma}\f_K(x)\cdot \u_{K}(x,0),\quad \forall t\geq 0,\label{uk2}
			\end{align}
			entailing $\norm{\u_K(t)}\to \infty$ as $t\to
                        \infty$  in case $\f_K \not =  0$. 
			\item[(ii)] If 
			\begin{align}
				\int_\Gamma \f_K(x, \u)\cdot \u\leq 0,
				\label{nega}			
				\end{align}
			for any $\u\in \L^2_\sigma(\Gamma)$, then 
			$$
			\norm{\u_K(t)}\leq \norm{\u_K(0)},\quad \forall t\geq 0.
			$$ 
			\item[(iii)] If
			\begin{align}
			\int_\Gamma \f_K(x, \u)\cdot \u\geq 0,
			\label{pos}\end{align}
			for any $\u\in \L^2_\sigma(\Gamma)$, then
			$$
			\norm{\u_K(t)}\geq \norm{\u_K(0)},\quad \forall t\geq 0.
			$$ 
		\end{itemize}
	\end{prop}
	\begin{osse}\rm ($\f_{K}$ antisymmetric)
		Notice that a necessary condition for \eqref{nega} or \eqref{pos} to hold (and $\f_K\not = 0$) is that $\f_{K}$ is antisymmetric with respect to $\u$, i.e., $\f_K(x,\u)=-\f_K(x,-\u)$.
	\end{osse}
\begin{osse}\rm (Weaker assumptions)
	We observe that, in case (iii) of  Proposition
        \emph{\ref{lemmadec}}  one can also assume the  weaker  assumption $\int_\Gamma \f_K(x, \u)\cdot \u\geq 0$,  for any $\u\in \L^2_\sigma(\Gamma)$ such that $\norm{\Pk\u}\geq R_0$, for some $R_0\geq0$. In this case, it is immediate to deduce by a similar argument that, if there exists $t_0\geq0$ such that $\norm{\u_K(t_0)}\geq R_0$, then 
	$\norm{\u_K(t)}\geq \norm{\u_K(t_0)}$ for any $t\geq t_0$. In
        particular, if $\norm{\u_K(0)}\geq R_0$, then
        $\norm{\u_K(t)}\geq \norm{\u_K(0)}\geq R_0$ for any $t\geq
        0$. 
        \label{relax}
\end{osse}

Our second general result is showing that, possibly depending on the
behavior of $\u_{K}$, the function $\u_{NK}$ has an exponential decay {(up to some additive uniform constant)} as
time goes to infinity, as long as the constant {$C_1,C_2$} appearing in the
Theorem \ref{thm1} are sufficiently small compared to the viscosity
lower bound $\nu_*$. This   reflects the fact that  
the interaction between $\u_K$ and $\u_{NK}$ is mostly  due to 
the transport term $\overline{\u}\cdot \nabla_\Gamma\overline{\u}$,
which does not appear in the energy estimate. The proof is presented
in Section \ref{sec:prel}.

\begin{prop}[Bounds on the non-Killing component $\u_{NK}$]
\label{pp}	Under assumptions \eqref{C1}--\eqref{C2} and
\eqref{uk1}, if additionally  there exist $C_5,\, C_6>0$ such that 
	\begin{align}
		\int_\Gamma{\f_{NK}(x,\u)\cdot \u}\leq  C_5 \norm{(\boldsymbol I-\P_{\mathcal{K}})\u}^2+ C_6 \norm{(\boldsymbol I-\P_{\mathcal{K}})\u},\quad \forall \u\in \Lts,
		\label{extra}
	\end{align}	
	with $\zeta:=\frac{2\nu_*}{C_P^2}-2  C_5 >0$ ($C_P$
         is  defined in \eqref{Korn3}), then there exists $\omega>0$ such that
	\begin{align}\label{gr}
		\norm{\u_{NK}(t)}^2\leq e^{-\zeta t}\norm{\u_{NK}(0)}^2+\omega,\quad \forall t\geq0,
	\end{align}
for any $\u_0\in \Lts$. 

\noindent On the other hand, under the additional assumption \eqref{nega},  assuming 
	\begin{align}
	\int_\Gamma{\f_{NK}(x,\u)\cdot \u}\leq  C_5 \norm{(\boldsymbol I-\P_{\mathcal{K}})\u}^2+ C_6 \norm{(\boldsymbol I-\P_{\mathcal{K}})\u}+ C_6 \norm{\Pk\u}^2,\quad \forall \u\in \Lts,
	\label{extra2}
\end{align}	
and again  $\zeta:=\frac{2\nu_*}{C_P^2}-2  C_5 >0$, there exists $\omega>0$ such that 
	\begin{align}\label{gr2}
	\norm{\u_{NK}(t)}^2\leq e^{-\zeta t}\norm{\u_{NK}(0)}^2+\omega(1+\norm{\u_K(0)}^2),\quad \forall t\geq0,
\end{align}
for any $\u_0\in \Lts$.
\end{prop}
\begin{osse}\rm (More general $\f_{NK}$)	Thanks to the results of Proposition \ref{lemmadec} case (ii),
        i.e., under the  additional  assumption \eqref{nega},
        we can assume more general hypotheses of the behavior of
        $\f_{NK}$ {(see \eqref{extra2})}, which is now allowed to depend also on $\u_K$ in
        the estimate from above. This means that we can assume some
         influence  of the Killing component on the non-Killing component of the trajectory.
	\end{osse}
	\begin{osse}\rm (Case $\f_{NK}= 0$)
		Note that, in the case $\f_{NK}= 0$, we deduce
		$$
		\norm{\u_{NK}(t)}^2\leq e^{-\frac{4\nu_*}{C_P^2}t}\norm{\u_{NK}(0)}^2,
		$$
		so that $\u_{NK}\to 0$ exponentially as $t\to \infty$,
                 see \eqref{energy} in  the proof in Section
                \ref{sec:prel},  having a right-hand side $0$ in
                this case.  
                {This is in agreement with the results in \cite{Simonett2}.}
	\end{osse}
	\begin{osse}\rm (Examples)
		Some possible examples of forcing terms satisfying the assumptions of Proposition \ref{lemmadec} and \eqref{pp} are, for instance,
		\begin{align*}
			&\f_1=\f_1(x)\in \L^2_\sigma(\Gamma),\\&
			\f_2^{\pm}=\f_2^{\pm}(x,\u)=\v\pm \P_\mathcal K\u,\quad \text{for some fixed }\v\in (\boldsymbol I-\P_\mathcal K)\L^2_\sigma(\Gamma),\\& 
			\f_3^{\pm}=\f_3^{\pm}(\u)=\pm\u,\\&
			\f^{\pm}_4=\f_4^{\pm}(x,\u)=(\boldsymbol I-\P_\mathcal K)\u\pm\P_\mathcal{K}(\norma{x-p}\P_\mathcal K\u), \quad \text{for some fixed }p\in \Gamma,\\&
			\f_5=\f_5(x,\u)=(\I-\Pk)(\vert x\vert \u)-\u.
		\end{align*}
		Note that in the case of $\f_4$ it holds
		\begin{align*}
			\f_{4,K}^\pm(x,\u):=\pm\vert x-p\vert \Pk\u,\quad \f_{4,NK}^\pm(x,\u):=(\I-\Pk)\u,
		\end{align*} 
		whereas, in the case $\f_5$,
		\begin{align*}
			\f_{5,K}(x,\u):= -\Pk\u,\quad \f_{5,NK}(x,\u):=(\I-\Pk)(\vert x\vert \u-\u),
		\end{align*}
		and it clearly satisfies both case (ii) of Proposition \ref{lemmadec} as well as \eqref{extra2}, since
		$$
		\int_\Gamma \f_{5,NK}(x,\u)\cdot \u\leq C\norm{(\I-\Pk)(\vert x\vert \u)}\norm{(\I-\Pk)\u}+C\norm{(\I-\Pk)\u}^2\leq C\norm{(\I-\Pk)\u}^2+C\norm{\Pk\u}^2,
		$$
		recalling $\norm{\u}^2=\norm{\Pk\u}^2+\norm{(\I-\Pk)\u}^2.$
		\label{remf}
	\end{osse}
	\begin{osse}\rm (Discussion of  the examples)
		Concerning the forcing terms proposed in Remark
                \ref{remf}, we see that case (i) of Proposition
                \emph{\ref{lemmadec}} corresponds to $\f_1$. Then
                $\f_2^-,\f_3^-,\f_4^-,\f_5$ are related to case (ii),
                whereas $\f_2^+,\f_3^+,\f_4^+$ fall into case
                (iii). Concerning $\f_2^\pm,\f_3^\pm$, we observe that
                they both correspond to $\f_{K}=\pm\u_K$, so that,
                more  precisely  we deduce
		$$
		\frac12\frac{d}{dt}\norm{\u_K}^2=\norm{\u_K}^2,
		$$ 
		i.e.,
		$$
		\norm{\u_K(t)}^2=e^{\pm 2t}\norm{\u_K(0)}^2,
		$$
		so that $\u_K(t)\to 0$ as $t\to \infty$ in the cases $\f_2^-,\f_3^-$, whereas $\norm{\u_K(t)}\to +\infty$ as $t\to\infty$ for $\f_2^+,\f_3^+$.
		
		\noindent In conclusion, in the case $\f_K= \boldsymbol 0$, by uniqueness (see Proposition \ref{dec}) we simply have $$\u_K(t)\equiv \u_K(0),\quad \forall t\geq0.$$
	\end{osse}

 Under   assumptions \eqref{uk1}, \eqref{nega}, and
\eqref{extra2}, we can show that any weak solution instantaneously
regularizes and, in particular, all its norms are uniformly bounded by
the initial energy. Namely, we have  the following  (see Section \ref{sec:prel} for the proof)
\begin{teor}[Instantaneous regularization]\label{insta}
	Let assumptions \eqref{C1}--\eqref{C2} hold, together with \eqref{uk1},  \eqref{nega}, and \eqref{extra2}.
	Then the unique weak solution $\u$ instantaneously regularizes, i.e., $\u$ is such that, for any $\tau>0$,
	\begin{align}
		& \boldsymbol{u}\in BUC([\tau,\infty );\boldsymbol{H}^{1}(\Gamma )),\label{regg1}\\&
		\u \in H^1(t,t+1;\L^2_\sigma(\Gamma))\cap L^2(t,t+1,\H^2_\sigma(\Gamma)),\quad \forall t\geq \tau.\label{regg2}
	\end{align}
	The solution $\u$ is uniformly bounded in the above spaces
         just in terms of   $\tau$, $\f$, $\Gamma$,
        $\norm{\u_0}$ and  of the  the parameters of the
        problem. Moreover, if  $\u_0\in \H^1_\sigma(\Gamma)$, one can
        choose $\tau=0$ in \eqref{regg1}--\eqref{regg2}. In this case,
        the  above-mentioned bounds also depend  on $\norm{\u_0}_{\H^1(\Gamma)}$.
\end{teor}
\begin{osse}\rm (Initial datum)
	We point out that the instantaneous regularization in \eqref{regg1}--\eqref{regg2} strongly depends on the initial energy of the complete datum $\u_0$, meaning that also the effect of the  $\L^2$-norm of the $\L^2$-projection of $\u_0$ on the Killing field space $\mathcal K$, i.e., $\norm{\Pk\u_0}$, has an influence on the $L^\infty_t\H_\sigma^1$-norm of the complete solution $\u$. {This means that, even in the case $\ff=\zero$, the Killing and non-Killing components of $\u$ cannot be decoupled when obtaining higher-order \textit{a priori }bounds, in contrast to what happens at the level of the basic energy estimate. This is of course due to the convective term $(\u\cdot\nabla_\Gamma)\u$.}  \label{effect}
\end{osse}

\subsection{The notion of attractor for the dynamical system $(S(t),\Lts)$}\label{notion}
The natural object to look for in the case of a dissipative dynamical
system is the global attractor, which is usually defined as either the
maximal compact invariant set, the minimal closed set which uniformly
attracts all bounded sets, or the set of points on complete bounded
trajectories. For many dissipative systems,  and in particular for
 the Navier--Stokes equations in 2D bounded domains (see, e.g., \cite{Foias2}), these definitions are equivalent. For these standard dissipative systems (i.e., with bounded absorbing sets), the global attractor $\mathcal{A}$ is compact, invariant, and attracting, entailing in particular that the following properties hold
\begin{enumerate}
	\item If $\mathcal{F}$ is a bounded invariant set, then
          $\mathcal{F}\subset \mathcal A$, i.e., $\mathcal A$ is the
          maximal  bounded   invariant set.
	\item If $\mathcal G$ is a closed attracting set, then $\mathcal A\subset \mathcal G$, i.e., $\mathcal A$ is the minimal closed attracting set.
\end{enumerate}
It is easy to see that these two properties entail uniqueness of the
global attractor $\mathcal A$, since  it is itself   closed and bounded.

 For the  dynamical system $(S(t),\Lts)$  in study, 
the definition  asks for a modification.  In particular, 
the  compactness  requirement  (and, more in general, 
the  boundedness)  must be dropped, in consideration of 
the possible choices of the forcing term $\f$.  In the following,
  we  introduce suitable  notions of attractors, which to
some extent recover the properties (1) and (2) above.  Note on the
other hand that uniqueness may be  lost, since  the  attractor might be unbounded. 

The  basic features  of the global attractor  which  we
might aim  at preserving  in our new notion are the invariance
property  and   the fact that it is the minimal closed set which attracts all bounded sets (i.e., property (2) above). This means that a good candidate is any set $\mathcal{C}\subset \Lts$ such that
\begin{align*}
&S(t)\mathcal C=\mathcal C\quad\forall t\geq 0,
\end{align*}  
and, if there exists a closed set $\mathcal{D}$ such that 
$$
\lim_{t\to +\infty}\text{\rm dist}(S(t)B,\mathcal D)=0,\quad \forall B\subset \Lts, \quad B\text{ bounded set},
$$
then $\mathcal{C}\subset \mathcal{D}$.

As a preliminary result, we can immediately show that, if we simply
choose $\f_K$ to be independent of $\u$, i.e., $\f_K:\Gamma \to
\T\Gamma$, as in the standard case of 2D Navier--Stokes equations on
bounded domains (\cite{Foias2}), such a set $\mathcal C$ is 
necessarily  empty. To this aim, we first state the following
general  result. 
\begin{lemm}[Emptyness of invariant sets]
	Consider the dynamical system $(S(t),\Lts)$. If there exists a closed set $\mathcal{D}\subset \Lts$ such that 
	\begin{align}
	\lim_{t\to +\infty}\text{\rm dist}(S(t)B,\mathcal D)=0,\quad
          \forall B  \ \text{bounded  in} \ \Lts,\label{bdd}
	\end{align}
	and such that
	\begin{align}
	\lim_{t\to\infty}\inf_{\u\in 	\mathcal{D}}\norm{S(t)\u}=+\infty,
	\label{inf}
	\end{align}
	then any invariant set $\mathcal C\subset \Lts$  with 
        $\mathcal C\subset \mathcal D$  is empty.\label{empty}
\end{lemm} 
\begin{proof}
 Assume by contradiction that $\mathcal C\not=\emptyset$. Since by the
 invariance property $S(t)\mathcal{C}=\mathcal{C}$, for any $t\geq0$
  and  $\u\in  \mathcal{C}$ there exists a sequence $\{\u_n\}_n\subset \mathcal{C}$ such that
	$$
	\u=S(n)\u_n,\quad \forall n\in\N.
	$$
 Hence, using the fact that   $\mathcal C\subset \mathcal D$,
	\begin{align*}
		\infty>\norm{\u}=\norm{S(n)\u_n}\geq \inf_{\v\in \mathcal C}\norm{S(n)\v}\geq \inf_{\v\in \mathcal D}\norm{S(n)\v}\to +\infty,
	\end{align*}
	a contradiction.
\end{proof}
Exploiting this lemma and Proposition \ref{dec}, we immediately 
obtain  
the following proposition,  which applies to the case of 
$\f_K$ being constant in time.  
\begin{prop}\label{const}
	Under the assumptions of Propositions
        \emph{\ref{dec}--\ref{pp}} (in particular,  by   assuming \eqref{extra}), if $\f_K\not= \boldsymbol 0$ satisfies case (i) of Proposition \emph{\ref{lemmadec}}, i.e., it is independent of $\u\in \Lts$, then any invariant set $\mathcal C\subset \Lts$ which is also the minimal closed set uniformly attracting all bounded sets, is empty.
\end{prop}
\begin{proof}
	We first introduce the set $\mathcal B\subset \Lts$ as
	\begin{align}
		\mathcal B:=\left\{\v\in \Lts:\ \norm{(\boldsymbol I-\P_\mathcal K)\v}\leq \sqrt{\frac12+\omega},\quad \int_\Gamma \P_{\mathcal{K}}\v\cdot \f_K\geq 0 \right\},\label{BB}
	\end{align}
	where $\omega>0$ is given in \eqref{gr}. Thanks to \eqref{fd} and \eqref{gr}  we can prove that for any bounded set $B\subset \Lts$ there exists a finite $t_B>0$ such that 
	\begin{align}
		S(t)B\subset \mathcal B,\quad \forall t\geq t_B,\label{absorbing}
	\end{align} 
	i.e., $\mathcal B$ is an absorbing set. To see this, we set $\u_{NK}(t):=(\boldsymbol I-\P_K) S(t)\u_0$, and $\u_K(t):=\P_{\mathcal{K}}S(t)\u_0$, for $\u_0\in B$. By \eqref{gr}, there exists $t_{B,1}>0$, depending only on $\norm{(\boldsymbol I-\P_{\mathcal{K}}
		)B}:=\sup_{\v\in (\boldsymbol I-\P_{\mathcal{K}}
		) B}\norm{\v}\leq \sup_{b\in B}\|b\|=: \norm{B}\leq  C_B$, such that
		$$
		\norm{\u_{NK}(t)}\leq \frac12+\omega,\quad \forall t\geq t_{B,1}.
		$$
	Moreover, since $B$ is a bounded set, also $\norm{\Pk B}\leq C_B$, so that 
	$$
	\norma{ \int_\Gamma \f_K\cdot \v} \leq \norm{\f_K}\norm{\v}\leq C_B\norm{\f_K},
	$$
	and thus, from \eqref{fd}, 
	$$
	\int_{\Gamma}\f_K\cdot \u_K\geq -C_B\norm{\f_K}+t\norm{\f_K}^2\to +\infty \quad\text{as }t\to \infty,
	$$
	entailing that there exists $t_{B,2}>0$, depending only on $\norm{B}$, such that 
	$$
	\int_\Gamma \f_K\cdot \u_K\geq 0,\quad \forall t\geq t_{B,2}.
	$$
	Choosing $t_B:=\max\{t_{B,1},t_{B,2}\}$, we have shown \eqref{absorbing}, since this time does not depend on the specific $\u_0\in B$, but only on $\norm{B}$.
	
	Now, we show that the absorbing (closed) set $\mathcal D:=\mathcal B$  satisfies the assumptions of Lemma \ref{empty}. The fact that $\mathcal{B}$ satisfies \eqref{bdd} is straightforward from its absorbing property. Concerning property \eqref{inf}, we have, for any $\u_0\in \mathcal{B}$, 
	\begin{align*}
			\norm{\u_K(t)}^2=\norm{\u_K(0)}^2+t^2\norm{\f_K}^4+2t\norm{\f_K}^2\int_{\Gamma}\f_K(x)\cdot \u_{K}(x,0)\geq t^2\norm{\f_K}^4,
	\end{align*} 
	since on $\mathcal B$ the last summand is always nonnegative. Therefore, this implies
	\begin{align*}
		\norm{\P_{\mathcal{K}}S(t)\u_0}\geq t \|\f_K\|^2, \forall \u_0\in \mathcal{B}.
	\end{align*}
 Since  $\norm{S(t)\u_0}\geq \norm{\P_{\mathcal{K}}S(t)\u_0}$,
 this entails 
		\begin{align*}
		\inf_{\u_0\in \mathcal{B}}\norm{S(t)\u_0}\geq t\norm{\f_K}^2\to +\infty \quad\text{as } t\to \infty,
	\end{align*}
	i.e., property \eqref{inf} holds.  Moreover, since $\mathcal C$
        is an invariant set which is the minimal closed set attracting
        all bounded sets, it holds  that  $\mathcal C\subset \mathcal D=\mathcal B$.
	
	Therefore, Lemma \ref{empty} shows that $\mathcal
        C=\emptyset$, i.e., any invariant set $\mathcal C\subset \Lts$
        which is also the minimal closed set uniformly attracting all
        bounded sets is empty.  This concludes  
        the proof.
	\end{proof}

         In constrast with the classical results for 
        the Navier--Stokes equations in 2D bounded domains, here the case $\f_K$ independent of the solution $\u$ does not allow to give a satisfying notion of attractor. Therefore we will assume from now on that $\f_K$, if not exactly equal to zero, also depends on $\u$.

In general, the task to find an attractor with good properties is not
easy, and strongly depends on the characteristics of the forcing term
$\f_K$.  Following  
\cite{Chepyzov},
we are interested in the set of \textit{bounded-in-the-past} complete trajectories 
\begin{align}
\mathcal{J}:=\{\xi(0)\in \Lts:\ \xi \text{ is a bounded-in-the-past complete trajectory for $S(t)$}\},\label{J}
\end{align}
where $\xi : \R \to \Lts$ is a bounded-in-the-past complete trajectory of $S$ if $S(t)\xi(s) = \xi(t +s)$
for all $t\geq 0$  and $s\in\R$, and $\xi((-\infty, 0])$ is a bounded subset of $\Lts$. 

In the standard case of dissipative systems,  where all
trajectories are bounded, 
this set coincides with the set 
\begin{align}
	\mathcal{I}:=\{\xi(0)\in \Lts:\ \xi \text{ is a bounded complete trajectory for $S(t)$}\},\label{I}
\end{align}
where $\xi((-\infty, +\infty))$ is a bounded subset of $\Lts$, and the
global attractor exactly coincides with this set. In our case of the
system $(S(t),\Lts)$, it is not \textit{a priori} ensured that all
trajectories are bounded in the future, and thus in general we can
only expect $\mathcal I \subset \mathcal J$. In general,  the
interest in $\mathcal J$ as a descriptor of the long-time behavior of
the system is given  
by the following trivial lemma (see for instance \cite[Proposition
3]{Bortolan}), which we prove  here for completeness.  
\begin{lemm}[Properties of $\mathcal J$]\label{obv}
	The set $\mathcal J$ satisfies the following properties:
	\begin{enumerate}
		\item[A.] $\mathcal J$ is invariant, i.e., $S(t)\mathcal J=\mathcal J$ for any $t\geq 0$.
			\item[B.]  If a closed set $\mathcal D\subset \Lts$ attracts all bounded sets, then $\mathcal J\subset \mathcal D$.
			\item[C.] If $A\subset \Lts$ is a nonempty  bounded invariant set, then $A\subset \mathcal I\subset \mathcal J$ and thus $\mathcal J$ is nonempty.
	\end{enumerate}
\end{lemm}
\begin{proof}
	The invariance property of $\mathcal J$ is  easily
        checked. Let  
        $\v\in \mathcal J$, so that there exists a complete trajectory
        $\xi$ bounded-in-the-past with $\xi(0)=\v$. Note that
        $\xi_r:=\xi(\cdot+r)$ for any $r\in \R$ is still a
        bounded-in-the-past complete trajectory, so that, for any
        $t\geq0$, $S(t)\v=\xi_t(0)$ and thus by definition
        $S(t)\mathcal J\subset \mathcal J$. Analogously, since
        $\xi_{-t}(0)\in \mathcal J$ for any $t\geq 0$, it holds 
        that  $\v=S(t)\xi_{-t}(0)\in S(t)\mathcal J$, and thus $\mathcal J\subset S(t)\mathcal J$, i.e., $\mathcal J$ is invariant. 
	
	To prove assertion  B.,  let us  fix   $\v\in
        \mathcal J$  and let $\xi$ be  a complete trajectory
         which is  bounded-in-the-past with $\xi(0)=\v$. Consider $B:=\xi((-\infty,0])$, which is a bounded set containing $\v$ {and such that $B\subset S(t)B$ for all $t\geq0$}, and observe that 
	$$
	\text{\rm dist}(B,\mathcal D){\leq}\text{\rm dist}(S(t)B,\mathcal
        D)\to 0,\quad \text{as} \ t\to\infty.
	$$
	Since $\mathcal D$ is closed, this means $\v\in B\subset \mathcal D$, entailing $\mathcal J\subset \mathcal D$.
	 
 To  prove property  C.,  we recall that, by \cite[Lemma 1.4]{Robinsonbook}, a set $A\subset \Lts$ is invariant if and only if for any $\v\in A$ there exists a complete trajectory $\xi$ such that $\xi(\R)\subset A$. Therefore if we assume $A\subset \Lts$ to be bounded and invariant, then clearly $A\subset \mathcal I$. 
\end{proof}

	\begin{osse}\rm ($\mathcal I$ is invariant)
It is immediate to verify that also $\mathcal I$ is an invariant set, i.e., $S(t)\mathcal I=\mathcal I$, for any $t\geq 0$. 
Note also that, given any complete bounded or bounded-in-the-past trajectory $\xi$, it clearly holds $\xi(\R)\subset \mathcal I$ or $\xi(\R)\subset \mathcal J$, respectively.
\end{osse}
\begin{osse}\rm (Case $\mathcal J\not =\emptyset$)
	We point out that if we prove that $\mathcal J$ is nonempty,
        closed, and attracts all bounded sets of $\Lts$, then by
        property B. of  Lemma \emph{\ref{obv}}  we deduce that
        it is also the minimal closed attracting set. Then $\mathcal
        J$ satisfies properties (1) and (2) described above for the
        classical global attractor for dissipative dynamical systems:
        the only property which is missing  is the compactness (in
        particular,  the   boundedness). This property will be
        retrieved in a weaker sense: we will prove  that $\mathcal
        J$ is in many cases  \emph{bounded compact}, i.e., the intersection of $\mathcal J$ with any closed and bounded set is compact.      	
\end{osse}

A first nontrivial  task   is to prove that $\mathcal J$ is
nonempty.  Note that this is generally false.   Indeed, for instance in the case of $\f_K$ independent of the velocity it is immediate to see that $\mathcal J=\emptyset$:
\begin{lemm}[Case $\mathcal J =\emptyset$]
		Under the assumptions of Propositions \emph{\ref{dec}--\ref{pp}}, if $\f_K\not = \boldsymbol 0$ satisfies case (i) of Proposition \emph{\ref{lemmadec}}, i.e., it is independent of $\u\in \Lts$, then $\mathcal J=\emptyset$.
\end{lemm}
\begin{proof}
Since the set $\mathcal{B}$ introduced in \eqref{BB} is absorbing, then, by property B. of Lemma \ref{obv}, it holds $\JJ\subset \mathcal{B}$. Since $\JJ$ is invariant and $\JJ\subset \mathcal B$, we can apply Lemma \ref{empty}, with $\mathcal{D}=\mathcal{B}$ and $\mathcal C=\JJ$ and conclude that $\JJ=\emptyset$.
\end{proof}
Our objective then becomes to first show that $\mathcal J$ is
nonempty, and then further  ascertain the properties of  the
set $\mathcal J$ (or $\mathcal I$, in some cases),  such as 
being the minimal closed set attracting all bounded sets. 

We distinguish two cases, corresponding to cases (ii) and (iii) of Proposition \ref{lemmadec}. They will be treated with a completely different approach.

\subsection{Attractors for bounded trajectories}
\label{sec:attr}
 We first focus on  
case (ii) of Proposition \ref{lemmadec}.  Namely, we assume
\eqref{nega}, that is  $$\int_\Gamma \f_K(x, \u)\cdot \u\leq 0,\quad\text{for any } \u\in \Lts.$$
In this situation, it is clear from Propositions \ref{dec}--\ref{pp} that any trajectory is bounded, i.e., given any $\u_0\in \Lts$,  $S([0,+\infty))\u_0$ is bounded in $\Lts$, entailing that 
$$
\mathcal J=\mathcal I.
$$
In this case, we can fully characterize the set $\mathcal J$ (which is proven to be nonempty), and we  call this set \textit{$\sigma$-attractor}, in analogy to what has been first introduced in \cite{Foias} to deal with systems with similar dissipativity properties. This name comes from the fact $\mathcal J$ can be obtained as the countable union of nonempty, compact, finite-dimensional (with respect to the fractal dimension) sets. We recall that the fractal dimension of a compact set $A\subset \Lts$ is defined as $$
\text{dim}_{\Lts}(A)=\limsup_{\epsilon\to0}\frac{\log N(\epsilon)}{-\log \epsilon},
$$
and $N(\epsilon)$ is the minimum number of $\epsilon$-balls of $ \Lts$ necessary to cover $A$.

We now introduce the following notation: for any $r\geq0$, we define 
\begin{align}
	\mathbb B_r:=\{\v\in \Lts:\ \norm{\P_{\mathcal{K}}\v}\leq r\},
	\label{Ar}
\end{align}
which is a complete metric space if endowed with the distance induced by the norm of $\Lts$, i.e., the $\L^2(\Gamma)$-norm.
We will show that, for any $r\geq0$, the dynamical system
$(S(t),\mathbb B_r)$ admits a finite-dimensional global attractor
$\mathcal A_r$,  namely, fulfilling the following properties 
\begin{itemize}
	\item $\mathcal A_r$ is nonempty, compact, connected,  and
          of finite fractal dimension,  
	\item $\mathcal A_r$ is invariant, i.e., $S(t)\mathcal A_r=\mathcal A_r$, for any $t\geq0$,
	\item $\mathcal A_r$ is attracting, i.e., for any $B\subset \mathbb B_r$ bounded it holds 
	${\rm dist}(S(t)B,\mathcal A_r)\to 0$ as $t\to \infty$.
\end{itemize}
Clearly, from these properties we can also deduce that $\mathcal A_r$ is the maximal bounded invariant set, as well as the minimal closed attracting set. Moreover, it also holds
\begin{align}
	\mathcal A_r=\{\xi(0)\in \mathbb B_r:\ \xi \text{ is a bounded complete trajectory for $S(t)$ in }\mathbb B_r \},
	\label{bdd_traj}
\end{align}
which clearly entails $\mathcal A_r\subset \mathcal A_s$ if $0\leq
r\leq s$. Indeed, if $\xi$ is a complete bounded trajectory in
$\mathcal A_r$,  one has  that $\xi(\R)\subset \mathcal A_r\subset  \mathbb B_r$, but, since $\mathbb B_r\subset \mathbb B_s$, then also $\xi$ is a complete bounded trajectory in $\mathbb B_s$, and thus  $\xi(\R)\subset \mathcal A_s$.

We can now state our main result of this section (which is proven in Section \ref{bdda}).
\begin{teor}[$\sigma$-attractor]\label{thm:323}
		Under the assumptions \eqref{C1}--\eqref{C2},
                \eqref{uk1}, and \eqref{extra2}, if $\f_K$ is such that 
		\begin{align}
		\int_\Gamma \f_K(x, \u)\cdot \u\leq 0,\quad\text{for any } \u\in \Lts,
		\label{neg}
		\end{align}
		then the nonempty set $\mathcal J$ defined in \eqref{J}, coinciding with $\mathcal I$ defined in \eqref{I}, is the $\sigma$-attractor and enjoys the following properties:
		\begin{enumerate}
			\item $\mathcal J=\bigcup_{r\geq 0}\mathcal
                          A_r=\bigcup_{n\in \N}\mathcal A_n$, where
                          $\mathcal A_r\subset\H^1_\sigma(\Gamma)$ are
                          the global attractors to the dynamical
                          system $(S(t), \mathbb B_r)$ for any $r\geq
                          0$. Namely, $\mathcal A_r$ are nonempty,
                          compact,  of finite fractal dimension,
                          invariant, and connected. 
\item $\mathcal J\subset \H^1_\sigma(\Gamma)$ has empty interior in $\Lts$.
			\item $\mathcal J$ attracts all bounded sets. More precisely, if $B\subset \Lts$ is a bounded set, then there exists $r>0$ such that $${\rm dist}(S(t)B,\JJ)\leq {\rm dist}(S(t)B, \mathcal A_r)\to 0$$ as $t\to \infty$.
			\item $\mathcal J$ is invariant, i.e., $S(t)\mathcal J =\mathcal J$ for all $t\geq0$. 
			\item $\mathcal J$ is  minimal among 
                          closed sets attracting  
                          all bounded sets, which means that if there exists a closed set $\mathcal{D}$ such that 
			$$
			\lim_{t\to +\infty}\text{\rm dist}(S(t)B,\mathcal D)=0,\quad \forall B\subset \Lts, \quad B\text{ bounded set},
			$$
			then $\mathcal{J}\subset \mathcal{D}$.
		\end{enumerate}
\label{th1}
\end{teor}
%
%
%
%
%
\begin{osse}\rm (Complete trajectories do not cross)
	We point out that, thanks to the backward-uniqueness result of
        Theorem \ref{backuni}, together with the well-posedness given
        in Theorem \ref{thm1}, since any complete trajectory $\xi$
        belongs to $\JJ\subset \H^1_\sigma(\Gamma)$, and thus it is a
        strong solution to \eqref{t1}--\eqref{tangential}, 
        if any two
        complete trajectories intersect, then they must be
        identical.  In particular,  two distinct  complete
        trajectories $\xi$  cannot  intersect. This means than
        one could  require   the semigroup $S(t)$ on $\JJ$ 
        to  be defined  for negative times,  as well. 
\end{osse}

As a by-product of  Theorem \ref{thm:323},  we also have the
existence of a set $\mathcal M$ which can be defined as a
$\sigma$-exponential attractor, in the sense that it is the {countable} union of
the exponential attractors $\mathcal M_m$ for each dynamical system
$(S(t),\BBB_m)$, $m\in \mathbb N$. We recall that an exponential attractor
$\mathcal M_r$ (which is possibly  not  unique) for the dynamical system $(S(t),\BBB_r)$, {$r\geq0$}, enjoys the properties:
\begin{itemize}
	\item $\mathcal M_r$ is compact and of finite fractal dimension $N_r$, possibly increasing with $r\geq0$,
	\item $\mathcal M_r$ is positively invariant, i.e., $S(t)\mathcal M_r\subset \mathcal M_r$, for any $t\geq0$,
	\item $\mathcal M_r$ is exponentially attracting, i.e., for any $B\subset \mathbb B_r$ bounded it holds 
	$${\rm dist}(S(t)B,\mathcal M_r)\leq Q_r(\norm{B}_\Lts)e^{-
          \gamma_r  t},$$ where $Q_r>0$ is an increasing function
        of $\norm{B}_	\Lts$, only depending on $r$,  and
        $\gamma_r>0$ is a universal constant depending only on
        $r$. Here,  we have  used the notation   $\norm{A}_{\Lts}=\sup_{\u\in A}\norm{\u}$, for any bounded set $A$.
\end{itemize}

 As $\mathcal A_r\subset \mathcal M_r$, if $ \mathcal M_r$ has
finite fractal dimension, so does $\mathcal A_r$.  All the following results are proven in Section \ref{bdda}. We then have the following 
\begin{teor}[$\sigma$-exponential attractor]\label{thm:325}
	Under the assumptions \eqref{C1}--\eqref{C2}, \eqref{uk1}, and \eqref{extra2}, if $\f_K$ is such that 
$$\int_\Gamma \f_K(x, \u)\cdot \u\leq 0,\quad\text{for any } \u\in \Lts,$$
then there exists a set $\mathcal M$ enjoying the following properties:
\begin{enumerate}
	\item ${\mathcal M=\bigcup_{m\in \mathbb N} \mathcal M_m}$, {where $\mathcal M_m$ are exponential attractors for the dynamical system $(S(t),\BBB_m)$, for any $m\in \mathbb N$},
\item $\mathcal M$ is positively invariant, i.e., $S(t)\mathcal M\subset \mathcal M$, for any $t\geq0$,
\item {$\mathcal M$ is exponentially attracting, i.e., for any bounded set $B\subset \mathbb B_m$, for some $m\in \mathbb N$, it holds 
$${\rm dist}(S(t)B,\mathcal M)\leq {\rm dist}(S(t)B,\mathcal M_m)\leq Q_m(\norm{B}_\Lts)e^{-\gamma_m t},\quad \forall t\geq0,$$ where $Q_m>0$ is an increasing function of $\norm{B}_	\Lts$, only depending on $m$, and $\gamma_m>0$ is a universal constant depending only on $m$.}
\end{enumerate}
\end{teor}
{\begin{osse}(Nonuniqueness of the $\sigma$-exponential attractor)
		We point out here that, since for any $m\in \mathbb N$ an exponential attractor might be in general not unique, also the $\sigma$-exponential attractor we construct in the theorem above is not unique. Additionally, also the uncountable union $\bigcup_{r\geq 0}\mathcal M_r$, where $\mathcal M_r$ is an exponential attractor for the dynamical system $(S(t),\BBB_r)$, {$r\geq0$}, has the same properties (1)-(3) given in Thorem \ref{thm:325}. Neverheless, since in general {it is not ensured that $\mathcal M_r\subset \mathcal M_s$ for any $0\leq r\leq s$, the $\sigma$-exponential attractor $\mathcal M$ given by Theorem \ref{thm:325} is possibly smaller (and thus more desirable), since it clearly holds $\mathcal M\subset \bigcup_{r\geq 0}\mathcal M_r$. }
\end{osse}}
In conclusion, we can further refine the results of Theorem \ref{th1}
if we give some specific structure to the forcing term $\f$. In
particular, we have the following

\begin{teor}[Case $\f_K=0$]\label{spe}
		Under the assumptions \eqref{C1}--\eqref{C2} and
                \eqref{extra2}, if $\f_K= \boldsymbol 0$ or
                $\f_K(\cdot,\u)=\Pk((\v\cdot \nabla_\Gamma) \u_K)$,
                for $\v\in \L^2_\sigma(\Gamma)\cap \L^\infty(\Gamma)$,
                then the $\sigma$-attractor $\JJ$, defined in Theorem
                \emph{\ref{th1}}, possesses a  so-called 
                \emph{pancake-like} structure, { \cite{Foias}}, i.e., 
			$$
			\JJ=\bigcup_{r\geq 0}\widetilde{\AA}_r,
			$$
			where, having defined the complete metric space (endowed with $\Lts$ topology)
			\begin{align}
				\label{tB}
				\widetilde{\BBB}_r:=\{\u\in \Lts:\ \norm{\Pk\u}=r\},
			\end{align}
			the compact, invariant, attracting, 
                        finite-fractal dimensional  set
                        $\widetilde{\AA}_r$ 
                        is the global attractor of the dynamical system $(S(t),\widetilde{\BBB}_r)$, for any $r\geq 0$. 
			
			Furthermore, $\Pk \JJ=\KK$ and $\JJ$ is
                        closed, bounded compact, and  of finite
                        fractal dimension,  
                        i.e., the intersection of $\mathcal J$ with
                        any closed and bounded set $B\subset \Lts$ is
                        compact  and of  finite fractal dimension. In particular, there exists $r\geq 0$ such that $\mathcal A\cap B\subset \mathcal A_r$, where $\AA_{r}$ is the corresponding global attractor for the dynamical system $(S(t),\BBB_{r})$.
			
			In conclusion, if additionally $\f_{NK}=\boldsymbol 0$, then 
			$
			\JJ=\KK$.
\end{teor}
\begin{osse}\rm (Advective field)
	Notice that the theorem above also holds in case of the
        presence of an {external }advective field $\v\in \L^2_\sigma(\Gamma)\cap
        \L^\infty(\Gamma)$, in which the $\L^2$-norm of $\u_K$ does
        not change, but $\u_K$ is not trivially constant in time. The
        $\u_K$ component on the attractor $\JJ$ is then a rearranged
        version of the initial component $\Pk\u_0$. An example of this
        kind is when $\f=(\v\cdot\nabla_\Gamma)\u$, $\v\in
        \L^2_\sigma(\gam)\cap \L^\infty(\gam)$, which also satisfies
        \eqref{C1}--\eqref{C2} and \eqref{extra2}, thanks to Korn's
        inequality \eqref{Korn2}. 
\end{osse}
	In the case $\f_K$ is such that there exists a bounded
        absorbing set, we can say much more on the attractor. Notice
        that in this case we do not strictly need assumption
        \eqref{neg} on $\f_K$, since we can directly operate on the
        dynamical system $(S(t),\Lts)$, without any restriction. In
        this case, the set $\JJ$ is exactly the global attractor in the standard definition. Indeed, we have the following.
\begin{teor}[Properties of $\mathcal J$, II]
	\label{interesting}
	Let the assumptions \eqref{C1}--\eqref{C2}, \eqref{uk1}, and \eqref{extra} hold. If $\f_K$ is such that there exists $C\subset \Lts$ closed and bounded, such that, for any bounded set $B\subset \Lts$, there exists $t_B\geq0$ so that 
	$S(t)B\subset C$, for any $t\geq t_B$ (i.e., $C$ is a bounded
        absorbing set for the dynamical system $(S(t),\Lts)$), then
        the unique global attractor $\JJ$ defined in \eqref{J},
        coinciding with $\mathcal I$  from  \eqref{I}, is such that 
	\begin{enumerate}
		\item $\mathcal J\subset C$ is compact, connected, and of finite fractal dimension,
		\item $\mathcal J$ is invariant, i.e., $S(t)\mathcal J=\mathcal J$, for any $t\geq0$,
		\item $\mathcal J$ is attracting, i.e., for any $B\subset \Lts$ bounded it holds 
		$${\rm dist}(S(t)B,\mathcal J)\to 0,$$ as $t\to \infty$.
		\item $$\JJ \subset \left\{\u\in\Lts:\ \norm{(\I-\Pk)\u}\leq \sqrt{\frac 12 +\omega}\right\}.$$
	\end{enumerate}
	Moreover, there exists an exponential attractor $\mathcal M$, which is a compact, finite-dimensional, positively invariant set such that, for any $B\subset \Lts$ bounded it holds 
	$${\rm dist}(S(t)B,\mathcal M)\leq C(\norm{B}_\Lts)e^{-\gamma t},$$ where $C>0$ depends on $\norm{B}_{\Lts}$, and $\gamma>0$ is a universal constant. 
	
	\noindent If also $\Pk S(t)B\to \boldsymbol 0$ as $t\to \infty$ for any set $B\subset \Lts$ such that $\Pk B$ is bounded, then 
	\begin{align}\label{best}
		\Pk\JJ=\{\boldsymbol 0\}.
	\end{align}
\end{teor}
\begin{osse}\rm (Affine $\f$)
	The case $\f(x,\u):=-\u+c_0\v$, with $\v\in \KK$ and
        $c_0\geq0$, corresponding to $\f_K(x,\u)=-\Pk\u+c_0\v$ is
        clearly one example for which Theorem \emph{\ref{interesting}} 
        applies.  Notice that in this case assumption \eqref{neg}
        is \textit{not} satisfied. To see that the theorem holds,
        fixing $\u(t)=S(t)\u_0$, for some $\u_0\in \Lts$,  and
          multiplying the equation for $\Pk\u(t)$ by $\Pk\u(t)$
         we have 
	$$
	\frac12\frac{d}{dt}\norm{\Pk\u}^2= -\norm{\Pk\u}^2+\int_\Gamma c_0\v\cdot \Pk\u\leq -\norm{\Pk\u}^2+\frac12\norm{\Pk\u}^2+\frac {c_0^2}{2}\norm{\v}^2,
	$$
	entailing 
	\begin{align}
		\frac12\frac{d}{dt}\norm{\Pk\u}^2+\frac12\norm{\Pk\u}^2\leq \frac12c_0^2\norm{\v}^2 
	\label{p}
	\end{align}
	and thus by Gronwall's Lemma $$\norm{\Pk\u}\leq e^{-\zeta_0 t}\norm{\Pk\u_0}+\omega_0,$$ for some $\zeta_0,\omega_0>0$. Then the assumption of Theorem \ref{interesting} is satisfied for any set $B\subset \Lts$ such that $\Pk B$ is bounded. Indeed, recalling that also \eqref{gr} holds, it is easy to see that the set 
	$$
	 C:=\left\{\u\in\Lts:\ \norm{\Pk\u}\leq \sqrt{\frac12+\omega_0},\quad \norm{(\I-\Pk)\u}\leq \sqrt{\frac12+\omega}\right\}
	$$
is a bounded absorbing set for the dynamical system.	

	 \noindent When $c_0=0$, the global attractor $\JJ$ can be
         further characterized as in \eqref{best}, since 
         $\omega_0$ from \eqref{p} is $0$ in this case,  and thus the $\Lts$-norm of $\Pk\u$ is subject to an exponential decay to zero.

\end{osse}
\begin{osse}\rm (Interactions of Killing and non-Killing components)
	We point out that the assumptions of Theorem \ref{interesting}
        are satisfied also in some cases when we assume some
        interaction of the non-Killing component   $\f_{NK}$
        on the Killing component of the trajectory. For instance, let
        us assume $\f$ satisfying \eqref{C1}--\eqref{C2}, and
        $\f_{NK}$ satisfying \eqref{extra}. If we choose, for
        instance, $$\f(x,\u):=-\u+\vert (\I-\Pk)\u\vert \v,$$ for some
        $\v\in \L^\infty(\Gamma)$, then
        $\f_{NK}(x,\u)=-(\I-\Pk)\u+(\I-\Pk)(\vert (\I-\Pk)\u\vert
        \v)$, whereas $\f_K(x,\u)=-\Pk\u+\Pk(\vert (\I-\Pk)\u\vert
        \v)$.  Property  \eqref{gr}  immediately follows,  whereas for the Killing component of the trajectory, we have by standard estimates, together with \eqref{gr},
	\begin{align*}
		\frac12\frac{d}{dt}\norm{\Pk\u}^2&=-\norm{\Pk\u}^2+\int_\Gamma \vert (\I-\Pk)\u\vert \v\cdot \Pk\u
		\\&
		\leq -\norm{\Pk\u}^2+\norm{(\I-\Pk)\u}\norm{\Pk\u}\norm{\v}_{\L^\infty(\Gamma)}\\&
		\leq -\frac12\norm{\Pk\u}^2+\frac12\norm{(\I-\Pk)\u}^2\norm{\v}_{\L^\infty(\Gamma)}^2\\&
		\leq -\frac12\norm{\Pk\u}^2+\frac12(e^{-\zeta t}\norm{(\I-\Pk)\u_0}^2+\omega)\norm{\v}_{\L^\infty(\Gamma)}^2,
	\end{align*}
	where we recall $\v\in\L^\infty(\Gamma)$.  By 
        Gronwall's Lemma,  this entails   that 
	\begin{align*}
		\norm{\Pk\u(t)}^2\leq e^{-\frac12 t}\norm{\Pk\u_0}^2+\omega_1(1+q(t)\norm{(\I-\Pk)\u_0}^2),\quad \forall t\geq 0,
	\end{align*}
	where $q(t)\to 0$ as $t\to \infty$ and $\omega_1>0$. This means that the set 
		$$
	C:=\left\{\u\in\Lts:\ \norm{\Pk\u}\leq \sqrt{\frac12+\frac32\omega_1},\quad \norm{(\I-\Pk)\u}\leq \sqrt{\frac12+\omega}\right\}
	$$
	is an absorbing set for the whole dynamical system $(S(t),\Lts)$, and thus the assumptions of Theorem \ref{interesting} hold.
\end{osse}
\begin{osse}\rm (Case $\KK=\{\boldsymbol 0\}$)
	 If  $\KK=\{\boldsymbol 0\}$  one obviously has  that $\Pk S(t)\u_0=\boldsymbol 0$ for any $\u_0\in\Lts$, so that also in this case the global attractor $\JJ$ is characterized by Theorem \ref{interesting}, and satisfies \eqref{best}. 
\end{osse}

\subsection{Attractors in case of possibly unbounded trajectories}
\label{unbd}In this section, we  discuss   the case  when  trajectories can be unbounded, and thus in general we do not have the coincidence between $\mathcal J$ and $\mathcal I$. In particular, we assume that $$\int_\Gamma \f_K(x, \u)\cdot \u\geq 0,\quad\text{for any } \u\in \Lts\text{ such that }\norm{\Pk\u}\geq R_0,\quad\text{for some }R_0>0,$$ 
corresponding to case (iii) of Proposition \ref{lemmadec} (see in
particular Remark \ref{relax}). In this case we can characterize the
set $\mathcal J$,  which we already know that is nonempty and
which we call {\it unbounded attractor} following  
\cite{Chepyzov},  see also the  
recent works \cite{Carvalho,Bortolan}. Still, we will see that many
properties we expect from a global attractor are here preserved. The
main drawback, which is due the fact that there might be unbounded
trajectories, is the attraction property, which holds only on the
bounded (as $t\to \infty$) trajectories departing from a given bounded
set (see property (4) in  Theorem \ref{th2}  below). The other
properties, especially concerning the invariance and the property of
being the minimal invariant attracting closed set (in a slightly
weaker version) are here present as for Theorem \ref{th1}. Also in
this case, $\mathcal J$ is bounded compact.  Still,  we 
cannot recover 
the  precise  characterization of $\mathcal J$ as the countable union
of compact or  finite-fractal-dimensional  sets.

To be precise, we have the following main result, which is proven in Section \ref{unbdd}.

\begin{teor}[Properties of $\mathcal  J$, III]\label{thm:332}
	Under the assumptions \eqref{C1}--\eqref{C2}, together with assumptions \eqref{uk1} and \eqref{extra}, if $\f_K$ is such that 
	$$\int_\Gamma \f_K(x, \u)\cdot \u\geq 0,\quad\text{for any } \u\in \Lts \text{ such that }\norm{\Pk\u}\geq R_0,\quad \text{ for some }R_0>0,$$
	then the nonempty set $\mathcal J$ defined in \eqref{J} is the unbounded attractor, satisfying:
	\begin{enumerate}
		\item $\mathcal J=\overline{\bigcap_{t\geq0}S(t)Q}$, where $Q$ is given by 
		$$
		Q:=\overline{\bigcup_{t\geq 0}S(t)\left\{\v\in \Lts: \norm{(\boldsymbol I-\P_{\mathcal{K}})\v}\leq \frac12+\omega\right\}},
		$$
		where $\omega>0$ is given in  Proposition  \emph{\ref{pp}}.
		\item $\mathcal J$ is closed and invariant, i.e., $S(t)\mathcal J =\mathcal J$ for all $t\geq0$.
		\item $\JJ\subset \H^1_\sigma(\Gamma)$ and it has empty interior in the topology of $\Lts$. 
		\item $\mathcal J$ is bounded compact, i.e., the intersection of $\mathcal J$ with any closed and bounded set $B\subset \Lts$ is compact.
		\item If for some $R\geq 0$, $B \in\Lts$ bounded, and $t_1 > 0$ the sets $S(t)B \cap \{\v\in \Lts:\	\norm{\P_\mathcal K \v}\leq R\}$ are nonempty for every $t \geq  t_1$, then
		$$
		\lim_{t\to\infty}
		{\rm dist}(S(t)B \cap  \{\v\in \Lts:\	\norm{\P_\mathcal K \v}\leq R\}, \mathcal J ) = 0,$$
		and $\mathcal J$ is the minimal closed set with the above property.
\item $\P_\mathcal K \mathcal J = \mathcal K$.
	\end{enumerate}
	\label{th2}
\end{teor}
 Since in general in this case  one has that  $\mathcal
 I\subset \JJ$ but  the two may  not coincide, we specify here
 the properties of $\mathcal I$,  as well.
 
 \begin{coro}[Properties of $\mathcal I$]
 	Under the same assumptions of Theorem \ref{th2} the set $\mathcal I$ defined in \eqref{I} is such that it attracts all bounded sets with bounded trajectories, i.e., if a bounded set $B\in \Lts$ is such that there exist $t_1 > 0$ and $R \geq R_0$ with $S(t)B	\subset H_R$ for $t \geq  t_1$, where $H_R:=\{\u\in Q:\ \norm{\Pk\u}\leq R\}$ and $Q$ is defined in Theorem \ref{th1}, then
 	$$
 	\text{\rm dist} 
        (S(t)B,	\mathcal I)\to 0\quad\text{as }t\to \infty.
 	$$
 	Moreover, if $\u\in \JJ\setminus \mathcal I$, then
 	$$
 	\lim_{t\to \infty}\norm{S(t)\u}=+\infty.
 	$$
 \end{coro}
 \begin{proof}
 	The proof of this corollary is simply an application of \cite[Theorem 7, Lemma 6]{Carvalho}, which holds as a consequence of the validity of Lemma \ref{A1}.
 \end{proof}
 In conclusion, if we make some further assumptions on the behavior of
 the forcing term $\f_K$, we can also obtain that the set $\mathcal I$
 is nonempty, which is in general not trivial to prove. In particular,
 we have (see Section \ref{unbdd} for the proof)
 
  \begin{teor}[{Case $\f_K(\cdot,\mathbf 0)\equiv\mathbf 0 $}]\label{t1b}
 	Under the same assumptions of Theorem \emph{\ref{th2}}, assume additionally that $\f_K(x, \boldsymbol 0)= \boldsymbol 0$ for any $x\in \Gamma$. Then, by defining as $\mathcal A_0$ the global attractor to the system $(S(t),(\I-\Pk)\Lts)$, we have
 $$
 \mathcal A_0\subset \mathcal I,
 $$
 i.e., the set of complete bounded trajectories $\mathcal I$ defined in \eqref{I} is nonempty.
 	\label{tt}
 \end{teor}
 \begin{osse}\rm 
 	The theorem above can be easily generalized by assuming that there exists $\u_A\in \Lts$ such that $\f_K(x,\u_A)= \boldsymbol 0$ for any $x\in \Gamma$.
 \end{osse}
 \begin{osse}\rm 
 	Recalling Remark \emph{\ref{remf}}, examples of forcing terms
        satisfying the assumptions of Theorem \ref{t1b} are $\f_2^+,\
        \f_3^+,\ \f_4^+$,  and  $ \f_5$.
 \end{osse}

\section{ Well-posedness and backward uniqueness:  Proofs of Section \ref{main}}
\label{sec:pr1}
\subsection{Proof of Theorem \ref{thm1}}
	The proof can be carried out in many ways. Since we need to
        consider both weak and strong solutions, we first prove the
        short-time existence of a unique strong solution exploiting
        the $L^p-L^q$ maximal regularity properties of the surface
        Stokes operator with  variable   viscosity (see
        \cite[Lemma 7.4]{AGP2024b}). Then, by means of energy estimates,
        we show the global existence of both weak and strong solutions
        to the problem.

         Henceforth, we use the symbol $C$ to indicate a positive
        constant, possibly depending on data. The  value of $C$ may
        change from line to line. 
        
	\subsubsection{Local well-posedness of strong solutions}
	Let us first assume that $\u_0\in \H^1_\sigma(\Gamma)$. Given $T\in(0,T_0)$, {where $T_0>0$ is chosen arbitrarily large}, we introduce the space
	$$
	Z_T:=L^2(0,T;\H^2(\Gamma))\cap H^1(0,T;\L^2_\sigma(\Gamma)),
	$$
	equipped with the norm
	\begin{align}
		& \Vert \f \Vert_{Z_T}:=\Vert \f \Vert_{L^2(0,T; \H^2(\Gamma))}+\Vert \f\Vert_{ H^{1}(0,T;\L^2_\sigma(\Gamma))}+\Vert \f(0) \Vert_{\H^1_\sigma(\Gamma)},
	\end{align}
	and we recall that, as in \cite[Lemma 2]{Saal1}, by standard embeddings, there exists $C(T_0)>0$ such that
	\begin{align}\label{BUC}
		\norm{\u}_{BUC([0,T];\L^2_\sigma(\Gamma))}\leq C(T_0)\norm{\u}_{X_T},
	\end{align}
    for any $T\in(0,T_0)$.
	We also denote $\YT:=L^2(0,T;\L^2_\sigma(\Gamma))$  and
          we introduce the space 
	$$
	\XT:=\{\u\in Z_T:\ \u(0)=\u_0\},
	$$
	and we rewrite system \eqref{t1}--\eqref{tangential} as follows 
	\begin{align*}
		\mathcal{L}(\u)=\mathcal{F}(\u).
	\end{align*}
 Here,  the  linear operator $\mathcal{L}:\XT\to\YT$  is
defined  as 
	\begin{align}
		\label{L}
		\mathcal{L}(\v):=\partial_t\v-2\P_0\P_\Gamma\mathrm{div}_{\Gamma}(\nu(x)\E_\Gamma(\v))+\omega \v
	\end{align}
	for some arbitrary $\omega>0$, where $\P_0$ is the 
        Leray--Helmholtz  projector defined in \eqref{Leray} 
        and  the possibly nonlinear operator
        $\mathcal{F}:\XT\to\YT$  given by 
	\begin{align} 
		\label{F}
		\mathcal{F}(\v):=-\P_0(\v\cdot \nabla_\Gamma)\v+\P_0\f(\cdot,\u)+\omega\v.
	\end{align}
	First, we notice that the operator $\mathcal F$ is well defined. Indeed, by standard embeddings and the continuity of the operator $\P_0$,
	$$
	\norm{\P_0(\v\cdot \nabla_\Gamma)\v}_{\YT}\leq C\norm{\v}_{L^\infty(0,T;\L^4(\Gamma))}\norm{\v}_{L^2(0,T;\H^2(\Gamma))}\leq C\norm{\v}^2_{\XT}<+\infty,
	$$
	as well as, from assumptions \eqref{C1}--\eqref{C2}, 
	\begin{align*}
		\norm{\P_0\f(\cdot,\v)}_{\YT}&\leq C\norm{\f(\cdot,\v)-\f(\cdot,\boldsymbol 0)}_{L^2(0,T;\L^2(\Gamma))}+C\norm{\f(\cdot,\boldsymbol 0)}_{L^2(0,T;\L^2(\Gamma))}\\&\leq  C\norm{\v}_{\XT}+CT<+\infty,
	\end{align*}
	so that $\mathcal F(\u)$ is well defined for any $\u\in\XT$. 
	
	We aim now at proving that, under the assumptions of Theorem \ref{thm1}, there is a constant $C(T, R)>0$ such that
	\begin{align}
		\Vert\mathcal{F}(\v_1) -\mathcal{F}(\v_2)\Vert_{\YT}
		\leq C(T, R)\Vert(\v_1 - \v_2)\Vert_{\XT},
		\label{CR}\end{align}
	for all $\v_i\in \XT$ with $\Vert\v_i\Vert_{\XT}
	\leq R$, $R>0$, and $i = 1, 2$. Furthermore it holds
	$C(T, R) \to0$ as ${T} \to0$.
	
	To see this, let us fix $\v_i\in \XT$ with $\Vert\v_i\Vert_{\XT}
	\leq R$, $R>0$, split $\mathcal F$ in its summands, and 
        use   Gagliardo--Nirenberg's inequality, the embedding
        $\H^1(\Gamma)\hookrightarrow \L^4(\Gamma)$, and \eqref{BUC}
         to write 
	\begin{align*}
		&\norm{\P_0((\v_1\cdot \nabla_\Gamma)\v_1-(\v_2\cdot \nabla_\Gamma)\v_2)}_{\YT}\\&\leq \norm{\P_0((\v_1\cdot \nabla_\Gamma)\v_1-(\v_2\cdot \nabla_\Gamma)\v_2)}_{L^2(0,T;\L^2(\Gamma))}\\&
		\leq \norm{\v_1}_{L^\infty(0,T;\L^4(\Gamma))}\norm{\v_1-\v_2}_{L^2(0,T;\W^{1,4}(\Gamma))}+\norm{\v_1-\v_2}_{L^\infty(0,T;\L^4(\Gamma))}\norm{\v_2}_{L^2(0,T;\W^{1,4}(\Gamma))}\\&
		\leq
          C\norm{\v_1}_{\XT}\norm{\v_1-\v_2}_{L^\infty(0,T;\H^1(\Gamma))}\norm{\v_1-\v_2}_{L^1(0,T;\H^2(\Gamma))}\\&\quad
          +\norm{\v_1-\v_2}_{ \XT}\norm{\v_2}_{L^\infty(0,T;\H^1(\Gamma))}\norm{\v_2}_{L^1(0,T;\H^2(\Gamma))}\\&
		\leq C(R)T^\frac12\norm{\v_1-\v_2}_{\XT}
	\end{align*}
	Then, considering the forcing term $\f$, it holds, recalling assumption \eqref{C2},
	\begin{align*}
		\norm{\P_0(\f(\cdot,\v_1)-\f(\cdot,\v_2))}_{\YT}&\leq C\norm{\f(\cdot,\v_1)-\f(\cdot,\v_2)}_{L^2(0,T;\L^2(\Gamma))}\\&\leq C(R)T\norm{\v_1-\v_2}_{\XT},
	\end{align*}
	which, together with the previous estimate, leads to \eqref{CR}.
	
	We now  focus  on the operator $\mathcal L$. The fact
        that, for any ${T}>0$, this  operator  is invertible from $\YT$
        to $\XT$ can be deduced from \cite[Lemma 7.8]{AGP2024b}, by
        setting, in the notations of  that  lemma, $\nu\equiv
        2$, $\widetilde{\vphi}_0=\nu(\cdot)\in W^{1,\infty}(\Gamma)$,
        $\rho\equiv 1$ (so that the assumption
        $\vert\widetilde{\vphi}_0\vert\leq 1$ is not needed here,
        since $\rho_0\equiv1$).  Therefore, by the Bounded-Inverse
        Theorem, for any $T>0$ 
        there exists $C(T)>0$, possibly depending on $T$, such that 
	$$
	\norm{\mathcal{L}^{-1}}_{\mathcal{L}(\YT,\XT)}\leq C(T),\quad \forall T>0.
	$$
	It is then enough to show that the constant above does not
        actually change with $T$.  This   can be obtained by a simple extension argument (see for instance the proof of \cite[Lemma 7]{AWe}), leading to the fact that there exits $C(T_0)$  such that
	\begin{align}
		\norm{\mathcal{L}^{-1}}_{\mathcal{L}(\YT,\XT)} \leq  C(T_0),\quad \forall 0<T< T_0.
		\label{map}\end{align}
	
	We can now complete the existence proof.  We aim at
        solving via a fixed-point  argument 
         the equation 
	$$
	\u=\mathcal{L}^{-1}\mathcal{F}\u \quad\text{in }\XT,
	$$
    {for some }
$T\in(0,T_0)$.	First, we consider a generic $\overline{\v}\in \XT$. Then we fix $R>0$ such that it holds
	$\overline{\v}\in\overline{B}_{ R}^{X_{T_0}}(0)$, where
        $\overline{B}_R^{\XT}(0)$ is the closed ball of $\XT$ of
        radius $R$ centered at  $0$.  Clearly $\overline{B}_R^{X_{T_0}}(0)\subset \overline{B}_R^{\XT}(0)$, since $T\leq T_0$ and the norm of $\XT$ is nondecreasing. We also set $R$ such that $$\frac R2>\norm{\mathcal{L}^{-1}\mathcal{F}\overline{\v}}_{X_{T_0}}.$$ Note that in this way $R$ does not depend on $T$, since $T\leq T_0$, and it is also finite, since we have shown that $\mathcal F$ is well defined from $X_{T_0}$ to $Y_{T_0}$. Then one fixes $0<T< T_0$ (possibly depending also on $R$) such that the operator $\mathcal{L}^{-1}\mathcal{F}$ is a $(1/4)$-contraction mapping from $\XT$ to $\XT$. This is possible thanks to \eqref{CR} and \eqref{map}, since it holds
	\begin{align*}
		\Vert\mathcal{L}^{-1}\mathcal{F}(\v_1) -\mathcal{L}^{-1}\mathcal{F}(\v_2)\Vert_{\XT}
		\leq \norm{\mathcal{L}^{-1}}_{\mathcal{L}(\YT,\XT)}C(T, R)\Vert(\v_1 - \v_2)\Vert_{\XT}\leq C(T_0) C(T, R)\Vert(\v_1 - \v_2)\Vert_{\XT},
	\end{align*}
	and thus we choose $T$ sufficiently small so that $C(T_0) C(T, R)\leq \frac14$. Then, thanks to the estimates, one shows that $\mathcal{L}^{-1}\mathcal{F}$ is well defined from $\overline{B}_R^{\XT}(0)$ to itself. Indeed, for any $\v\in \overline{B}_R^{\XT}(0)$ we have
	\begin{align*}
		\norm{\mathcal{L}^{-1}\mathcal{F}\v}_{\XT}&\leq 	\norm{\mathcal{L}^{-1}\mathcal{F}\v-\mathcal{L}^{-1}\mathcal{F}\overline{\v}}_{\XT}+	\norm{\mathcal{L}^{-1}\mathcal{F}\overline{\v}}_{\XT}\\&
		\leq
          \frac14\norm{\v-\overline{\v}}_{\XT}+\norm{\mathcal{L}^{-1}\mathcal{F}\overline{\v}}_{X_{T_0}}
          < R,
	\end{align*}
	since $\v,\overline{\v}\in \overline{B}_R^{\XT}(0)$.
	
	Thus, by Banach fixed point theorem applied on $\mathcal
        {L}^{-1}\mathcal F: \overline{B}_R^{\XT}(0)\to
        \overline{B}_R^{\XT}(0)$, there exists a unique solution
        $\u\in \overline{B}_R^{\XT}(0)\subset \XT$ to the problem
        under study. By a standard argument it is also easy to show
        that the solution  $\u\in\overline{B}_R^{\XT}(0)\subset\XT$ we
        just found is unique in $\XT$. Indeed, let us assume that
        there exists another solution $\v\in  \XT$. Then consider
        $\tilde{R}>0$ larger than $R$ used in the previous argument,
        so that $\v\in\overline{B}_{\tilde{R}}^{\XT}(0)$. Then, 
        by   repeating the same argument we deduce that there exists $T_1(\tilde{R})\in(0,T_0)$ such that the solution in $\overline{B}_{\tilde{R}}^{X_{T_1}
	}(0)$ is unique on $[0,T_1(\tilde{R})]$ and coincides with $\v_{\vert [0,T_1]}$. Thus, since also $\u_{\vert [0,\tT_1]}\in \overline{B}_{\tilde{R}}^{X_{T_1}
	}(0)$,  it is immediate to infer $\v_{\vert [0,T_1]}=\u_{\vert [0,T_1]}$. If $T_1>T$ we are done, otherwise we can restart the same argument on the interval $[T_1,T_2]$, for some $T_2\in(T_1,T_0)$, and repeat the iterative continuation argument until we see that the identity holds on the entire interval $[0, T]$. The procedure terminates in a finite number of steps, since the time step in the contraction argument only depends on the radius of the ball containing $\v$, i.e., $\tilde{R}$. 
	
	In conclusion, the pressure can be retrieved by standard arguments exploiting De Rahm's Theorem.  
	\subsubsection{Existence of a global strong solution}
	In order to prove the existence of a global strong solution,
        we can use a standard  continuation  
        argument. Given $\u_0\in \H^1_\sigma(\Gamma)$, we know from the previous section that there exists a local strong solution $\u$. Let us assume that the maximal time $T_m$ of existence is finite. Then it holds
	\begin{align*}
		\u\in H^1_{loc}(0,T_m;\L^2_\sigma(\Gamma))\cap L^2_{loc}(0,T_m;\H^2(\Gamma)),
	\end{align*}
	entailing $\u\in BUC([0,T_m);\H^1(\Gamma))$.
    
   \noindent  By means of suitable energy estimates we now show that $\u\in C([0,T_m];\H^1(\Gamma))$, which is a contradiction, since then we can define $\u(T_m)$ and extend the maximal existence interval of the solution outside $[0,T_m]$. 
	
        First, we multiply \eqref{t1} by $\u$ and integrate over $\Gamma$. After some integration by parts we obtain the energy inequality
	\begin{align*}
			\frac12\frac{d}{dt}\norm{\u}^2+2\nu_*\int_\Gamma \norma{\E_\Gamma(\u)}^2\leq \frac12\frac{d}{dt}\norm{\u}^2+2\int_\Gamma \nu\norma{\E_\Gamma(\u)}^2=\int_\Gamma \f(\cdot,\u)\cdot \u.
	\end{align*}
	Recalling assumptions \eqref{C1}--\eqref{C2}, we obtain, by Young's inequality,
	\begin{align*}
	\int_\Gamma \f(\cdot,\u)\cdot \u\leq \norm{\f(\cdot,\u)-\f(\cdot,\boldsymbol 0)}\norm{\u}+\norm{\f(\cdot,\boldsymbol 0)}\norm{\u}\leq  C_1\norm{\u}^2+C_2\norm{\u}\leq C(\norm{\u}^2+1),
	\end{align*}
	so that we infer \begin{align}
			\frac12\frac{d}{dt}\norm{\u}^2+2\nu_*\int_\Gamma \norma{\E_\Gamma(\u)}^2\leq C(\norm{\u}^2+1),\quad \forall t<T_m.
			\label{energia}
	\end{align}
 From  Korn's inequality \eqref{Korn1} we infer by Gronwall's Lemma that
	\begin{align}
		\norm{\u}_{L^\infty(0,T_m;\L^2_\sigma(\Gamma))}+\norm{\u}_{L^2(0,T_m;\H^1(\Gamma))}\leq C(T_m).
		\label{reg1}
	\end{align}
	We can now pass to show higher-order estimates, namely, we multiply \eqref{t1} by $\partial_t\u$ and obtain, after integration by parts,
	\begin{align*}
		\frac{d}{dt}\int_\Gamma \nu \norma{\E_\Gamma(\u)}^2+\norm{\partial_t\u}^2+\int_\Gamma (\u \cdot\nabla_\Gamma)\u\cdot \partial_t\u=\int_\Gamma \f(\cdot,\u)\cdot \partial_t\u.
	\end{align*} 
	Then, we have, by Gagliardo--Nirenberg's and Korn's inequalities, recalling \eqref{reg1} and $\nu\geq \nu_*$,
	\begin{align*}
	&	\norma{\int_\Gamma (\u \cdot\nabla_\Gamma)\u\cdot \partial_t\u}\\&\leq \norm{\u}_{\L^4(\Gamma)}\norm{\nabla_\Gamma\u}_{\L^4(\Gamma)}\norm{\partial_t\u}\\&\leq C\norm{\u}^\frac12\norm{\u}_{\H^1(\Gamma)}\norm{\u}_{\H^2(\Gamma)}^\frac12\norm{\partial_t\u}
	\\&\leq 
	C\norm{\u}^\frac12(\norm{\u}+\norm{\E_\Gamma(\u)})\norm{\u}_{\H^2(\Gamma)}^\frac12\norm{\partial_t\u}\\&
	\leq  C(T_m)   (1+\norm{\E_\Gamma(\u)})\norm{\u}_{\H^2(\Gamma)}^\frac12\norm{\partial_t\u}\\&
	\leq C(T_m)\norm{\E_\Gamma(\u)}^2\left(\int_\Gamma \nu\norma{\E_\Gamma(\u)}^2\right)+\frac12\norm{\u}_{\H^2(\Gamma)}^2+\frac14\norm{\partial_t\u}^2+C(T_m).
	\end{align*} 
	Moreover, recalling \eqref{C1}--\eqref{C2} and \eqref{reg1}
         it holds 
	\begin{align*}
		&\norma{\int_\Gamma \f(\cdot,\u)\cdot \partial_t\u}\leq \norm{\f(\cdot,\u)}\norm{\partial_t\u}\\&\leq \norm{\f(\cdot,\u)-\f(\cdot,\boldsymbol 0)}\norm{\partial_t\u}+\norm{\f(\cdot,\boldsymbol 0)}\norm{\partial_t\u}\\&
		\leq C\norm{\u}\norm{\partial_t\u}+\norm{\f(\cdot,\boldsymbol 0)}\norm{\partial_t\u}\\&
		\leq C(T_m)+\frac14\norm{\partial_t\u}^2.
	\end{align*}
	Therefore, we can write
	\begin{align}
			\frac{d}{dt}\int_\Gamma \nu \norma{\E_\Gamma(\u)}^2+\frac12\norm{\partial_t\u}^2\leq C(T_m)\norm{\E_\Gamma(\u)}^2\left(\int_\Gamma \nu\norma{\E_\Gamma(\u)}^2\right)+C(T_m)+\frac12\norm{\u}_{\H^2(\Gamma)}^2.\label{time}
	\end{align}
	In order to close the estimate, we  use  \cite[Lemma 7.4]{AGP2024b}. Indeed, we can rewrite equation \eqref{t1} in weak formulation as
	\begin{align*}
		2\int_\Gamma \nu \E_\Gamma(\u):\E_\Gamma(\boldsymbol \eta)+\omega\int_\Gamma \u \cdot \boldsymbol \eta=\int_\Gamma \widetilde{\f}\cdot \boldsymbol \eta,\quad \forall \boldsymbol \eta \in \L^2_\sigma(\Gamma)\cap \H^1(\Gamma),
\end{align*}
	 for some $\omega>0$ and for 
	 $$
	 \widetilde{\f}:=\f(\cdot,\u)-(\u\cdot \nabla_\Gamma)\u-\partial_t\u.
	 $$
	Thus, we can apply  \cite[Lemma 7.4]{AGP2024b} with, in the notations of the lemma, $t=0$, $\Gamma(0)=\Gamma$, $\nu\equiv 2$, $\varphi_0=\nu\in W^{1,\infty}(\Gamma)$, $\f=\widetilde{\f}$, and obtain, similarly as above, recalling \eqref{C1}--\eqref{C2}, and \eqref{reg1},
	\begin{align*}
	&\norm{\u}_{\H^2(\Gamma)}\leq C(\omega)\norm{\widetilde{\f}}\\&\leq C(\norm{\f(\cdot,\u)}+\norm{(\u\cdot\nabla_\Gamma)\u}+\norm{\partial_t\u})\\&
	\leq C(\norm{\f(\cdot,\u)-\f(\cdot,\boldsymbol 0)}+\norm{\f(\cdot,\boldsymbol 0)}+\norm{\u}_{\L^4(\Gamma)}\norm{\nabla_\Gamma\u}_{\L^4(\Gamma)}+\norm{\partial_t\u})\\&
	\leq C\norm{\u}+C\norm{\f(\cdot,\boldsymbol 0)}+C\norm{\u}_{\L^4(\Gamma)}\norm{\u}_{\H^1(\Gamma)}^\frac12\norm{\u}_{\H^2(\Gamma)}^\frac12 +C_0\norm{\partial_t\u}	\\&
	\leq C+C_0\norm{\partial_t\u}+\frac12 \norm{\u}_{\H^2(\Gamma)},
		\end{align*}
	for some $C,C_0>0$. This entails, taking the squares,
	$$
	\norm{\u}_{\H^2(\Gamma)}^2\leq C+4C_0^2\norm{\partial_t\u}^2,
	$$
	so that, summing this inequality multiplied by
        $\gamma:=\frac1{16C_0^2}$  to inequality  \eqref{time}, we get
	\begin{align}
		\label{final}
			\frac{d}{dt}\int_\Gamma \nu \norma{\E_\Gamma(\u)}^2+\gamma \norm{\u}_{\H^2(\Gamma)}^2+\frac14\norm{\partial_t\u}^2\leq C(T_m)\norm{\E_\Gamma(\u)}^2\left(\int_\Gamma \nu\norma{\E_\Gamma(\u)}^2\right)+C(T_m),\quad \forall t<T_m.
				\end{align}
				Recalling that, due to \eqref{reg1}, $\u \in L^2(0,T_m;\H^1_\sigma(\Gamma))$, we can apply Gronwall's Lemma and infer
				\begin{align}
					\label{more_reg}
					\norm{\u}_{L^2(0,T_m;\H^2(\Gamma))}+\norm{\u}_{L^\infty(0,T_m;\H^1_\sigma(\Gamma))}+\norm{\u}_{H^1(0,T_m;\L^2_\sigma(\Gamma))}\leq C(T_m),
				\end{align}
				which clearly entails also 
				\begin{align*}
					\norm{\u}_{C([0,T_m];\H^1(\Gamma)\cap \L^2_\sigma(\Gamma))}\leq C(T_m),
				\end{align*}
				allowing to conclude the 
                                continuation  
                                argument, as outlined at the beginning of the proof. This concludes the argument to show that there exists a global in time strong solution.
	\subsubsection{Existence of a global weak solution} 
	The proof of  the  existence of a global weak solution
        is now straightforward.  Let us fix  $\u_0\in \L^2_\sigma(\Gamma)$. Due to the density of $\H^1_\sigma(\Gamma)$ in $\L^2_\sigma(\Gamma)$, as $\varepsilon\to0$, we can approximate the initial datum by a sequence $\{\u_0^\varepsilon\}\subset \H^1_\sigma(\Gamma)$ such that $\u_0^\varepsilon\to \u_0$ in $\L^2_\sigma(\Gamma)$. Then, for any $\varepsilon>0$ we have just shown that there exists a global strong solution $\u^\varepsilon$, which also satisfies the following energy estimate (see \eqref{energia})
	\begin{align}
		\frac12\frac{d}{dt}\norm{\u^\varepsilon}^2+2\nu_*\int_\Gamma \norma{\E_\Gamma(\u^\varepsilon)}^2\leq C(\norm{\u^\varepsilon}^2+1),\quad \forall t\geq 0.
		\label{energia2}
	\end{align}
	Applying Gronwall's Lemma  and  recalling that $\norm{\u_0^\varepsilon}\leq C$ uniformly in $\varepsilon$, we infer the following uniform estimates for any $T>0$
	\begin{align}
	\norm{\u^\varepsilon}_{L^\infty(0,T;\L^2_\sigma(\Gamma))}+\norm{\u^\varepsilon}_{L^2(0,T;\H^1_\sigma(\Gamma))}\leq C(T),\quad \forall \varepsilon>0.
	\label{unif1}
	\end{align}
We now obtain an estimate for $\partial_t\u^\varepsilon$, namely we consider $\v\in \H^1_\sigma(\Gamma)$ and observe that 
\begin{align*}
\int_\Gamma \partial_t\u^\varepsilon\cdot \v=-\int_\Gamma (\u^\varepsilon\cdot \nabla_\Gamma)\u^\varepsilon\cdot \v+\int_\Gamma \f(\cdot,\u^\varepsilon)\cdot\v.
\end{align*} 	
 By   Gagliardo--Nirenberg's inequalities and \eqref{unif1} 
we get 
\begin{align*}
&\norma{\int_\Gamma (\u^\varepsilon\cdot \nabla_\Gamma)\u^\varepsilon\cdot \v}=\norma{\int_\Gamma (\u^\varepsilon\cdot \nabla_\Gamma)\v\cdot \u^\varepsilon}\\&\leq \norm{\u^\varepsilon}_{\L^4(\Gamma)}^2\norm{\v}_{\H^1(\Gamma)}\leq C\norm{\u^\varepsilon}\norm{\u^\varepsilon}_{\H^1(\Gamma)}\norm{\v}_{\H^1(\Gamma)}\leq C\norm{\u^\varepsilon}_{\H^1(\Gamma)}\norm{\v}_{\H^1(\Gamma)}.
	\end{align*}
Moreover, recalling assumptions \eqref{C1}--\eqref{C2}  we can
handle the $\f$-term as follows 
\begin{align*}
	&\norma{\int_\Gamma \f(\cdot,\u^\varepsilon)\cdot\v}\leq \norm{\f(\cdot,\u^\varepsilon)}\norm{\v}\\&
	\leq  \norm{\f(\cdot,\u^\varepsilon)-\f(\cdot,\boldsymbol 0)}\norm{\v}+\norm{\f(\cdot,\boldsymbol 0)}\norm{\v}\\&
	\leq C\norm{\u^\varepsilon}\norm{\v}+\norm{\f(\cdot,\boldsymbol 0)}\norm{\v}
	\leq C(1+\norm{\u_\varepsilon})\norm{\v}.
\end{align*} 	
	Therefore, since by \eqref{unif1} we know that $\norm{\u_\varepsilon}_{L^2(0,T;\H^1_\sigma(\Gamma))}\leq C(T)$ uniformly in $\varepsilon$, we infer that, for any $T>0$, 
	\begin{align*}
		\norm{\partial_t\u}_{L^2(0,T; \H^1_\sigma(\Gamma)')}\leq C(T),\quad \forall \varepsilon>0.
	\end{align*}
	This result, together with \eqref{unif1} allows to deduce by
        standard compactness arguments that there exists
        $\u:\Gamma\times[0,\infty)\to \T\Gamma$ such that, up to 
        not relabeled  subsequences, for any $T>0$, 
	\begin{align*}
		&\u^\varepsilon\overset{*}{\rightharpoonup} \u,\quad \text{in }L^\infty(0,T;\L^2_\sigma(\Gamma)),\\&
			\u^\varepsilon \rightharpoonup \u,\quad\text{ in }L^2(0,T;\H^1_\sigma(\Gamma)),\\&
			\partial_t\u^\varepsilon\rightharpoonup \partial_t\u,\quad \text{ in }L^2(0,T;{\H^1_\sigma}(\Gamma)'),
	\end{align*}
	as $\varepsilon\to 0$. This also entails, by Aubin-Lions Lemma,
	$$
	\u^\varepsilon\to \u,\quad \text{ in }L^2(0,T;\H^s(\Gamma)),\quad \forall s\in[0,1),\quad\text{and almost everywhere in }\Gamma\times[0,T].
	$$
	These convergences are enough to pass to the limit in the equations satisfied by $\u^\varepsilon$ and conclude that $\u$ is a global weak solution to \eqref{t1}--\eqref{tangential}.

	\subsubsection{Weak uniqueness and continuous-dependence estimate}
	
	Thanks to the first part of Theorem \ref{thm1}, given $\u_{0,1},\u_{0,2}\in \L^2_\sigma(\Gamma)$ there exist two weak solutions $\u_1,\u_2$ on $[0,\infty)$ to \eqref{t1}--\eqref{tangential} departing from those initial data. Their regularity is enough to perform rigorously the next computations, leading to \eqref{contdep2}. In particular, let us notice that $\u:=\u_1-\u_2$ satisfies
	\begin{align}
		&\label{t2}\int_\Gamma\langle\partial_t\u,\v\rangle_{\H^1_\sigma(\Gamma)',\H^1_\sigma(\Gamma)}+\int_\Gamma(\u_1\cdot \nabla_\Gamma)\u\cdot \v+\int_\Gamma(\u\cdot \nabla_\Gamma)\u_2\cdot \v+2\int_\Gamma \nu\boldsymbol\varepsilon_\Gamma(\u):\E_\Gamma(\v)\\&=\int_\Gamma(\f(\cdot,\u_1)-\f(\cdot,\u_2))\cdot \v,\quad \forall \v\in \H^1_\sigma(\Gamma)\\nonumber\\&
		\u(0)=\u_{0,1}-\u_{0,2},
		\label{tangential2222}
	\end{align} 
	Let us now set $\v=\u$ in \eqref{t2}.  Recalling  that
        $\nu\geq \nu_*>0$,  we obtain 
	\begin{align*}
		\frac12\frac{d}{dt}\norm{\u}^2+2\nu_*\norm{{\boldsymbol \varepsilon}_\Gamma(\u)}^2\leq-\int_\Gamma (\u\cdot\nabla_\Gamma) \u_2\cdot \u +\int_\Gamma (\f(\cdot,\u_1)-\f(\cdot,\u_2))\cdot \u.
	\end{align*}
	Note that, by H\"{o}lder's, Gagliardo--Nirenberg's
        inequalities,  and  Korn's  inequalities we have
        that  
	\begin{align*}
		\norma{\int_\Gamma (\u\cdot\nabla_\Gamma) \u_2\cdot \u}&\leq \norm{\u}_{\L^4(\Gamma)}\norm{\nabla_\Gamma \u_2}\norm{\u}_{\L^4(\Gamma)}\\&\leq C\norm{\u}\norm{\u}_{\H^1(\Gamma)}\norm{\u_2}_{\H^1(\Gamma)}\\&\leq C\norm{\u}(\norm{\u}+\norm{{\boldsymbol \varepsilon}_\Gamma(\u)})\norm{\u_2}_{\H^1(\Gamma)}\\&
		\leq 
		C(1+\norm{\u_2}_{\H^1(\Gamma)}^2)\norm{\u}^2+\nu_*\norm{{\boldsymbol \varepsilon}_\Gamma(\u)}^2.
	\end{align*}
	 Moreover,   exploiting assumption \eqref{C2} on $\f$, we deduce, again by H\"{o}lder's, Gagliardo--Nirenberg's, and Korn's inequalities  
	\begin{align*}
		\norma{\int_\Gamma (\f(\cdot,\u_1)-\f(\cdot,\u_2))\cdot \u}&\leq \norm{\f(\cdot,\u_1)-\f(\cdot,\u_2)}\norm{\u}
		\\&
		\leq C\norm{\u}^2.
	\end{align*}
	Putting all these estimates together we end up with
	\begin{align*}
		\frac12\frac d{dt}\norm{\u}^2+\nu_*\norm{{\boldsymbol \varepsilon}_\Gamma(\u)}^2\leq C(1+\norm{\u_2}_{\H^1(\Gamma)}^2)\norm{\u}^2,
	\end{align*}
	recalling that the regularity of weak solutions implies that
        $\u_i\in L^2(0,T;\H^1_\sigma(\Gamma))$, $i=1,2$.  Using
         the embedding $\H^1(\Gamma)\hookrightarrow \L^4(\Gamma)$, we can apply Gronwall's lemma and deduce \eqref{contdep2}.
\subsubsection{Instantaneous regularization}
	An immediate consequence of the above results is the instantaneous regularization of a weak solution. Let us assume $\u_0\in \L^2_\sigma(\Gamma)$. Then there exists a unique global weak solution $\u$ departing from that datum. Let us fix an arbitrary $\tau>0$. Since $\u\in L^2(0,T;\H^1_\sigma(\Gamma))$ for any $T>0$, there exists $\tau_1\in(0,\tau]$ such that $\u(\tau_1)\in \H^1_\sigma(\Gamma)$. Therefore, by the first part of Theorem \ref{thm1}, there exists a unique global strong solution $\u_1$ departing from $\u(\tau_1)$. Since the solutions are unique, it holds $\u\equiv \u_1$ on $[\tau_1,+\infty)$, and thus the solution $\u$ becomes instantaneously strong. Let us fix $\tau>0$. In order to prove \eqref{regulari} we can repeat, for any $t\geq \tau$, the estimates leading to \eqref{final}, obtaining, for any $T>0$,
		\begin{align}
		\label{finalb}
		\frac{d}{dt}\int_\Gamma \nu \norma{\E_\Gamma(\u)}^2+\gamma \norm{\u}_{\H^2(\Gamma)}^2+\frac14\norm{\partial_t\u}^2\leq C(T)\norm{\E_\Gamma(\u)}^2\left(\int_\Gamma \nu\norma{\E_\Gamma(\u)}^2\right)+C(T),\quad \forall t\geq \tau,
	\end{align}
	where the constants $C$ depend on $T$ due to the energy estimates leading to \eqref{regga}. Multiplying by $t$ this inequality we are led to 
		\begin{align*}
		&\frac{d}{dt}t\int_\Gamma \nu \norma{\E_\Gamma(\u)}^2+\gamma t\norm{\u}_{\H^2(\Gamma)}^2+\frac t4\norm{\partial_t\u}^2\\&\leq C(T)t\norm{\E_\Gamma(\u)}^2\left(\int_\Gamma \nu\norma{\E_\Gamma(\u)}^2\right)+C(T)t\frac12+\int_\Gamma \nu \norma{\E_\Gamma(\u)}^2\\&
		\leq C(T)T\norm{\E_\Gamma(\u)}^2\left(\int_\Gamma \nu\norma{\E_\Gamma(\u)}^2\right)+C(T)T+\frac12\int_\Gamma \nu \norma{\E_\Gamma(\u)}^2,\quad \forall t\geq \tau.
	\end{align*}
	Since then $\norm{\u}_{L^2(0,T;\H^1_\sigma(\Gamma))}\leq C(T)$ by \eqref{regga}, we can apply Gronwall's Lemma and deduce \eqref{regulari}, recalling also Korn's inequality. The proof is concluded.
\subsection{Proof of Theorem \ref{backuni}}
In order to prove the theorem, we follow a somehow similar argument as
in \cite{Foias}.  There,  one needs to show that the Dirichlet
quotient $\lambda := \frac{\norm{\u}_{\H^1(\Gamma)}^2}{\norm{\u}^2}$
grows at most exponentially in time, exploiting the dissipativity
properties of the system. In this case, due to the lack of
dissipativity,  a suitable quantity to consider is  
\begin{align}
	\label{qty}
	\Lambda:=\frac{\norm{\sqrt{2\nu}\E_\Gamma(\u)}^2}{\norm{\u}^2}=\frac{\norm{\sqrt{2\nu}\E_\Gamma(\u)}^2}{\norm{\u_K}^2+\norm{\u_{NK}}^2},\quad \forall \u\in \H^1_\sigma(\Gamma),\quad \u\not\equiv0,
\end{align}
where we set $\u_K:=\Pk\u$ and $\u_{NK}:=(\I-\Pk)\u$. Note that this is not a Dirichlet quotient, since $\Lambda=0$ for any $\u\in \KK$, and in general $\KK$ is nontrivial.

Let us consider the two solutions $\u_i$, $i=1,2$ given in the
statement, and introduce $\u=\u_1-\u_2$,   and  $p=p_1-p_2$, so that they solve
\begin{align}
	\partial_t\u +(\u\cdot\nabla_\Gamma)\u_1+(\u_2\cdot\nabla_\Gamma)\u-2\P_\Gamma\divg(\nu\E_\Gamma(\u))+\nabla_\Gamma p=\f(\cdot,\u_1)-\f(\cdot,\u_2),\label{a0}
\end{align}
where we recall that we assumed $\nu\in W^{1,\infty}(\Gamma)$ such
that $\nu\geq \nu_*>0$. Applying the  Helmholtz  projector $\P_0$, we can also write 
\begin{align}
	\partial_t\u +\P_0(\u\cdot\nabla_\Gamma)\u_1+\P_0(\u_2\cdot\nabla_\Gamma)\u+\A_S\u=\P_0\f(\cdot,\u_1)-\P_0\f(\cdot,\u_2),\label{a1}
\end{align}
where we set $\A_S\u:=-2\P_0\P_\Gamma\divg(\nu\E_\Gamma(\u))$ as the (modified) Stokes operator. We observe that, after an integration by parts, it holds 
\begin{align}
	\int_\Gamma{\A_S(\v)}\cdot \v=\norm{\sqrt{2\nu}\E_\Gamma(\v)}^2,\quad \forall \v\in \H^1_\sigma(\Gamma).
	\label{intbyparts}
\end{align}
Now, assume by contradiction that there exists $t_*\in [0,T_*)$ such
that $\u_1(t_*)\not=\u_2(t_*)$, i.e., there exists a set of positive
Lebesgue measure on which the two solutions are different. Then
$\norm{\u(t_*)}>0$, and, recalling that, as   $\u_1,\u_2\in
C([0,\infty);\H^1_\sigma(\Gamma))$, this entails that there exists an
interval $I_\delta:=[t_*,t_*+\delta)$, $\delta>0$, such that
$\norm{\u(t)}>0$ for any $t\in I_\delta$. Assuming $I_\delta$ to be
the largest  of such intervals  with this property, it must be $\norm{\u(t_*+\delta)}=0$. We will now consider $t\in I_\delta$ and show that we obtain a contradiction.

Let us also introduce the quantity
\begin{align*}
	L(t):=-\frac12 \ln(\norm{\u(t)}^2),\quad t\in I_\delta,
\end{align*}
which is clearly well defined on $I_\delta$, and observe that we have the following identities, which are obtained by multiplying (and integrating) \eqref{a0} by $\u$ and \eqref{a1} by $\A_S\u$, respectively:
\begin{align}
	\label{e1}
	&\frac12\frac{d}{dt}\norm{\u}^2=-\int_\Gamma (\u\cdot\nabla_\Gamma )\u_1\cdot \u -\norm{\sqrt{2\nu}\E_\Gamma(\u)}^2+\int_\Gamma (\f(\cdot,\u_1)-\f(\cdot,\u_2))\cdot \u,\\&\nonumber
	 \frac12\frac d{dt}\norm{\sqrt{2\nu}\E_\Gamma(\u)}^2\\&	\label{e2}=-\int_\Gamma (\u\cdot\nabla_\Gamma )\u_1\cdot \A_S\u-\int_\Gamma (\u_2\cdot\nabla_\Gamma )\u\cdot \A_S\u-\norm{\A_S\u}^2+\int_\Gamma (\f(\cdot,\u_1)-\f(\cdot,\u_2))\cdot \A_S\u.
\end{align}
Now, let us consider the time derivative of $L(t)$: we have
\begin{align*}
	\frac d{dt}L(t)=-\frac{\frac12\frac d{dt}\norm{\u}^2}{\norm{\u}^2}=\frac{\int_\Gamma (\u\cdot\nabla_\Gamma )\u_1\cdot \u}{\norm{\u}^2}+\Lambda(t)-\frac{\int_\Gamma (\f(\cdot,\u_1)-\f(\cdot,\u_2))\cdot \u}{\norm{\u}^2},
\end{align*}  
and  are left with estimating  the first and the last terms
in the right-hand side  in terms of   
$\Lambda$. First observe that, by Korn's inequality \eqref{Korn1},  we have
\begin{align}
	\norm{\u}_{\H^1(\Gamma)}\leq C(\norm{\u}+\norm{\E_\Gamma(\u)})=C(\norm{\u_K}+\norm{\u_{NK}}+\norm{\E_\Gamma(\u)})\leq C(\norm{\u_K}+\norm{\E_\Gamma(\u)}).
	\label{essential}
\end{align}
Then, by H\"{o}lder's and Gagliardo--Nirenberg's inequalities, we have, recalling \eqref{essential} and $\nu\geq \nu_*>0$,
\begin{align*}
&\frac{\int_\Gamma (\u\cdot\nabla_\Gamma )\u_1\cdot \u}{\norm{\u}^2}\leq \frac{1}{\norm{\u}^2}\norm{\u}_{\L^4(\Gamma)}^2\norm{\u_1}_{\H^1(\Gamma)}\\&\leq C\frac{\norm{\u}_{\H^1(\Gamma)}\sqrt{\norm{\u_K}^2+\norm{\u_{NK}}^2}}{\norm{\u_K}^2+\norm{\u_{NK}}^2}\norm{\u_1}_{\H^1(\Gamma)}\leq C\frac{\norm{\u_K}+\norm{\E_\Gamma(\u)}}{\sqrt{\norm{\u_K}^2+\norm{\u_{NK}}^2}}\norm{\u_1}_{\H^1(\Gamma)}\\&
\leq C\norm{\u_1}_{\H^1(\Gamma)}^2+\frac{\norm{\u_K}^2}{{\norm{\u_K}^2+\norm{\u_{NK}}^2}}+\frac{1}{2\nu_*}\frac{\norm{\sqrt{2\nu}\E_\Gamma(\u)}^2}{{\norm{\u_K}^2+\norm{\u_{NK}}^2}}\\&
\leq C\norm{\u_1}_{\H^1(\Gamma)}^2+1+\frac1{2\nu_*}\Lambda(t)\leq C(T)+\frac1{2\nu_*}\Lambda(t),
\end{align*}
where $C(T)$ depends on some time horizon $T>T^*$, since $\u_1\in L^\infty(0,T;\H^1_\sigma(\Gamma))$ for any $T>0$.

Concerning the terms related to the forcing $\f$, recalling assumption
\eqref{C2} and by Gagliardo--Nirenberg's and Young's inequalities,
 we have 
\begin{align*}
	&\norma{\frac{\int_\Gamma (\f(\cdot,\u_1)-\f(\cdot,\u_2))\cdot \u}{\norm{\u}^2}}\leq \frac{C}{\norm{\u}^2}\norm{\u}^2=C.
\end{align*}
To sum up, we have obtained that 
\begin{align}
	\frac{d}{dt}L(t)\leq C(T)+C\Lambda(t),\quad \forall t\in I_\delta.\label{L1}
\end{align}
 Let us now establish some  control on $\Lambda(t)$. 
To this aim, similar to \cite{Foias}, recalling \eqref{intbyparts}, we
observe that 
\begin{align*}
	&\norm{\A_S\u-\Lambda(t)\u}^2=\norm{\A_S\u}^2+\frac{\norm{\sqrt{\nu}\E_\Gamma(\u)}^4}{\norm{\u}^2}-2\frac{\norm{\sqrt{2\nu}\E_\Gamma(\u)}^2}{\norm{\u}^2}\int_\Gamma \A_S\u\cdot\u\\&
	=\norm{\A_S\u}^2-\frac{\norm{\sqrt{2\nu}\E_\Gamma(\u)}^4}{\norm{\u}^2}=\norm{\A_S\u}^2-\Lambda(t)\norm{\sqrt{2\nu}\E_\Gamma(\u)}^2.
\end{align*}
Then, we can compute the time derivative of $\Lambda$.  By 
H\"{o}lder's inequality and exploiting \eqref{e1}--\eqref{e2} we get
\begin{align*}
&	\frac d{dt}\Lambda(t)=\frac{\frac{d}{dt}\norm{\sqrt{2\nu}\E_\Gamma(\u)}^2}{\norm{\u}^2}-\Lambda(t)\frac{\frac d {dt}\norm{\u}^2}{\norm{\u}^2}\\&
=\frac{2}{\norm{\u}^2}\left(-\int_\Gamma (\u\cdot\nabla_\Gamma )\u_1\cdot \A_S\u-\int_\Gamma (\u_2\cdot\nabla_\Gamma )\u\cdot \A_S\u-\norm{\A_S\u}^2+\int_\Gamma (\f(\cdot,\u_1)-\f(\cdot,\u_2))\cdot \A_S\u\right)\\&
-\frac{2\Lambda(t)}{\norm{\u}^2}\left(-\int_\Gamma (\u\cdot\nabla_\Gamma )\u_1\cdot \u-\int_\Gamma (\u_2\cdot\nabla_\Gamma )\u\cdot \u -\norm{\sqrt{2\nu}\E_\Gamma(\u)}^2+\int_\Gamma (\f(\cdot,\u_1)-\f(\cdot,\u_2))\cdot \u\right)\\&
=\frac{2}{\norm{\u}^2}\left((\Lambda(t)\norm{\sqrt{2\nu}\E_\Gamma(\u)}^2-\norm{\A_S\u}^2)+\int_\Gamma (\u\cdot\nabla_\Gamma )\u_1\cdot (\Lambda(t)\u-\A_S\u)\right.\\& \left.+\int_\Gamma (\u_2\cdot\nabla_\Gamma )\u\cdot (\Lambda(t)\u-\A_S\u)  +\int_\Gamma (\f(\cdot,\u_1)-\f(\cdot,\u_2))\cdot (\A_S\u-\Lambda(t)\u) \right)\\&
\leq \frac{2}{\norm{\u}^2}\left(-\norm{\Lambda(t)\u-\A_S\u}^2+\norm{\u}_{\L^4(\Gamma)}\norm{\u_1}_{\W^{1,4}(\Gamma)}\norm{\Lambda(t)\u-\A_S\u}\right.\\&\left.+\norm{\u_2}_{\L^\infty(\Gamma)}\norm{\u}_{\H^1(\Gamma)}\norm{\Lambda(t)\u-\A_S\u}+\norm{\f(\cdot,\u_1)-\f(\cdot,\u_2)}\norm{\Lambda(t)\u-\A_S\u} \right)\\&
=\frac{2}{\norm{\u}^2}\left(-\norm{\Lambda(t)\u-\A_S\u}^2+\mathcal I_0+\mathcal I_1+\mathcal I_2\right).
\end{align*}
To estimate $\mathcal I_0$, by Young's and Gagliardo--Nirenberg's
inequalities, together with \eqref{essential},  we have 
\begin{align*}
	\mathcal I_0&\leq C\norm{\u}\norm{\u}_{\H^1(\Gamma)}\norm{\u_1}_{\H^2(\Gamma)}\norm{\u_1}_{\H^1(\Gamma)}+\frac14\norm{\Lambda(t)\u-\A_S\u}^2\\&
	\leq C(T)\norm{\u}^2\norm{\u_1}_{\H^2(\Gamma)}^2+C(\norm{\u_K}^2+\norm{\E_\Gamma(\u)}^2)+\frac14\norm{\Lambda(t)\u-\A_S\u}^2\\&\leq 
	C(T)\norm{\u}^2(1+\norm{\u_1}_{\H^2(\Gamma)}^2)+C\norm{\sqrt{2\nu}\E_\Gamma(\u)}^2)+\frac14\norm{\Lambda(t)\u-\A_S\u}^2,
\end{align*}
where we used $\norm{\u_K}\leq \norm{\u}$ as well as $\norm{\u_1}_{L^\infty(0,T;\H^1_\sigma(\Gamma))}\leq C(T)$, choosing $T>T^*$.

 We estimate $\mathcal I_1$  analogously,  by using 
Agmon's, Young's,  Korn's  inequalities, and
\eqref{essential},  as 
\begin{align*}
	\mathcal I_1&\leq C\norm{\u_2}_{\H^1(\Gamma)}^\frac12\norm{\u_1}_{\H^2(\Gamma)}^\frac12\sqrt{\norm{\u_K}^2+\norm{\E_\Gamma(\u)}^2}\norm{\Lambda(t)\u-\A_S\u}\\&
	\leq C(T)\norm{\u_2}_{\H^2(\Gamma)}(\norm{\sqrt{2\nu}\E_\Gamma(\u)}^2+\norm{\u}^2)+\frac14\norm{\Lambda(t)\u-\A_S\u}^2,
\end{align*}
where again $\norm{\u_2}_{L^\infty(0,T;\H^1_\sigma(\Gamma))}\leq C(T)$, choosing $T>T^*$.

In conclusion, recalling \eqref{C2}, we have, by Gagliardo--Nirenberg's and Young's inequalities, together with \eqref{essential},
\begin{align*}
	\mathcal I_2&\leq C\norm{\u}\norm{\Lambda(t)\u-\A_S\u}\\&
	\leq C\norm{\u}^2+\frac14\norm{\Lambda(t)\u-\A_S\u}^2.
\end{align*}

To sum up, we have obtained 

\begin{align*}
	\frac d{dt}\Lambda(t)+ \frac12
  \frac{\norm{\Lambda(t)\u-\A_S\u}^2}{\| u \|^2}\leq C(T)(1+\norm{\u_2}_{\H^2(\Gamma)})\Lambda(t)+C(T)(1+\norm{\u_1}_{\H^2(\Gamma)}^2),\quad \forall t\in I_\delta.
\end{align*}
Recalling that, for $T>T^*$, we have $\u_i\in L^2(0,T;\H^2(\Gamma))$, $i=1,2$, we can apply Gronwall's Lemma on the interval $I_\delta$ and obtain that
\begin{align*}
	\sup_{t\in[t_*,t_*+\delta)}\Lambda(t)\leq \Lambda(t_0)e^{C(T)}+C(T),
\end{align*}
entailing that $\Lambda$ is bounded on $I_\delta$. From this,
recalling \eqref{L1}, we deduce by integrating over $I_\delta$ that
also $L(t)$ is bounded on $[t_*,t_*+\delta)$. In particular, we can
deduce that also $L(t_*+\delta)$ is bounded, contradicting the fact
that $\norm{\u(t_*+\delta)}=0$.  This concludes the contradiction
argument and   $t_*\in[0,T^*)$ with the property that $
\norm{\u(t_*)}>0$  does not exist. 

\section{ Preliminary results on long-time behavior:  Proofs
  of Section \ref{longtt}}\label{sec:prel}

\subsection{Proof of Proposition \ref{dec}}
	In the case $\mathcal{K}=\{\boldsymbol 0\}$ the proof is trivial. We thus consider the case of nontrivial Killing field space $\mathcal{K}$ of dimension $n\geq1$.
	We see by summing up \eqref{t11} and \eqref{K1} that it holds,
        by uniqueness,  $\overline{\u}=\u_{K}+\u_{NK}$  (recall
        that $\E_\Gamma(\u_K)=\boldsymbol 0$ by definition of Killing
        vector field), and thus $\u_{NK}$ {satisfying \eqref{t11}} can be obtained as the
        difference between $\overline{\u}$ and $\u_K$, as soon as we
        prove that $\u_K$ exists. In order to prove the proposition we
        then need first to show that $\u_K$ exists and is uniquely
        defined. To see this, notice that, from its definition,
        $\u_K\in \mathcal K$ for  all  $t\geq0$.  Thus,
          considering the basis $\{\v_i\}$ of $\mathcal{K}$, it can be written as $\u_K(t)=\sum_{i=1}^n\alpha_i(t)\v_i$, where $\alpha_i$ can be assumed to be in $C^1([0,\infty))$. Let us multiply \eqref{K1} by $\v_i$ and integrate over $\Gamma$. Since the basis $\{\v_j\}$ is orthonormal, this gives a system of ODEs in $\alpha_i$ of the form
	\begin{align}
		\partial_t\alpha_i=\int_\Gamma\f_K\left(x,\overline{\u}\right)\cdot \v_i, \quad i=1,\ldots,n.
		\label{ode}
	\end{align}
	%
	The right-hand side of \eqref{ode} is continuous in time. Indeed, thanks to the regularity of $\overline{\u}$, $\overline{\u}(s)\to\overline{\u}(t)$ as $s\to t$, $t\geq0$, in $\Lts$, and thus, by the continuity of $\Pk$,
	\begin{align*}
		&\norma{\int_\Gamma(\f_K(\cdot,\overline{\u}(t))-\f_K(\cdot,\overline{\u}(s)))\cdot \v_i}\\&\leq \norm{\f(\cdot,\overline{\u}(t))-\f(\cdot,\overline{\u}(s))}\norm{\v_i}\\&
		\leq 
		C\norm{\overline{\u}(t)-\overline{\u}(s)}\to 0\quad\text{as }s\to t.
	\end{align*}
	Therefore, by standard theory there exists $t_*>0$ and a
        unique vector $(\alpha_i)_i\in C^1([0,t_*];\mathbb R^n)$ such
        that \eqref{ode} is satisfied. This corresponds to a unique
        solution $\u_K$ to \eqref{K1}--\eqref{K2}. Moreover, by
        multiplying \eqref{K1} by $\u_K$, applying  the 
        Cauchy-Schwartz inequality, and recalling
        \eqref{C1}--\eqref{C2},  we infer  that
	\begin{align*}
		\frac 12\partial_t\sum_{i=1}^n\norma{\alpha_i}^2&\leq \norm{\f(\cdot, \overline{\u})-\f(\cdot,\boldsymbol 0)}\sqrt{\sum_{i=1}^n\norma{\alpha_i}^2}+\norm{\f(\cdot,\boldsymbol 0)}\sqrt{\sum_{i=1}^n\norma{\alpha_i}^2}\\&\leq C\sum_{i=1}^n\norma{\alpha_i}^2+C(1+\norm{\overline\u}^2),
	\end{align*} 
	entailing that $t_*=+\infty$, i.e., the solution $\u_K\in
        \mathcal K$ is indeed global. It is now trivial to deduce that
        $\u_{NK}=\overline{\u}-\u_{K}$ {satisfies \eqref{t11}},  is uniquely determined, and  belongs to $(\I-\Pk)\Lts$, by
        simply multiplying \eqref{t11} by $\z\in \mathcal K$,
        integrating over $\Gamma$, and recalling (see again
        \cite{Simonett}) that $\int_\Gamma
        (\overline{\u}\cdot\nabla_\Gamma)\overline{\u}\cdot
        \z=0$. Since also $\divg\z=0$, we deduce that  one has 
	$$
	\partial_t\int_{\Gamma}\u_{NK}\cdot \z=0,\quad \forall \z\in \mathcal K.
	$$
	 On the other hand,   we also have that $\int_\Gamma
        \u_{NK}(0)\cdot \z=\int_\Gamma (\boldsymbol
        I-\P_{\mathcal{K}})\u_0\cdot \z=0$ for any $\z\in \mathcal
        K$.  This entails that  $\u_{NK}(t)\in (\I-\Pk)\Lts$ for any $t\geq0$. Since $\L^2_\sigma(\Gamma)=\mathcal K\oplus (\I-\Pk)\Lts$,
	we have just  proved  the assertion of the proposition, by the uniqueness of this decomposition, i.e., $\u_K=\Pk\overline{\u}$ and $\u_{NK}=(\I-\Pk)\overline{\u}$.
\subsection{Proof of Proposition \ref{lemmadec}}
In order to  prove  \eqref{exp}, it is enough to multiply equation \eqref{K1} by $\u_K$ to obtain, by the assumptions on $\f_K$ and Proposition \ref{dec},
	$$
	\frac12\frac{d}{dt}\norm{\u_K}^2=\int_\Gamma
        \f_K(x,\overline{\u})\cdot \u_K=\int_\Gamma
        \f_K(x,\overline{\u})\cdot \overline{\u}\leq  C_3 
        \norm{\u_K}^2+ C_4 ,
	$$ 
	where we recall that $\overline{\u}=S(t)\u_0$ is the entire
        unique global solution corresponding to $\u_0$. This entails
        the result by Gronwall's Lemma.
        
	For case (i), we multiply \eqref{K1} by $\f_K$ and integrate
        over $\Gamma$  getting 
	\begin{align*}
		\frac{d}{dt}\int_\Gamma \u_K\cdot \f_K=\norm{\f_K}^2,
	\end{align*}
	so that 
	$$
	\int_\Gamma \u_K(t)\cdot \f_K=\int_\Gamma \u_K(0)\cdot \f_K+t\norm{\f_K}^2
	$$
         which is \eqref{fd}. 
	Moreover, if we consider the subspace $\mathcal{K}\cap \f_K^{\perp_{\L^2}}$, it holds
	\begin{align*}
		\frac{d}{dt}\int_\Gamma \u_K\cdot \v=0,\quad \forall \v\in \mathcal{K}\cap \f_K^{\perp_{\L^2}},
	\end{align*}
	entailing \eqref{uk1}. Summing up \eqref{fd} squared with
        \eqref{uk1}  one obtains  \eqref{uk2}.
	
 In order to check case  (ii),  one multiplies 
\eqref{K1} with $\u_K$  and integrates on $\Gamma$ getting  
	$$
	\frac12 \frac{d}{dt}\norm{\u_K}^2=\int_\Gamma \f_K(x,\overline{\u})\cdot \u_K=\int_\Gamma \f_K(x,\overline{\u})\cdot \overline{\u}\leq 0.
	$$
	
	Analogously, for the case (iii) we have
	$$
	\frac12 \frac{d}{dt}\norm{\u_K}^2=\int_\Gamma \f_K(x,\overline{\u})\cdot \u_K=\int_\Gamma \f_K(x,\overline{\u})\cdot \overline{\u}\geq 0.
	$$
 The assertion follows by integrating in time.  

\subsection{Proof of Proposition \ref{pp}}
 Multiply  \eqref{t11} by $\u_{NK}$ and integrate over $\Gamma$. Observe that (see \cite{Simonett})
	$$
	\int_\Gamma (\overline{\u}\cdot\nabla_\Gamma)\overline{\u}\cdot \u_{NK}=\int_\Gamma (\overline{\u}\cdot\nabla_\Gamma)\overline{\u}\cdot \overline{\u}=0,
	$$
	since $\int_\Gamma
        (\overline{\u}\cdot\nabla_\Gamma)\overline{\u}\cdot {\u}_K=0$
         from the divergence-free property of  
        $\overline{\u}$. Therefore, from \eqref{extra} and Proposition
        \ref{dec}, recalling that $\nu\geq \nu_*>0$  we get that 
	\begin{align}
		&\nonumber \frac12\ddt\norm{\u_{NK}}^2+2\nu_*\norm{\E_\Gamma(\u_{NK})}^2\\&
		\leq { C_5 }\norm{\u_{NK}}^2+ C_6 \norm{\u_{NK}}\nonumber\\&
		\leq 
		2  C_5 \norm{\u_{NK}}^2+C.\label{energy}
	\end{align}
         From  Korn's inequality \eqref{Korn3} we deduce from \eqref{energy} that 
	\begin{align}
		& \frac12\ddt\norm{\u_{NK}}^2+\left(\frac{2\nu_*}{C_P^2}-2  C_5 \right)\norm{\u_{NK}}^2\leq C,\label{energy2}
	\end{align}
	which gives the desired result \eqref{gr} by Gronwall's Lemma,
        since by assumption $\zeta:=\frac{2\nu_*}{C_P^2}-2 C_5>0$.
	
Under the additional assumption \eqref{nega} on $\f_K$, we have 
	$$
	\norm{\u_K (t) }\leq \norm{\u_K(0)},\quad \forall t\geq0,
	$$
	so that,  working now with  \eqref{extra2} and performing similar energy estimates as above we obtain
	\begin{align}
		&\nonumber \frac12\ddt\norm{\u_{NK}}^2+2\nu_*\norm{\E_\Gamma(\u_{NK})}^2\\&
		\leq { C_5}\norm{\u_{NK}}^2+ C_6\norm{\u_{NK}}+ C_6\norm{\u_K}^2\nonumber\\&\nonumber \leq  C_5\norm{\u_{NK}}^2+ C_6\norm{\u_{NK}}+ C_6\norm{\u_K(0)}^2
		\\&
		\leq 
		2  C_5\norm{\u_{NK}}^2+C(1+\norm{\u_K(0)}^2).\label{energy3}
	\end{align}
 Again   by Korn's inequality  we have 
	\begin{align}
		& \frac12\ddt\norm{\u_{NK}}^2+\left(\frac{2\nu_*}{C_P^2}-2  C_5\right)\norm{\u_{NK}}^2\leq C(1+\norm{\u_K(0)}^2),\label{energy2b}
	\end{align}
	allowing to deduce \eqref{gr2} by Gronwall's Lemma.

 \subsection{Proof of Theorem \ref{insta}} 
	The instantaneous regularization of (global) weak solutions is already pointed out in Theorem \ref{thm1}. Let us consider a weak solution $\u$ departing from an initial datum $\u_0\in \L^2_\sigma(\Gamma)$. From any positive time onward, this solution is also strong. We now aim at finding uniform estimates for all positive times.  First, from Proposition \ref{lemmadec} and \eqref{gr2} we have
	\begin{align*}
		\norm{\u  (t) }^2&=\norm{\u_{NK} (t) }^2+\norm{\u_{K} (t) }^2\\&\leq
		e^{-\zeta t}\norm{\u_{NK}(0)}^2+\omega(1+\norm{\u_K(0)}^2)+\norm{\u_K(0)}^2 \\&
		\leq C\norm{\u_0}^2,\quad \forall t\geq 0,
	\end{align*} 
	entailing that 
	\begin{align}
		\norm{\u}_{L^\infty(0,\infty;\L^2_\sigma(\Gamma))}\leq C,\label{unif}
	\end{align}
	where $C$ depends on the $\L^2$-norm of $\u_0$. Then, repeating the same estimates leading to \eqref{energy}, we get
	\begin{align}
		2\nu_*\norm{\E_\Gamma(\u_{NK})}^2
		\leq 
		2 C_5  \norm{\u_{NK}}^2+C,\label{energybb}
	\end{align}
	so that, integrating over $[t,t+1]$, $t\geq0$,  using  Korn's
        inequality \eqref{Korn1} we obtain
	$$
	\norm{\u_{NK}}_{L^2(t,t+1;\H^1_\sigma(\Gamma))}\leq C,\quad \forall t\geq 0.
	$$
	Now recall that,  again  by Korn's inequality \eqref{Korn2}, it holds
	\begin{align*}
		\norm{\u_K}_{\H^1(\Gamma)}\leq C\norm{\u_K}, 
	\end{align*}
	so that, together with \eqref{unif}, we immediately end up with 
	\begin{align}
		\norm{\u}_{L^2(t,t+1;\H^1_\sigma(\Gamma))}\leq C,\quad \forall t\geq 0,\label{H1}
	\end{align}
	where $C$ only depends on $\norm{\u_0}$, $\f$, $\Gamma$, and
        the parameters of the problem.

        We can now deal with higher-order estimates. Let us fix an
        arbitrary $\tau>0$. The solution $\u$ is strong on
        $[\tau,\infty)$, so that, by repeating the same computations
        leading to \eqref{final}  and  exploiting \eqref{unif} we get 
	\begin{align}
	\label{final2}
	\frac{d}{dt}\int_\Gamma \nu \norma{\E_\Gamma(\u)}^2+\gamma \norm{\u}_{\H^2(\Gamma)}^2+\frac14\norm{\partial_t\u}^2\leq C\norm{\E_\Gamma(\u)}^2\left(\int_\Gamma \nu\norma{\E_\Gamma(\u)}^2\right)+C,\quad \forall t\geq \tau.
        \end{align}
         Indeed, the constant in estimate \eqref{final} depends on
        the final time. Nonetheless, given the information that we now
        have, the argument towards \eqref{final} could be refined in
        order to prove that no dependence on $T$ is actually needed in
        that constant. 
Recalling \eqref{H1} and $\nu\geq \nu_*>0$, we can apply the uniform Gronwall's Lemma (\cite{Temam}) to infer
\begin{align*}
	\norm{\E_\Gamma(\u)}_{L^\infty(\tau,\infty;\L^2(\Gamma))}+\norm{\u}_{L^2(t,t+1;\H^2(\Gamma))}+\norm{\partial_t\u}_{L^2(t,t+1,\L_\sigma^2(\Gamma))}\leq C(\tau),\quad \forall t\geq \tau,
\end{align*}
where $C$ depends on $\tau$, $\norm{\u_0}$, $\f$, $\Gamma$, and the parameters of the problem.
Again by Korn's inequality, this also implies, thanks to \eqref{unif},
\begin{align*}
		\norm{\u}_{L^\infty(\tau,\infty;\H^1_\sigma(\Gamma))}\leq C,
\end{align*}
entailing also \eqref{regg1} by standard embeddings.

In conclusion, if we assume $\u_0\in \H^1_\sigma(\Gamma)$, by Theorem \ref{thm1} there exists a unique strong solution $\u$ departing from this initial datum. Then, having fixed $\tau=1$, we know that estimates \eqref{regg1}--\eqref{regg2} hold on $[1,\infty)$. Moreover, from Theorem \ref{thm1} we infer that 
\begin{align*}
		\norm{\u}_{L^\infty((0,1);\H^1_\sigma(\Gamma))}+\norm{\u}_{L^2((0,1);\H^2(\Gamma))}+\norm{\partial_t\u}_{L^2((0,1),\L_\sigma^2(\Gamma))}\leq C(\norm{\u_0}_{\H^1_\sigma(\Gamma)}),
\end{align*}
and thus, putting together the estimates, we deduce that
\eqref{regg1}--\eqref{regg2} also hold on $[0,\infty)$, as long as we
recall that the constants appearing also depend  on the $\H^1$-norm of
$\u_0$.  This concludes the proof. 
 \section{ Attractors for bounded trajectories:  Proofs of Section \ref{sec:attr}}\label{bdda}
  \subsection{Proof of Theorem \ref{th1}}
 We begin the proof by fixing $r\geq0$ and considering the complete
 metric space $\mathbb B_r$ defined in \eqref{Ar}.  Recall that 
 $$
 \norm{A}_\Lts:=\sup_{\u\in A}\norm{\u},\quad \forall A\subset \Lts.
 $$
 Thanks to Proposition \ref{pp}, in the notation of the (restricted) dynamical system $(S(t),\BBB_r)$, we deduce that 
 
 \begin{align*}
 \norm{(\I-\P_\mathcal K)S(t)B}_{\Lts}\leq \sqrt{e^{-\zeta t}\norm{(\I-\P_\mathcal K)B}_{\Lts}^2+\omega(1+r^2)},\quad \forall t\geq 0,	
 \end{align*}
 for any bounded $B\subset \BBB_r$, where the constants $\zeta,\omega$
 do not depend on the size of $B$, i.e.,  on 
 $\norm{B}_{\Lts}$. Moreover, due to the assumptions on $ \f_K$,
  exploiting Proposition \ref{lemmadec} point (ii) we also have
 \begin{align*}
 	\norm{\Pk S(t)B}_\Lts\leq \norm{\Pk B}_\Lts\leq  r,\quad \forall t\geq 0,
 \end{align*}
 for \textit{any} set $B\subset \BBB_r$. This shows in particular that the map $S(t)$  satisfies
 \begin{align*}
 	S(t):\ \BBB_r\to \BBB_r, \quad \forall t\geq 0.
 \end{align*} 
 Therefore, we can consider the dynamical system $(S(t),\BBB_r)$. 
 Next, we check that  
 \begin{align}
 	\mathcal B_0^r:= \left\{\u\in \BBB_r:\ \norm{(\I-\Pk)\u}\leq \sqrt{\frac12+\omega(1+r^2)}\right\}\subset \BBB_r,
 	\label{B0r}
 \end{align}  
 which is an absorbing set for $(S(t),\BBB_r)$. Indeed, from the estimates above we have that, given $B\subset \BBB_r$ bounded set, there exists $t_B:=\ln(2\norm{(I-\Pk)B}_\Lts)/\zeta$, depending only on the size of $B$ (i.e., $\norm{B}_{\Lts}$), so that
 $$
 S(t)B\subset \BB_0^r,\quad \forall t\geq t_B.
 $$
 Since also $\BB_0^r$ is bounded, there exists $t_{\BB_0^r}\geq0$, depending only on $\zeta,\omega,r$ through $\norm{(\I-\Pk)\BB_0^r}_\Lts$, such that
 \begin{align}
 S(t)\BB_0^r\subset \BB_0^r,\quad \forall t\geq t_{\BB_0^r}.
 \label{absorbing1}
 \end{align}
 We can now observe that, thanks to the instantaneous-regularization property \eqref{regg1} of any solution $\u$, for the fixed $\tau=t_{\BB_0^r}$ it holds (see Theorem \ref{insta}) 
 $$
 \sup_{t\geq t_{\BB_0^r}}\norm{S(t)\BB_0^r}_{\H^1(\Gamma)}\leq C(t_{\BB_0^r})\norm{\BB_0^r}_{\Lts}\leq C(t_{\BB_0^r},r,\omega,\zeta), 
 $$
 where the dependencies of $C$ come from the definition of $\BB_0$.
 Therefore, we can introduce the set 
 \begin{align}
 \BB_1^r:=\{\u\in \BB_0^r\cap \H^1(\Gamma):\ \norm{\u}_{\H^1(\Gamma)}\leq C(t_{\BB_0^r}, r,\omega,\zeta)\},
 \label{B1r}
 \end{align}
 so that 
 $$S(t)\BB_0^r\subset \BB_1^r,\quad \forall t\geq t_{\BB_0^r},$$
 entailing that $\BB_1^r\subset \BB_r$ is also a bounded absorbing set. Since this set is also compact, by a standard application of the theory for dissipative dynamical systems on the complete metric space $\BBB_r$ (see, e.g., \cite{Temam}) we infer that there exists the (unique) global attractor $\AA_r\subset \mathcal B_1^r\subset \H^1_\sigma(\Gamma)$ for the dynamical system $(S(t),\BBB_r)$, such that  
 \begin{itemize}
 	\item $\mathcal A_r$ is nonempty, compact and connected,
 	\item $\mathcal A_r$ is invariant, i.e., $S(t)\mathcal A_r=\mathcal A_r$, for any $t\geq0$,
 	\item $\mathcal A_r$ is attracting, i.e., for any $B\subset \mathbb B_r$ bounded it holds 
 	${\rm dist}(S(t)B,\mathcal A_r)\to 0$ as $t\to \infty$.
 \end{itemize}
 Moreover, it also holds
 \begin{align}
 	\mathcal A_r=\{\xi(0)\in \mathbb B_r:\ \xi \text{ is a bounded complete trajectory for $S(t)$ in }\mathbb B_r \}.
 	\label{bdd_traj2}
 \end{align}
 This result can be obtained for any $r\geq 0$, so that we have constructed a family $\{\AA_r\}_{r\geq 0}\subset \Lts$ of global attractors for the dynamical systems $(S(t),\BBB_r)$. As already noticed, the characterization property \eqref{bdd_traj2} allows to deduce, since $\BBB_{r_1}\subset \BBB_{r_2}$ for any $r_1\leq r_2$, that $\mathcal A_{r_1}\subset \mathcal{A}_{r_2}$ for any $r_1\leq r_2$, i.e., the family $\{\AA_r\}_{r\geq 0}\subset \Lts$ is increasing.
 
 Our aim is now to  link  this family to the attractor of the
 \textit{full} dynamical system $(S(t),\Lts)$, which is expected to be
 the set $\mathcal J$ (coinciding  in this case  with $\mathcal I$, as already observed). We thus define the following 
 $$
 \widetilde{\AA}:={\bigcup_{r\geq 0}\AA_r}={\bigcup_{ r \in \N}\AA_r},
 $$
 where the last identity is due to the fact that the family of global attractors is increasing with  $r\geq 0$.
 
 First, we show that $\widetilde{\AA}$ is attracting for $(S(t),\Lts)$. Indeed, let us consider a bounded set $B\subset \Lts$. There exists $r>0$ sufficiently large such that $B\subset \BBB_r$. From the attracting property of the corresponding set $\AA_r\subset \widetilde{\AA}$, it then holds 
 \begin{align*}
 	\text{\rm dist}_\Lts(S(t)B, \widetilde{\AA})\leq \text{\rm dist}_\Lts(S(t)B, \AA_r)\to 0,\quad \text{as }t\to \infty. 
 \end{align*}
 As a consequence, $\tA$ is attracting.  As 
 $\overline{\tA}^{\Lts}$ is  obviously  closed, we immediately deduce from Property B. of Lemma \ref{obv} that $\mathcal J\subset \overline{\tA}^{\Lts}$. 
 
 Now, by construction any bounded set $\AA_r$ is fully invariant, in the sense that $S(t)\AA_r=\AA_r$ for any $t\geq 0$. Then, by Lemma \ref{obv} property C., $\AA_r\subset \mathcal J$ for any $r\geq 0$, and thus 
 $$
\bigcup_{r\geq 0}\AA_r\subset \mathcal J.
 $$
 Clearly, this immediately entails that $\JJ$ is nonempty.
 In order to conclude the identification and show that $\JJ=\tA$, we thus need to further show that $\JJ\subset \tA$. Notice that we cannot in general show that $\JJ$ is closed, which would also entail, from what observed above, that $\overline{\tA}^{\Lts}=\tA$. We will see some specific cases in Theorem \ref{spe} in which this result is true.

 To show $\JJ\subset \tA$, let us  fix  $\u\in \JJ$.  By
  construction there exists a complete bounded trajectory $\xi$
 such that $\xi(0)=\u$. Since the trajectory is bounded, it clearly
 holds $\xi(\R)\subset \BBB_r$ for some $r\geq0$. From property
 \eqref{bdd_traj2} of the corresponding set $\AA_r$, this also means
 that $\xi(0)=\u\in \AA_r$, entailing that $\JJ\subset \tA$.  We
 then conclude that  
 $\JJ=\tA$.

  Having  
 identified the two sets, we can proceed  with  the
 proof. First, since $\AA_r\subset \H^1_\sigma(\Gamma)$ for any
 $r\geq0$, then also $\JJ\subset \H^1_\sigma(\Gamma)$ and thus, since
 $\H^1_\sigma(\Gamma)\hookrightarrow \Lts$ compactly, $\JJ$ must have
 empty interior in $\Lts$,  which is  property 2  of  Theorem \ref{th1}.

Moreover,  as  $\tA$  attracts and $\tA=\JJ$ we also have
Property 3  of Theorem \ref{th1}. Property 4 is then a consequence of Proposition \ref{obv}, whereas Property 5 is again a consequence of Proposition \ref{obv} point B.
 
  In order to conclude the proof of Theorem \ref{th1}, we need to show
  that each set $\AA_r$ is of finite fractal dimension.  We prove
  that  for any $r\geq0$ there exists an exponential attractor
  $\mathcal M_r$ for the dynamical system $(S(t),\BBB_r)$, whose
  properties are  the following 
  \begin{itemize}
 	\item $\mathcal M_r$ is compact and of finite fractal dimension $N_r$, possibly increasing with $r\geq0$.
 	\item $\mathcal M_r$ is positively invariant, i.e.,
          $S(t)\mathcal M_r\subset \mathcal M_r$,  for any $t\geq0$,
 	\item $\mathcal M_r$ is exponentially attracting, i.e., for any $B\subset \mathbb B_r$ bounded it holds 
 	$${\rm dist}(S(t)B,\mathcal M_r)\leq C(\norm{B}_\Lts)e^{-\gamma_r t},$$ where $C>0$ depends on the size of $B$ (i.e., $\norm{B}_{\Lts}$), and $\gamma_r>0$ is a universal constant depending only on $r$. 
 \end{itemize}
 Since $\mathcal M_r$ is closed and attracting, it holds $\AA_r\subset \mathcal M_r$ and thus $\AA_r$ is of finite fractal dimension as well.
 
 In order to prove the existence of exponential attractors, we need some preliminary lemmas. First, for any $r\geq0$, recalling the definition of $\BB_1^r$ in \eqref{B1r}, we know that there exists $t_1^r=t_1^r(r,\zeta,\omega)$ such that 
 $$
 S(t)\BB_1^r\subset \BB_1^r,\quad \forall t\geq t_1^r.
 $$ 
 Therefore, we can introduce the set
 \begin{align*}
 	\mathbb{C}_r:=\overline {\bigcup_{t\geq t_1^r}S(t)\mathcal{B}_1^r}^{\Lts}\subset \BB_1^r,
 \end{align*}
 which is compact, positively invariant, and absorbing. Let us then prove the following.
 \begin{lemm}
Under assumptions \eqref{C1}--\eqref{C2}, given $\u_{0,1},\u_{0,2}\in \Lts$, it holds 
 	\begin{align}
 		\norm{\Pk \left(S(t)  \u_{0,1} -S(t) \u_{0,2} \right)}_{\H^1(\Gamma)}	\leq C\norm{\Pk \left(S(t)  \u_{0,1} -S(t)  \u_{0,2}  \right)},\quad \forall t\geq 0,
 		\label{regKill}
 	\end{align}
 	where $C>0$ only depends on $\Gamma$.
 	
 \noindent Moreover, given $\u_{0,1},\u_{0,2}\in \mathbb C_r$, for
  any  $r\geq0$  and  $T> 0$ there exists $C=C(T,r,\zeta,\omega)>0$ such that
 	\begin{align}
 		\Vert {\boldsymbol \varepsilon}_\Gamma\left(S(t)\uu_{0,1}-S(t)\uu_{0,2} \right)\Vert^2\leq \frac{C}{t}\Vert \uu_{0,1}-\uu_{0,2} \Vert^2,\quad \forall t\in(0,T],
 		\label{expo2}
 	\end{align}
 \end{lemm}
 \begin{proof}
 Thanks to Theorem \ref{thm1}, since $\mathbb C_r\subset \H^1_\sigma(\Gamma)$, given $\u_{0,1},\u_{0,2}\in \Lts$ there exist (unique) strong solutions $\u_1,\u_2$ on $[0,\infty)$ to \eqref{t1}--\eqref{tangential} corresponding to these initial data. Their regularity is enough to perform rigorously the next computations, leading to \eqref{expo2}. In particular, let us notice that $\u:=\u_1-\u_2$ satisfies
\begin{align}
	&\label{t22}\partial_t\u+(\u_1\cdot \nabla_\Gamma)\u+(\u\cdot \nabla_\Gamma)\u_2-2\P_\Gamma\divg(\nu\boldsymbol\varepsilon_\Gamma(\u))+\nabla_\Gamma p=\f(\cdot,\u_1)-\f(\cdot,\u_2),\\&
	\divg\u=0,\\&
	\u(0)=\u_{0,1}-\u_{0,2},
	\label{tangential2}
\end{align}
where $p$ as a suitable zero-integral-mean pressure, corresponding to $p_1-p_2$.  
First, to show \eqref{regKill}, it is enough to observe that, by Korn's inequality \eqref{Korn2}, 
\begin{align}
	\norm{\u}_{\H^1(\Gamma)}\leq C\norm{\u},\quad \forall \u\in \mathcal K,
\end{align}
thanks to the fact that ${\boldsymbol \varepsilon}_\Gamma(\u)=\boldsymbol 0$ for any $\u\in \mathcal K$. Let us now multiply \eqref{t22} by $\partial_t\u$ and integrate over $\Gamma$. After an integration by parts, this gives
\begin{align}
&\label{**}\frac d{dt}\int_\Gamma\nu(x)\norma{\E_\Gamma(\u)}^2+\norm{\partial_t\u}^2\\
  &\quad =\int_\Gamma (\u_1\cdot \nabla_\Gamma)\u\cdot \partial_t\u+\int_\Gamma (\u\cdot \nabla_\Gamma)\u_2\cdot \partial_t\u+\int_\Gamma(\f(\cdot,\u_1)-\f(\cdot,\u_2))\cdot \partial_t\u.\nonumber
	\end{align} 
	Now, by H\"{o}lder's and Korn's inequalities  and by
         recalling the embedding $\H^2(\Gamma)\hookrightarrow \L^\infty(\gam)$, we have 
	\begin{align*}
		&\norma{\int_\Gamma (\u_1\cdot \nabla_\Gamma)\u\cdot \partial_t\u}\\&\leq \norm{\u_1}_{\L^\infty(\gam)}\norm{\nabla_\Gamma \u}\norm{\partial_t\u}\\&\leq C\norm{\u_1}_{\H^2(\Gamma)}(\norm{\u}+\norm{{\boldsymbol \varepsilon}_\Gamma(\u)})\norm{\partial_t\u}
		\\&
		\leq C\norm{\u_1}_{\H^2(\Gamma)}^2\norm{\E_\Gamma(\u)}^2+C\norm{\u_1}_{\H^2(\Gamma)}^2\norm{\u}^2+\frac14\norm{\partial_t\u}^2.
	\end{align*}
In a similar way, by Gagliardo--Nirenberg's and Korn's inequalities we have
\begin{align*}
	\norma{\int_\Gamma (\u\cdot \nabla_\Gamma)\u_2\cdot \partial_t\u}&\leq \norm{\u}_{\L^4(\Gamma)}\norm{\nabla_\Gamma\u_2}_{\L^4(\Gamma)}\norm{\partial_t\u}
	\\&
	\leq 
	\norm{\u}^\frac12\norm{\u}_{\H^1(\Gamma)}^\frac12\norm{\u_2}_{\W^{1,4}(\Gamma)}\norm{\partial_t\u}\\&
	\leq C\norm{\u}^\frac12(\norm{\E_\Gamma(\u)}^\frac12+\norm{\u}^\frac12)\norm{\u_2}_{\W^{1,4}(\Gamma)}\norm{\partial_t\u}\\&
	\leq C(1+\norm{\u_2}_{\W^{1,4}(\Gamma)}^4)\norm{\u}^2+C\norm{\E_\Gamma(\u)}^2+\frac14\norm{\partial_t\u}^2.
\end{align*}
 Eventually,  
by assumption \eqref{C2} on $\f$ we have 
	\begin{align*}
		\norma{\int_\Gamma(\f(\cdot,\u_1)-\f(\cdot,\u_2))\cdot
          \partial_t\u}&\leq
                         \norm{\f(\cdot,\u_1)-\f(\cdot,\u_2)}\norm{\partial_t\u}\\&\leq
          C\norm{\u}\norm{\partial_t\u}\\&\leq C
          \norm{\u}^2+ \frac14 \norm{\partial_t\u}^2.
		\end{align*}
	As a consequence, we can sum up all the estimates and,
        recalling  that  $\nu\geq \nu_*>0$  and 
        multiplying equation  \eqref{**}   by $t$  we get 
	\begin{align}
		\nonumber&\frac d{dt}t\int_\Gamma\nu(x)\norma{{\boldsymbol \varepsilon}_\Gamma(\u)}^2+\frac18t\norm{\partial_t\u}^2\\&\leq\nonumber Ct(1+\norm{\u_1}_{\H^2(\Gamma)}^2)\int_\Gamma\nu(x)\norma{{\boldsymbol \varepsilon}_\Gamma(\u)}^2\\&+ Ct(1+\norm{\u_1}_{\H^2(\Gamma)}^2+\norm{\u_2}_{\W^{1,4}(\Gamma)}^4)\norm{\u}^2+C\int_\Gamma\nu(x)\norma{{\boldsymbol \varepsilon}_\Gamma(\u)}^2.\label{ree2}
	\end{align} 
	Now, recalling \eqref{contdep2} we deduce
	\begin{align}
	\sup_{t\in[0,T]}\norm{\u(t)}^2+2\nu_*\int_0^T \norm{\E_\Gamma(\u(t))}^2dt\leq C(T,r,\zeta,\omega)\norm{\u_{0,1}-\u_{0,2}}^2, 
	\label{ff}
	\end{align}
	so that 
		\begin{align}
		\nonumber&\frac d{dt}t\int_\Gamma\nu(x)\norma{\E_\Gamma(\u)}^2+ \frac14  t\norm{\partial_t\u}^2\\&\leq\nonumber Ct(1+\norm{\u_1}_{\H^2(\Gamma)}^2)\int_\Gamma\nu(x)\norma{{\boldsymbol \varepsilon}_\Gamma(\u)}^2\\&+ Ct(1+\norm{\u_1}_{\H^2(\Gamma)}^2+\norm{\u_2}_{\W^{1,4}(\Gamma)}^4)\norm{\u_{0,1}-\u_{0,2}}^2+\int_\Gamma\nu(x)\norma{{\boldsymbol \varepsilon}_\Gamma(\u)}^2,\label{ree}
	\end{align} 
	and, applying Gronwall's Lemma, since, by Theorem \ref{thm1}, $\u_i$, $i=1,2$, are strong solutions, and thus, by interpolation, $\u_i\in L^2(0,T;\H^2(\Gamma))\cap \L^\infty(0,T;\H^1(\Gamma))\cap L^4(0,T;\W^{1,4}(\Gamma))$, we obtain
	\begin{align*}
		t\norm{\E_\Gamma(\u)(t)}^2&\leq C(T,r,\zeta,\omega)\norm{\u_{0,1}-\u_{0,2}}^2+C\int_0^T\norm{\E_\Gamma(\u(t))}^2dt\\&\leq C(T,r,\zeta,\omega)\norm{\u_{0,1}-\u_{0,2}}^2,\quad \forall t\in(0,T],
	\end{align*}
	where we  also used  \eqref{ff}.  This concludes
         the proof of the lemma.
 \end{proof}
 Moving from  the continuous-dependence estimate \eqref{contdep2} and \eqref{regKill}--\eqref{expo2}, we can infer the following smoothing estimate
\begin{align}
	\norm{S(t)\u_{0,1}-S(t)\u_{0,2}}_{\H^1(\Gamma)}\leq \frac{C(T,r,\zeta,\omega)}{\sqrt t}\norm{\u_{0,1}-\u_{0,2}},\quad \forall t\in(0,T].\label{smoothing}
\end{align}
 We can now continue the proof  of Theorem \ref{thm:323}, following, for instance, \cite{Zelik}. Observe that, thanks to the uniform regularity of strong solutions given in Theorem \ref{insta}, it holds
 \begin{align*}
 	\norm{\partial_t S(t)\u_0}_{L^2(t,t+1;\Lts)}\leq
   C(r,\omega,\zeta,\tau),\quad \forall \u_0\in \mathbb C_r,\quad
   \forall t\geq  \tau .
 \end{align*}
Setting $ \u  (t)=S(t)  \u_0  $, with $ \u_0  \in\mathbb{C}_r$, we have, for any given $T>0$,
 \begin{align}
 	\label{lipt}
 	\Vert \u(t)-\u(s)\Vert\leq \int_s^t\Vert \partial_t\u(\tau)\Vert d\tau\leq \vert t-s\vert^{1/2}\left(\int_s^t\Vert \partial_t\u(\tau)\Vert^2 d\tau\right)^{1/2}\leq C(T,r,\omega,\zeta)\vert t-s\vert^{1/2}
 \end{align}
 for any $s,t\in[0,T]$, i.e., $t \mapsto S(t) u_0 $ is
 $(1/2)$-H\"{o}lder continuous in $[0,T]$, with $C$ depending only on
 $T,r,\omega,\zeta$. This, together with the continuous-dependence
 estimate \eqref{contdep2} allows to deduce that $S$ is $(1/2)$-H\"{o}lder continuous in $[0,T]\times \mathbb C_r$, for any
 $T>0$. Let us now fix $t_*>0$. Thanks to the smoothing property
 \eqref{smoothing} valid at $t=t_*>0$, the discrete dynamical system
 generated by the iterations of $(S(t_*),\mathbb C_r)$ possesses an
 exponential attractor $\mathcal{M}^*_r\subset \mathbb C_r$ (see,
 e.g., \cite[Thm.(2)7]{Zelik}). Moreover \eqref{contdep2} and
 \eqref{lipt} entail  that 
 $$
 S:[0,t_*]\times \mathbb{C}_r\to \mathbb{C}_r, \quad S(t, \u_0 ):=S(t)  \u_0 ,
 $$
 is H\"{o}lder continuous, when $\mathbb{C}_r$ is endowed with the $\Lts$ topology. Therefore, we can define
 $$
 \mathcal{M}_r:=\bigcup_{t\in[0,t_*]}S(t)\mathcal{M}_r^*\subset \mathbb{C}_r,
 $$
 and, following  again  \cite{Zelik}, show that
 $\mathcal{M}_r$ is an exponential attractor for $S(t)$ on
 $\mathbb{C}_r$. Since $\mathbb{C}_r$ is also a compact absorbing set,
 the basin of exponential attraction of $\mathcal{M}_r$ is the whole
 phase space $\BBB_r$. This means that $\mathcal{M}_r$ is an
 exponential attractor on $\BBB_r$. Note that the finite fractal
 dimension $N_r$ of $\mathcal M_r$ only depends on $r$  and that  there exists an increasing function $Q_r$ and $\gamma_r$, only depending on $r$ (and also $\zeta,\omega$), such that, for any bounded set $B\subset \BBB_r$,
 \begin{align*}
 	\text{\rm dist}(S(t)B,\mathcal M_r)\leq Q_r(\norm{B}_\Lts)e^{-\gamma_r t},\quad 	\forall t\geq 0.
 \end{align*}
 This comes from the properties \eqref{smoothing} and \eqref{lipt},
 since the constants involved only depend on $T,r,\omega$,  and  $\zeta$. 
 
 Since then $\mathcal A_r\subset \mathcal M_r$ for any $r\geq0$, we deduce that $\AA_r$ is also of finite fractal dimension. This concludes the proof of Theorem \ref{th1}.  
 \subsection{Proof of Theorem   \ref{thm:325}}
 In the previous section we have shown the existence of an exponential
 attractor $\mathcal M_r$  for  the system $(S(t),\BBB_r)$,
 for any $r\geq0$.  Let us  define {$\mathcal
 M:=\bigcup_{m\in \mathbb N}\mathcal M_m$}. Since $\mathcal M_r$ is positively
 invariant for any $r\geq0$,  we have 
 $$
 S(t)\mathcal M\subset \bigcup_{m\in \mathbb N} S(t)\mathcal M_m\subset \mathcal M, \quad \forall t\geq0,
 $$
 i.e., also $\mathcal M\subset \Lts$ is positively invariant. In
 conclusion, the exponential attraction stated in property (3)
 directly comes from the exponential attraction of each {$\mathcal M_m$
 in $\BBB_m$, for any $m\in \mathbb N$}.
 
 \subsection{Proof of Theorem \ref{spe}}
To prove this  statement,  we need first to show that
$\f_K=\Pk((\v\cdot\nabla_\Gamma) \u_K)$ satisfies assumptions
\eqref{C1}--\eqref{C2} and \eqref{uk1} ( as   $\f$  is
assumed to satisfy  these assumptions,  also $\f_{NK}$ will
satisfy them). Assumption \eqref{C1} is trivially  checked,  since 
\begin{align}
\int_\Gamma \Pk((\v\cdot\nabla_\Gamma) \u_K)\cdot \u=\int_\Gamma (\v\cdot\nabla_\Gamma) \u_K\cdot \u_K=0,
\label{ll}
\end{align}
for any $\u\in \L^2(\Gamma)$, since $\v\in \L_\sigma^2(\Gamma)$. Note
that $\u_K=\Pk\u$ belongs also to $\H^1_\sigma(\Gamma)$, since, by
Korn's inequality \eqref{Korn2},
$\norm{\u_K}_{\H^1_\sigma(\Gamma)}\leq C\norm{\u_K}\leq
C\norm{\u}$. Again by Korn's inequality, also assumption \eqref{C2} is
verified, since $\f_K$ is linear in $\u$ and $\v\in
\L^\infty(\Gamma)$.  Eventually,  assumption \eqref{uk1} is
easily verified, again thanks to \eqref{ll}.  Therefore, Theorem
\ref{th1}  concerning   the $\sigma$-attractor $\JJ$ also
applies in this case. Clearly, the same  holds if  $\f_K= \boldsymbol 0$. We now aim at refining the results of the theorem. 

In the case $\f_K=  \boldsymbol 0$, it is immediate to deduce from \eqref{K1}--\eqref{K2} that 
$$
\partial_t \u_K(t)=\boldsymbol 0,
$$   
 so that $\Pk S(t)\u_0=\u_K(t)=\u_K(0)=\Pk\u_0$, for any $\u_0\in \Lts$ and any $t\geq0$, i.e., the Killing component of any trajectory is always constant. In particular, it holds \begin{align}
 	\norm{\Pk S(t)\u_0}=\norm{\Pk\u_0},\quad \forall t\geq0.\label{conserved}
 \end{align}
On the other hand, when $\f_K=\Pk((\v\cdot\nabla_\Gamma) \u_K)$, the $\L^2$-norm of $\u_K$ does not change in time, since we have 
\begin{align*}
	\frac12\frac d {dt} \norm{\u_K}^2=\int_\Gamma (\v\cdot\nabla_\Gamma) \u_K\cdot \u_K=0,
\end{align*}
 since  $\v\in \L^2_\sigma(\Gamma)$. Thus, also in this case \eqref{conserved} holds. We can thus consider the two cases with the same approach.

Let us now introduce, for any fixed $r\geq0$, the (closed) set
\begin{align}
	\label{tB1}
	\widetilde{\BBB}_r:=\{\u\in \Lts:\ \norm{\Pk\u}=r\},
\end{align}
as in \eqref{tB1}, which is a complete metric space if endowed with
the $\Lts$ topology. Note that in this case the family
$\{\widetilde{\BBB}_r\}_{r\geq0}$ is not  monotone.  Then, by \eqref{conserved}, the map $S(t)$ is such that
$$
S(t):\quad \widetilde{\BBB}_r\to \widetilde{\BBB}_r,
$$ 
and thus we can define the dynamical system
$(S(t),\widetilde{\BBB}_r)$. Repeating the same arguments as in the
proof of Theorem \ref{thm1}, by simply substituting $\BBB_r$ with
$\widetilde{\BBB}_r$, we can show that, for any $r\geq0$, there exists
a compact, invariant, attracting set $\widetilde{\AA}_r$ of finite
fractal dimension, which is the global attractor for the dynamical
system $(S(t),\widetilde{\BBB}_r)$. Therefore, it immediately follows
from the properties of $\JJ$ that $\bigcup_{r\geq
  0}\widetilde{\AA}_r\subset \JJ$. For the reverse inclusion, let us
consider $\u\in \JJ$. Then, there exists a complete bounded trajectory
$\xi$ such that $\xi(0)=\u$. Since the trajectory is bounded and by
\eqref{conserved} it holds $\norm{\Pk\xi(t)}=\norm{\Pk\u}$ for any
$t\geq0$,  we have that  $\xi(\R)\subset \widetilde{\BBB}_r$
for $r=\norm{\Pk\u}$. From property \eqref{bdd_traj2} of the
corresponding set $\widetilde{\AA}_r$, which is the global attractor
 of  system $(S(t),\widetilde{\BBB}_r)$, this also means that
$\xi(0)=\u\in \widetilde{\AA}_r$, entailing that $\JJ\subset
\bigcup_{r\geq 0}\widetilde{\AA}_r$.  We  thus conclude that
$\JJ=\bigcup_{r\geq 0}\widetilde{\AA}_r$, i.e., $\JJ$ has  the
desired  pancake-like structure. 

In order now to prove that $\JJ$ is bounded compact and bounded finite
dimensional, we need to consider a closed and bounded set $B\subset
\Lts$ and  prove that  $\JJ\cap B	\subset \AA_r$, for
some $r\geq0$,  where $\AA_r$ is the global attractor for the system
$(S(t),\BBB_r)$ introduced in Theorem \ref{th1}. Indeed, this 
will entail  that also $\JJ\cap B$ is compact and of finite
fractal dimension.  To prove $\JJ\cap B	\subset \AA_r$ for some
$r\geq0$, let us observe that, since $B$ is bounded, there exists
$r\geq 0$ such that $\JJ\cap B\subset \BBB_r$.  We  then consider $\u\in \JJ\cap \BBB_r$. By definition of $\JJ$ there exists a complete bounded trajectory $\xi$ such that $\xi(0)=\u$. By \eqref{conserved} it holds $\norm{\Pk\xi(t)}=\norm{\Pk\u}\leq r$, for any $t\geq0$, so that clearly $\xi(\R) \subset  \BBB_r$. Then, by property \eqref{bdd_traj2} of the corresponding global attractor ${\AA}_r$ of system $(S(t),\BBB_r)$, we deduce $\u\in \AA_r$, entailing $\JJ\cap B\subset \AA_r$.

 As we have proved   that $\JJ$ is bounded closed, it is immediate to prove that $\JJ$ is closed. Indeed, let us consider a sequence $\{\u_n\}_{n\in\N}\subset \JJ$ such that $\u_n\to \u^*$ in $\Lts$ as $n\to \infty$. Then the set $B=\overline{\{\u_n\}_{n\in\N}}^\Lts$ is closed and bounded, so that $\JJ\cap B$ is also closed. Since $\{\u_n\}_{n\in\N}\subset \JJ\cap B$, this means that $\u^*\in \JJ\cap B\subset\JJ$, entailing that $\JJ$ is closed in $\Lts$.

To prove that $\Pk\JJ=\KK$, we refer to the proof of Theorem
\ref{th2}. Indeed, both $\f_K=\boldsymbol 0$ and
$\f_K=\Pk((\v\cdot\nabla_\Gamma) \u_K)$  satisfy the 
assumptions  there,  and thus property 5 in the statement of
Theorem \ref{th2} gives the  identification.  

To conclude the proof we only need to consider the case when,
additionally, $\f_{NK}=\boldsymbol 0$.  In this case, by  repeating the same proof leading to Proposition \ref{pp}, we immediately see that, for any $\u_0\in \Lts$ it holds (see also \cite{Simonett2})
\begin{align*}
	\norm{(\I-\Pk)S(t)\u_0}^2\leq e^{-\zeta t}\norm{(\I-\Pk)\u_0}^2\to 0\quad\text{as }t\to \infty,
\end{align*} 
so that the set $$\BB_s:=\{\u\in \Lts:\ \norm{(\I-\Pk)\u}\leq s\},\quad s> 0,$$
is an absorbing set, for any $s>0$. Therefore, since $\JJ$ is the
minimal closed set attracting bounded sets, it holds  that 
$$
\JJ\subset \bigcap_{s>0}\BB_s,
$$
entailing that $\norm{(\I-\Pk)\u}\leq s$ for any $s>0$ and any $\u\in \JJ$. Therefore, it must be $\norm{(\I-\Pk)\u}=0$ for any $\u\in \JJ$, and thus $\JJ=\KK$, concluding the proof. 
\subsection{Proof of Theorem \ref{interesting}}
The proof of this theorem is analogous to the one of Theorem \ref{th1}. In particular, in this case we can directly consider the map
$$
S(t):\quad \Lts\to \Lts,\quad \forall t\geq 0,
$$
since by assumption there exists $C\subset \Lts$ which is a bounded
absorbing set, so that the system is dissipative. Then, we can define
an absorbing $\Lts$-ball $\BB_r^\Lts$ of radius $r\geq0$ sufficiently
large such that $C\subset\BB_r^\Lts$, and repeat the very same proof
of Theorem \ref{th1}, substituting the absorbing ball $\BB_r^0$ with
$\BB_r^\Lts$ and the phase space $\BBB_r$ with the whole $\Lts$. Note
that in this case we are assuming the additional hypothesis
\eqref{extra} on $\f_{NK}$,  while omitting  \eqref{extra2}, since $\f_K$ does not satisfy \eqref{nega}, but \eqref{pos}. In this way we can retrieve the existence of a global attractor $\AA\subset \Lts$ to the system, which is nonempty, compact, invariant, attracting, and of finite fractal dimension. Since $C$ is an absorbing set, the global attractor is also contained in $C$. Clearly it also holds that $\AA$ coincides with $\JJ$ defined in \eqref{J}, by the natural property \eqref{bdd_traj} of any global attractor (in this case the phase space is directly the whole $\Lts$), i.e., $\AA$ only contains the complete and bounded trajectories in $\Lts$.  Note also that the existence of an exponential attractor $\mathcal M$ is a by-product of the aforementioned proof, from which we deduce the finite dimensionality of $\JJ$. 
Since any trajectory is bounded by assumption, it also holds $\JJ=\mathcal I$, where $\mathcal I$ is defined in \eqref{I}. Note also that, since we are assuming \eqref{extra}, estimate \eqref{gr} holds and thus the set 
$$
\BB:=\left\{\u\in\Lts:\ \norm{(\I-\Pk)\u}\leq \sqrt{\frac 12 +\omega}\right\}
$$
is an absorbing set for the system, entailing $\JJ\subset \mathcal B$.

 To conclude,
 if we assume that $\Pk S(t)B\to 0$ as $t\to \infty$ for any set $B\subset \Lts$ such that $\Pk B$ is bounded, then we can define the absorbing sets 
\begin{align}
	C_s:=\left\{\u\in\Lts:\ \norm{\Pk\u}\leq s\right\},\quad \forall s>0.
\end{align}
Indeed, by assumption, for any $s>0$ and any bounded set $B\subset
\Lts$ there exists $t_{s,B}>0$ such that $S(t)B\subset C_s$ for any
$t\geq t_{s,B}$. Therefore, by the minimal attracting property of
$\JJ$, $\JJ\subset \bigcap_{s>0}C_s$. This entails that, for any
$\u\in\Lts$, $\norm{\Pk\u}=0$, proving  \eqref{best}.

\section{ Attractors for unbounded trajectories:  Proofs of
  Section \ref{unbd}}\label{unbdd}

\subsection{Proof of Theorem \ref{th2}}
To prove the  statement,  we 
 apply  Lemma \ref{A1},  so that we only need to check its
assumptions. 
Let us define the set  
$$
Q_0:=\left\{\v\in \Lts: \norm{(\boldsymbol I-\P_{\mathcal{K}})\v}\leq \frac12+\omega\right\},
$$
where $\omega$ is  given  in \eqref{gr}.  Moroever, let
$$
Q:=\overline{\bigcup_{t\geq 0}S(t)Q_0}.
$$
Recall that in this case we are assuming \eqref{extra} for $\f_{NK}$.
Thanks to estimate \eqref{gr}, it is then immediate to deduce that $Q$
is an absorbing set, i.e., for any bounded set $B\subset \Lts$ there
exists $t_B>0$ such that $S(t)B\subset Q$, for any $t\geq t_B$. We now
need to show that $Q$ is positively invariant. Let us fix $t\geq0$ and
consider $\u\in S(t)Q$. This means that there exists $\u_0\in Q$ such
that $S(t)\u_0=\u$. Therefore, from the definition of $Q$ there exists
a sequence $\{t_n\}_n$, $t_n\geq0$ and a sequence $\{\u_n\}_n\subset
Q_0$ such that $S(t_n)\u_n \to \u_0$ as $n\to \infty$.  Thus, by the
continuity properties of the semigroup  $S(t)$,  we have 
$$
S(t)S(t_n)\u_n=S(t+t_n)\u_n \to S(t)\u_0=\u,\quad\text{as }n\to \infty.
$$
Since $\{S(t+t_n)\u_n\}_{n}\subset \bigcup_{t\geq 0}S(t)Q_0$, this entails that $\u\in Q$, and thus $S(t)Q\subset Q$ for any $t\geq0$, i.e., $Q$ is positively invariant. This means that assumption (H1) of Lemma \ref{A1} is satisfied by this choice of $Q$, and by setting $D_1=D_2=\frac12+\omega$. 

 Concerning  assumption (H2), we recall that, by assumption,
$\f$ is chosen so that we  are  in case (iii) of Proposition \ref{lemmadec} (see also Remark \ref{relax}). Then, for any $\u_0\in\Lts$ such that $\norm{\Pk\u_0}\geq R_0$, we have
\begin{align}\label{exploding}
\norm{\Pk S(t)\u_0}\geq \norm{\Pk\u_0},\quad \forall t\geq 0.
\end{align}
We can thus define the set 
\begin{align*}
	H_R:=\{\u\in Q:\ \norm{\Pk\u_0}\leq R\},\quad \forall R\geq R_0.
\end{align*}
By setting $S(R)=R$ and $R_1=R_0$, we see that $\{\u\in\Lts:\
\norm{\Pk\u}\leq S(R)\}\cap Q\subset H_R$, and,  obviously, 
$S(R)=R\to \infty$ as $R\to \infty$. Moreover, for every $R\geq
R_1=R_0$ it holds $H_R \subset \{\u\in \Lts:\ 	\norm{\Pk\u}\leq
R\}$. In conclusion, thanks to \eqref{exploding} (notice that this
condition justifies the assumptions on  $\f_K$),  we have
$$
S(t)(Q\setminus H_R)\subset Q\setminus H_R,\quad \forall t\geq 0,\quad \forall R\geq R_1=R_0.
$$
Indeed, $Q\setminus H_R=\{\u\in Q:\ \norm{\Pk\u}>R\}$ and, by \eqref{exploding}, $\norm{\Pk S(t)\u}\geq \norm{\Pk \u}>R$, for any $\u\in Q\setminus H_R$. Moreover, since $Q$ is positively invariant, 
$S(t)(Q\setminus H_R)\subset Q$ for any $t\geq 0$. This means that
$S(t)(Q\setminus H_R)\subset Q\setminus H_R$, as desired.  Hence,
 assumption (H2) of Lemma \ref{A1} is satisfied.

 Eventually  assumption (H3)  of Lemma \ref{A1}  is easily verified. Indeed, by Theorem \ref{thm1} (see, in particular, \eqref{regulari}) we know that, for any $\u_0\in \Lts$
$$
\norm{S(\cdot)\u_0}_{L^\infty(\frac t2,2t;\H^1_\sigma(\Gamma))}\leq C(\norm{\u_0},t,T),\quad \forall t>0.
$$
Let us fix $t>0$. Then, for any $B\subset \Lts$ bounded, it holds
$$
\sup_{\u_0\in B}\norm{S(t)\u_0}_{\H^1_\sigma(\Gamma)}\leq C(\norm{B}_{\Lts},t),
$$
where the constant $C>0$ only depends on $t$ and the  $\L^2$-diameter
  of $B$.  The   ball $K(t,B):=\{\u\in
\H^1_\sigma(\Gamma):\ \norm{\u}_{\H^1_\sigma(\gam)}\leq
C(\norm{B}_{\Lts},t) \}$ is compact in $\Lts$ and 
	\begin{align*}
{S(t)B}\subset K(t,B).
\end{align*}
Since both $t>0$ and the bounded set $B$ are arbitrary, this means that the semigroup $S(t)$ is generalized asymptotically compact, verifying assumption (H3) of Lemma \ref{A1}.

Therefore, having verified all the assumptions, Lemma \ref{A1} can be
applied,  which  
proves assertions 1--2  and   4--6 in the statement of
Theorem \ref{th2}.  To prove   property 3  let  us
consider $\u_0\in \JJ$. Since $S(t)\JJ=\JJ$ for any $t\geq0$ there
exists $\u_1\subset \JJ$ such that $\u_0=S(1)\u_1$.  From  the
regularization  property  \eqref{regulari}, it holds  that
 $S(1)\u_1\in \H^1_\sigma(\Gamma)$, entailing that $\u_0\in \H^1_\sigma(\Gamma)$. This means $\JJ\subset \H^1_\sigma(\Gamma)$ and thus $\JJ$ must have empty interior in $\Lts$. 

\subsection{Proof of Theorem \ref{t1b}}
To  check the statement  
we observe that,   by restricting  the semigroup to $(S(t),
(\I-\Pk)\Lts)$, since $\f_K(x,\boldsymbol 0)= \boldsymbol 0$ for any $x\in \Gamma$, we 
obtain  by uniqueness that $\Pk S(t)\u_0\equiv \boldsymbol 0$ for any
$\u_0\in (\I-\Pk)\Lts$.  This implies that $S(t)=(\I-\Pk)S(t)$ on 
$(\I-\Pk)\Lts$, hence  
$$
S(t): \ (\I-\Pk)\Lts\to (\I-\Pk)\Lts,
$$
for any $t\geq 0$. Therefore,  by  exploiting the
 regularization  properties of the system as well as the dissipative
estimate \eqref{gr}, we can argue as in the proof of Theorem \ref{th1}
 and   deduce that there exists a compact absorbing set of the
same form as \eqref{B1r}, but with $r=0$.  By  the standard
theory of dynamical systems,  this entails  that there exists the unique global attractor $\mathcal A_0$ for the dynamical system $(S(t), (\I-\Pk)\Lts)$. This attractor is nonempty, compact, invariant, and attracting. Moreover, $\AA_0$ can be shown (see property \eqref{bdd_traj}) to be composed of complete and bounded trajectories. Therefore we immediately infer $\mathcal A_0\subset \mathcal I$, where $\mathcal  I$ is defined in \eqref{I}, concluding the proof.

pendix

\section{A Lemma on the existence of the unbounded attractor}\label{sec:appendix}
 We conclude with an Appendix presenting a technical lemma, following some ideas
in  
\cite{Carvalho}.  The proof is essentially mutated from
\cite{Carvalho}, up to  minor modifications.  
	\begin{lemm}\label{A1}
		Under the assumptions \eqref{C1}--\eqref{C2}, let us
                consider the dynamical system $(S(t),\Lts)$ defined in
                Section \emph{ \ref{longtt}}. If 
		\begin{enumerate}[label=(\subscript{H}{{\arabic*}})]
		\item There exist $D_1, D_2 > 0$ and a closed set
                  $Q\subset \Lts$ with $$\{\u\in \Lts:\	\norm{(\I -
                    \Pk)\u}	 \leq D_1\} \subset  Q \subset \{
                  \u\in \Lts:\ 	\norm{(\I -
			\Pk)\u}	 \leq D_2\}$$ such that $Q$ is an absorbing set, i.e., for any bounded set $B\subset \Lts$ there exists $t_B>0$ such that 
			$$
			S(t)B\subset Q,\quad \forall t\geq t_B,
			$$
			and positively invariant, i.e., for every
		$B\subset \Lts$ there exists $t_B > 0$  such that
		$S(t)Q \subset  Q$ for every
		$t \geq  t_B$;
		\item There exist two constants $R_0$ and $R_1$ with $0 < R_0 \leq R_1$ and an increasing family of
		closed and bounded sets $\{H_R\}_{R\geq R_0}$ with $H_R \subset Q$, such that
		\begin{itemize}
		\item for every $R \geq  R_1$ we can find $S(R) \geq R_0$ such that $\{\u\in\Lts:\ \norm{\Pk \u}	 \leq S(R)\} \cap Q \subset  H_R$ and
		moreover $\lim_{R\to\infty} S(R)=+\infty$,
		\item for every $R \geq R_1$ we have $H_R \subset \{\u\in \Lts:\ 	\norm{\Pk\u}\leq R\}$, 
		\item  for every $R \geq  R_1$, $S(t)(Q \setminus H_R) \subset Q \setminus H_R $ for every $t\geq0$;
		\end{itemize}
		\item The semigroup $\{S(t)\}_{t\geq0}$ is generalized asymptotically compact, i.e., if for every
	bounded set	$B\subset \Lts$ and every $t > 0$ there exists a compact set $K(t, B) \subset \Lts$ and $\varepsilon(t, B) \to 0$ as $t\to\infty$ such that
	\begin{align}
		S(t)B\subset O_\varepsilon(K(t,B)),
		\label{ascmp}
	\end{align}
i.e., for any $\u\in S(t)B$ it holds $\text{\rm dist}_\Lts(\u,K(t,B))\leq \varepsilon$,
		\end{enumerate}
		then the nonempty set $\mathcal J$ defined in \eqref{J} is the unbounded attractor, satisfying:
		\begin{enumerate}
			\item $\mathcal J=\overline{\bigcap_{t\geq0}S(t)Q}$,
			\item $\mathcal J$ is closed and invariant, i.e., $S(t)\mathcal J =\mathcal J$ for all $t\geq0$,
			\item $\mathcal J$ is bounded compact, i.e., the intersection of $\mathcal J$ with any closed and bounded set $B\subset \Lts$ is compact,
			\item If for some $R\geq 0$, $B \in\Lts$ bounded, and $t_1 > 0$ the sets $S(t)B \cap \{\v\in \Lts:\	\norm{\P_\mathcal K \v}\leq R\}$ are nonempty for every $t \geq  t_1$, then
			$$
			\lim_{t\to\infty}
			{\rm dist}(S(t)B \cap  \{\v\in \Lts:\	\norm{\P_\mathcal K \v}\leq R\}, \mathcal J ) = 0,$$
			and $\mathcal J$ is the minimal closed set with the above property,
			\item $\P_\mathcal K \mathcal J = \mathcal K$.
		\end{enumerate}
	\end{lemm}
	\begin{proof}
	Concerning the proof, we refer to \cite[Theorem
        3]{Carvalho}. In the notations of  such  theorem, we have $X=\Lts$, $E^+=\KK$, {which is a finite dimensional subspace of $X$}, and $E^-=(\I-\Pk)\Lts$, {which is clearly closed}. To be precise, in the aforementioned proof there is a stronger assumption on the semigroup $S(t)$, namely
		$$
		S\in C([0,+\infty)\times \Lts;\Lts), 
		$$ 
		i.e., it is jointly time-space continuous. Actually, this {additional }assumption can be relaxed to 
			\begin{align}
		S\in C([0,+\infty)\times \KK;\Lts), 
		\label{KK}\end{align}
		i.e., that $S$ is jointly continuous if we restrict the phase space to $\KK$. Indeed, the only part in the proof of \cite[Theorem 3]{Carvalho} which exploits this regularity is to show that $\JJ$ is nonempty and $\Pk\JJ=\KK$, i.e., \cite[Lemma 3]{Carvalho}. In this case, all the proof is based on the fact that, having defined
			$$
			B_R:=\{\u\in \KK:\ \norm{\u}<R+1\},
			$$
			for any $R\geq R_0$, we need
			$$
			S\in C([0,T]\times \overline{B}_R;\Lts),\quad \forall T>0,\quad \forall R\geq R_0,
			$$
			and this is ensured by  the weaker 
                        \eqref{KK},  as well.  
			
			In our case, the dynamical system $(S(t),\Lts)$ under consideration satisfies this assumption and thus \cite[Theorem 3]{Carvalho} can be applied. Indeed, note that, by Korn's inequality \eqref{Korn2}, 
				$$
				\norm{\u}_{\H^1(\Gamma)}\leq C\norm{\u},\quad \forall \u\in\KK,
				$$
				and thus, by Theorem \ref{thm1} we have that the semigroup $S(t)$, when restricted to $\KK$, is such that $\partial_t S(t)\u_0\in L^2(0,T;\Lts)$ for any $T\geq0$ and any $\u_0\in \KK$ (see \eqref{strong}), where the bounding constants only depend on $\norm{\u_0}_{\H^1(\Gamma)}$, the parameters of the problem, $\f$, and $\Gamma$. This shows, as in \eqref{lipt}, that the map $S(\cdot)\u_0: [0,T]\to \Lts$ is $(1/2)$-H\"{o}lder continuous for any $T>0$ and any fixed $\u_0\in \KK$. Together with the continuous dependence estimate \eqref{contdep2} (in which the constants appearing depend on $T>0$, $\norm{\u_{0,i}}$, $i=1,2$, the parameters of the problem, $\f$, and $\Gamma$), this allows to show that 
				$$
				S\in C^{\frac12}([0,T]\times  D;\Lts),
				$$
				for any $T>0$ and any $D\subset \KK$ closed bounded set, i.e., $S$ is locally $(1/2)$-H\"{o}lder continuous in $[0,+\infty)\times \KK$, which entails \eqref{KK} and concludes the proof.
                              \end{proof}

\section*{Acknowledgments} 
This research was funded in
whole or in part by the Austrian Science Fund (FWF) projects  10.55776/ESP552, 10.55776/F65,  10.55776/I5149,
10.55776/P32788, 
as well as by the OeAD-WTZ project CZ 09/2023. AP is  a member of Gruppo Nazionale per l’Analisi Matematica, la Probabilità e le loro Applicazioni (GNAMPA) of
Istituto Nazionale per l’Alta Matematica (INdAM), and gratefully acknowledges support  from the Alexander von Humboldt Foundation. For
open-access purposes, the authors have applied a CC BY public copyright
license to any author-accepted manuscript version arising from this
submission. Part of this research was conducted during a visit to the Mathematical Institute at Tohoku University, whose warm hospitality is gratefully acknowledged.

\bibliography{Bibliography}

\def\ocirc#1{\ifmmode\setbox0=\hbox{$#1$}\dimen0=\ht0 \advance\dimen0 by1pt\rlap{\hbox to\wd0{\hss\raise\dimen0 \hbox{\hskip.2em$\scriptscriptstyle\circ$}\hss}}#1\else {\accent"17 #1}\fi}
\begin{thebibliography}{10}

\bibitem{AGG}
H.~Abels, H.~Garcke, and G.~Gr{\"{u}}n.
\newblock Thermodynamically consistent, frame indifferent diffuse interface models for incompressible two-phase flows with different densities.
\newblock {\em Math. Models Methods Appl. Sci.}, 22(3):1150013 (40 pages), 2012.

\bibitem{AGP2024}
H.~Abels, H.~Garcke, and A.~Poiatti.
\newblock Diffuse interface model for two-phase flows on evolving surfaces with different densities: global well-posedness.
\newblock {\em Submitted}, 2024.

\bibitem{AGP2024b}
H.~Abels, H.~Garcke, and A.~Poiatti.
\newblock Diffuse interface model for two-phase flows on evolving surfaces with different densities: local well-posedness.
\newblock {\em Submitted}, 2024.

\bibitem{AGP}
H.~Abels, H.~Garcke, and A.~Poiatti.
\newblock Mathematical analysis of a diffuse interface model for multi-phase flows of incompressible viscous fluids with different densities.
\newblock {\em J. Math. Fluid Mech.}, 26(2):Paper No. 29, 51, 2024.

\bibitem{AWe}
H.~Abels and J.~Weber.
\newblock Local well-posedness of a quasi-incompressible two-phase flow.
\newblock {\em J. Evol. Equ.}, 21(3):3477--3502, 2021.

\bibitem{Bernardi}
R.~Agroum, S.~Mani~Aouadi, C.~Bernardi, and J.~Satouri.
\newblock Spectral discretization of the {N}avier-{S}tokes equations coupled with the heat equation.
\newblock {\em ESAIM Math. Model. Numer. Anal.}, 49(3):621--639, 2015.

\bibitem{Arnaudon}
M.~Arnaudon and A.~B. Cruzeiro.
\newblock Lagrangian {N}avier-{S}tokes diffusions on manifolds: variational principle and stability.
\newblock {\em Bull. Sci. Math.}, 136(8):857--881, 2012.

\bibitem{Carvalho}
J.~Banaśkiewicz, A.~Carvalho, J.~Garcia-Fuentes, and P.~Kalita.
\newblock Autonomous and non-autonomous unbounded attractors in evolutionary problems.
\newblock {\em Journal of Dynamics and Differential Equations}, pages 1--54, 12 2022.

\bibitem{Tartar}
C.~Bardos and L.~Tartar.
\newblock Sur l'unicit\'e{} r\'etrograde des \'equations d'\'evolution.
\newblock {\em C. R. Acad. Sci. Paris S\'er. A-B}, 273:A1239--A1241, 1971.

\bibitem{Foias}
A.~Biswas, C.~Foias, and A.~Larios.
\newblock On the attractor for the semi-dissipative {B}oussinesq equations.
\newblock {\em Ann. Inst. H. Poincar\'e{} C Anal. Non Lin\'eaire}, 34(2):381--405, 2017.

\bibitem{Bortolan}
M.~C. Bortolan and J.~da~Silva.
\newblock Sufficient conditions for the existence and uniqueness of maximal attractors for autonomous and nonautonomous dynamical systems.
\newblock {\em Journal of Dynamics and Differential Equations}, pages 1--30, 11 2022.

\bibitem{Brandner}
P.~Brandner, A.~Reusken, and P.~Schwering.
\newblock On derivations of evolving surface {N}avier-{S}tokes equations.
\newblock {\em Interfaces Free Bound.}, 24(4):533--563, 2022.

\bibitem{Robinsonbook}
A.~N. Carvalho, J.~A. Langa, and J.~C. Robinson.
\newblock {\em Attractors for infinite-dimensional non-autonomous dynamical systems}, volume 182 of {\em Applied Mathematical Sciences}.
\newblock Springer, New York, 2013.

\bibitem{Chan}
C.~H. Chan, M.~Czubak, and M.~M. Disconzi.
\newblock The formulation of the {N}avier-{S}tokes equations on {R}iemannian manifolds.
\newblock {\em J. Geom. Phys.}, 121:335--346, 2017.

\bibitem{Chepyzov}
V.~Chepyzhov and A.~Goritskii.
\newblock {U}nbounded attractors of evolution equations.
\newblock {\em Advances in Soviet Mathematics}, 10:85--128, 1992.

\bibitem{WeakGlobalNSAttractor}
A.~Cheskidov and C.~Foias.
\newblock On global attractors of the 3{D} {N}avier-{S}tokes equations.
\newblock {\em J. Differential Equations}, 231(2):714--754, 2006.

\bibitem{Constantin88}
P.~Constantin and C.~Foias.
\newblock {\em Navier-{S}tokes equations}.
\newblock Chicago Lectures in Mathematics. University of Chicago Press, Chicago, IL, 1988.

\bibitem{Escauriaza}
L.~Escauriaza, G.~Seregin, and V.~Sver\'{a}k.
\newblock Backward uniqueness for parabolic equations.
\newblock {\em Arch. Ration. Mech. Anal.}, 169(2):147--157, 2003.

\bibitem{Foias2}
C.~Foias and R.~Temam.
\newblock The connection between the {N}avier-{S}tokes equations, dynamical systems, and turbulence theory.
\newblock In M.~G. Crandall, P.~H. Rabinowitz, and R.~E. Turner, editors, {\em Directions in Partial Differential Equations}, pages 55--73. Academic Press, 1987.

\bibitem{Galdi94}
G.~P. Galdi.
\newblock {\em An introduction to the mathematical theory of the {N}avier-{S}tokes equations. {V}ol. {II}}, volume~39 of {\em Springer Tracts in Natural Philosophy}.
\newblock Springer-Verlag, New York, 1994.
\newblock Nonlinear steady problems.

\bibitem{Girault86}
V.~Girault and P.-A. Raviart.
\newblock {\em Finite element methods for {N}avier-{S}tokes equations}, volume~5 of {\em Springer Series in Computational Mathematics}.
\newblock Springer-Verlag, Berlin, 1986.
\newblock Theory and algorithms.

\bibitem{GPPV}
M.~Grasselli, N.~Parolini, A.~Poiatti, and M.~Verani.
\newblock Non-isothermal non-{N}ewtonian fluids: the stationary case.
\newblock {\em Math. Models Methods Appl. Sci.}, 33(9):1747--1801, 2023.

\bibitem{Iskauriaza}
L.~Iskauriaza, G.~A. Ser\"{e}gin, and V.~Shverak.
\newblock {$L_{3,\infty}$}-solutions of {N}avier-{S}tokes equations and backward uniqueness.
\newblock {\em Uspekhi Mat. Nauk}, 58(2(350)):3--44, 2003.

\bibitem{V1}
T.~Jankuhn, M.~A. Olshanskii, and A.~Reusken.
\newblock Incompressible fluid problems on embedded surfaces: modeling and variational formulations.
\newblock {\em Interfaces Free Bound.}, 20(3):353--377, 2018.

\bibitem{Koba}
H.~Koba, C.~Liu, and Y.~Giga.
\newblock Energetic variational approaches for incompressible fluid systems on an evolving surface.
\newblock {\em Quart. Appl. Math.}, 75(2):359--389, 2017.

\bibitem{Koba2}
H.~Koba, C.~Liu, and Y.~Giga.
\newblock Errata to ``{E}nergetic variational approaches for incompressible fluid systems on an evolving surface'' [ {MR}3614501].
\newblock {\em Quart. Appl. Math.}, 76(1):147--152, 2018.

\bibitem{Zelik}
A.~Miranville and S.~Zelik.
\newblock Attractors for dissipative partial differential equations in bounded and unbounded domains.
\newblock In {\em Handbook of differential equations: evolutionary equations. {V}ol. {IV}}, Handb. Differ. Equ., pages 103--200. Elsevier/North-Holland, Amsterdam, 2008.

\bibitem{Myers}
S.~B. Myers.
\newblock Isometries of 2-dimensional riemannian manifolds into themselves.
\newblock {\em Proc. Natl. Acad. Sci. USA}, 22(5):297--300, 1936.

\bibitem{Olshsup}
M.~A. Olshanskii.
\newblock {On equilibrium states of fluid membranes}.
\newblock {\em Physics of Fluids}, 35(6):062111, 06 2023.

\bibitem{Priebe}
V.~Priebe.
\newblock Solvability of the {N}avier-{S}tokes equations on manifolds with boundary.
\newblock {\em Manuscripta Math.}, 83(2):145--159, 1994.

\bibitem{greenbook}
J.~Pr\"uss and G.~Simonett.
\newblock {\em Moving interfaces and quasilinear parabolic evolution equations}, volume 105 of {\em Monographs in Mathematics}.
\newblock Birkh\"auser/Springer, [Cham], 2016.

\bibitem{Simonett}
J.~Pr\"{u}ss, G.~Simonett, and M.~Wilke.
\newblock On the {N}avier-{S}tokes equations on surfaces.
\newblock {\em J. Evol. Equ.}, 21(3):3153--3179, 2021.

\bibitem{Saal1}
J.~Saal.
\newblock Strong solutions for the {N}avier-{S}tokes equations on bounded and unbounded domains with a moving boundary.
\newblock In {\em Proceedings of the {S}ixth {M}ississippi {S}tate--{UBA} {C}onference on {D}ifferential {E}quations and {C}omputational {S}imulations}, volume~15 of {\em Electron. J. Differ. Equ. Conf.}, pages 365--375. Southwest Texas State Univ., San Marcos, TX, 2007.

\bibitem{Sakai}
T.~Sakai.
\newblock {\em Riemannian geometry}, volume 149 of {\em Translations of Mathematical Monographs}.
\newblock American Mathematical Society, Providence, RI, 1996.
\newblock Translated from the 1992 Japanese original by the author.

\bibitem{Seregin}
G.~Seregin and V.~Sver\'{a}k.
\newblock The {N}avier-{S}tokes equations and backward uniqueness.
\newblock In {\em Nonlinear problems in mathematical physics and related topics, {II}}, volume~2 of {\em Int. Math. Ser. (N. Y.)}, pages 353--366. Kluwer/Plenum, New York, 2002.

\bibitem{Simonett2}
G.~Simonett and M.~Wilke.
\newblock {$H^\infty $}-calculus for the surface {S}tokes operator and applications.
\newblock {\em J. Math. Fluid Mech.}, 24(4):Paper No. 109, 23, 2022.

\bibitem{Temam84}
R.~Temam.
\newblock {\em Navier-{S}tokes equations}, volume~2 of {\em Studies in Mathematics and its Applications}.
\newblock North-Holland Publishing Co., Amsterdam, third edition, 1984.
\newblock Theory and numerical analysis, With an appendix by F. Thomasset.

\bibitem{Temam}
R.~Temam.
\newblock {\em Infinite-dimensional dynamical systems in mechanics and physics}.
\newblock Springer, Berlin-Heidelberg-New York, 1988.

\end{thebibliography}
\bibliographystyle{abbrv}

\end{document}